\documentclass[preprint,12pt,a4paper]{elsarticle}

\usepackage{setspace}

\journal{Comput. Methods Appl. Mech. Eng.}

\usepackage{amsmath,amsfonts,amssymb,amscd,amsthm}
\usepackage{stmaryrd}
\usepackage{dsfont}
\usepackage{epsfig}
\usepackage{comment}
\usepackage{xcolor}
\usepackage{color}
\usepackage{soul}
\usepackage{psfrag}
\usepackage{colordvi}
\usepackage{cancel} 
\usepackage{multirow}
\usepackage{soul}
\usepackage[T1]{fontenc}
\usepackage[utf8]{inputenc}
\usepackage{graphicx}
\DeclareGraphicsExtensions{.png,.pdf}

\includecomment{notaNPM}
\excludecomment{notaNPM}

\usepackage[normalem]{ulem}
\newcommand{\Rem}[1]{\textcolor{black}{#1}}  			

\newcommand{\aHDG}{a_{\rm hdg}}
\newcommand{\at}[2]{\left.{#1}\right|_{#2}}

\newcommand{\bC}{\bm{C}}
\newcommand{\bm}[1]{\text{\boldmath $#1$\unboldmath}}
\newcommand{\bn}{\bm{n}}

\newcommand{\btau}{\bm{e}_{\bm{\varphi}}^\mp}
\newcommand{\btauZ}{\chi_{\bm{\varphi}}^\mp}
\newcommand{\btauC}{\psi_{\bm{\varphi}}^\mp}
\newcommand{\bv}{\bm{v}}

\newcommand{\bvarphi}{\bm{\varphi}}

\newcommand{\CTb}{C_2}
\newcommand{\curl}{\bm{\nabla}\!\!\times\!}
\newcommand{\dO}{\ \text{d} \Omega}
\newcommand{\dG}{ {\rm d}\Gamma}
\newcommand{\edge}{e}
\newcommand{\elem}{K}
\newcommand{\ellHDG}{\ell_{\rm hdg}}
\newcommand{\eu}{e_h}
\newcommand{\eq}{\bm{e}_h}
\newcommand{\euA}{\varepsilon_h}
\newcommand{\eqA}{\bm{\varepsilon}_h}
\newcommand{\grad}{\bm{\nabla}}
\newcommand{\GD}{\Gamma_{\rm D}}
\newcommand{\GN}{\Gamma_{\rm N}}
\newcommand{\jump}[1]{\ensuremath{[\![#1]\!]} }

\newcommand{\Lag}{L^\mp}
\newcommand{\mean}[1]{\ensuremath{\{\!\{#1\}\!\}}}
\newcommand{\norm}[1]{\vertiiis{#1}}

\newcommand{\numel}{n_{\rm el}}

\newcommand{\proj}{\Pi}
\newcommand{\projk}{\proj_K^{\hat{q}}}
\newcommand{\projg}{\proj_\edge^{\bar{q}}}
\newcommand{\q}{\bm{q}}
\newcommand{\qA}{\bm{\zeta}}
\newcommand{\qh}{\bm{\tilde{q}}_h}
\newcommand{\qhA}{\bm{\tilde{\zeta}}_h}
\newcommand{\qhAP}{\bm{\tilde{\zeta}}_h^\pi}
\newcommand{\qhAZ}{\bm{\tilde{\zeta}}_h^0}
\newcommand{\qhHDG}{\bm{\q}_h}
\newcommand{\qhHDGA}{\bm{\zeta}_h}
\newcommand{\qhHDGT}{\widehat{\bm{\q}}_h}
\newcommand{\qhHDGTA}{\widehat{\bm{\zeta}}_h}
\newcommand{\qhP}{\bm{\tilde{q}}_h^\pi}
\newcommand{\qhPA}{\bm{\tilde{\varphi}}_h}
\newcommand{\qPA}{\bm{\varphi}}
\newcommand{\qhZ}{\bm{\tilde{q}}_h^0}

\newcommand{\sh}{\tilde{s}_h}
\newcommand{\slb}{s_h^-}
\newcommand{\sub}{s_h^+}
\newcommand{\source}{f}
\newcommand{\test}{\mathcal{W}}
\newcommand{\testeh}{M_h^p}
\newcommand{\testh}{\mathcal{W}_h^p}
\newcommand{\testq}{\bm{\mathcal{V}}}
\newcommand{\testqh}{\bm{\mathcal{V}}_h^p}
\newcommand{\trac}{g_{_{\rm N}}}
\newcommand{\triang}{\mathcal{T}_h}
\newcommand{\triangE}{\mathcal{E}_h}
\newcommand{\uD}{g_{_{\rm D}}}
\newcommand{\uDe}{\hat{g}_{_{\rm D}}}
\newcommand{\uA}{\xi}
\newcommand{\uh}{\tilde{u}_h}
\newcommand{\uhA}{\tilde{\xi}_h}
\newcommand{\uhHDG}{u_h}
\newcommand{\uhHDGA}{\xi_h}
\newcommand{\uhHDGT}{\widehat{u}_h}
\newcommand{\uhPA}{\tilde{\phi}_h}
\newcommand{\uPA}{\phi}
\newcommand{\vertiii}[1]{{\left\vert\kern-0.25ex\left\vert\kern-0.25ex\left\vert #1 \right\vert\kern-0.25ex\right\vert\kern-0.25ex\right\vert}}
\newcommand{\vertiiis}[1]{{\vert\kern-0.25ex\vert\kern-0.25ex\vert #1 \vert\kern-0.25ex\vert\kern-0.25ex\vert}}

\newtheorem{Remark}{Remark}
\newtheorem{Theorem}{Theorem}

\makeatletter
\def\@author#1{\g@addto@macro\elsauthors{\normalsize%
    \def\baselinestretch{1}%
    \upshape\authorsep#1\unskip\textsuperscript{%
      \ifx\@fnmark\@empty\else\unskip\sep\@fnmark\let\sep=,\fi
      \ifx\@corref\@empty\else\unskip\sep\@corref\let\sep=,\fi
      }%
    \def\authorsep{\unskip,\space}%
    \global\let\@fnmark\@empty
    \global\let\@corref\@empty  
    \global\let\sep\@empty}%
    \@eadauthor={#1}
}
\makeatother

\begin{document}

\begin{frontmatter}



\title{\Rem{A posteriori goal-oriented bounds for the Poisson problem using potential and equilibrated flux reconstructions: application to the hybridizable discontinuous Galerkin method}
}


\author{N. Par\'es\corref{cor1}\fnref{label1}}
\author{N. C. Nguyen\fnref{label3}}
\author{P. D\'iez\fnref{label1,label2}}
\author{J. Peraire\fnref{label3}}

\address[label1]{Laboratori de C\`alcul Num\`eric (LaC\`aN),
Universitat Polit\`ecnica de Catalunya, Jordi Girona 1-3, 
E 08034 Barcelona, Spain}
\address[label2]{International Centre for Numerical Methods in Engineering, CIMNE, Barcelona, Spain}
\address[label3]{Center for Computational Engineering, Department of Aeronautics and Astronautics, Massachusetts Institute of Technology, 77 Massachusetts Avenue, Cambridge, MA 02139, USA}

\cortext[cor1]{nuria.pares@upc.edu, URL: http://www-lacan.upc.edu, Tel: +34 934137314, Fax: +34 934011825
}



\begin{abstract}
We present a general framework to compute upper and lower bounds for linear-functional outputs of the exact solutions of the Poisson equation 
\Rem{
based on reconstructions of the field variable and flux for both the primal and adjoint problems. The method is devised from a generalization of the complementary energy principle and the duality theory.}
Using duality \Rem{theory}, the computation of bounds is reduced to \Rem{finding} independent potential and equilibrated flux reconstructions. A generalization of this result is also introduced allowing to derive alternative guaranteed bounds from nearly-arbitrary $\mathcal{H}(\rm{div};\Omega)$ flux reconstructions (only zero-order equilibration is required). \Rem{
This approach is applicable to} any numerical method used to compute the solution. In this work, the proposed approach is applied to derive bounds for the hybridizable discontinuous Galerkin (HDG) method.
An attractive feature of the proposed approach is that 
superconvergence on the bound gap is achieved, yielding accurate bounds even for very coarse meshes. Numerical experiments are presented to illustrate the performance and convergence of the bounds for the HDG method in both uniform and adaptive mesh refinements.
\end{abstract}

\begin{keyword}
exact/guaranteed/strict bounds for quantities of interest, output bounds, goal-oriented error estimation, adaptivity, potential and equilibrated flux reconstructions, hybridizable discontinuous Galerkin method (HDG).
\end{keyword}

\end{frontmatter}



\section{Introduction}

In many applications in computational science and engineering, the numerical approximations are used to accurately assess some target quantities or \emph{quantities of interest}. That is, to provide information on specific features of the true solution $u$, usually given by a linear functional $s=\ell^O(u)$. The approximations are computed using the numerical solution $u_h$, namely $s_h=\ell^O(u_h)$. In this context, it is crucial to assess the quality of the approximated outputs.

Numerous advances in goal-oriented error estimation have been done in recent years. The most well-established techniques provide approximations or bounds for the error in the computed numerical approximation $\ell^O(u)-\ell^O(u_h)$ and produce error indicators to drive goal-oriented mesh adaptivity, see for instance
\cite{MR1487947,MR1665351,MR1885308,MR2107001,MR2857237,MR3085344,MR3327021,MR4008052,MR3639952,MR3683070,MR3725770}. However, in practical applications, two other parallel lines of research are worth mentioning. 
The first one consists of techniques aimed at obtaining more accurate
approximations of the quantities of interest \cite{MR2009374,MR3725770,MR3028504,MR3741980}. In this case, the numerical approximation $u_h$ is used to either compute a new more accurate approximation $\tilde{u}_h$ yielding a more accurate approximation for the quantity of interest $\tilde{s}_h=\ell^O(\tilde{u}_h)$ or to directly compute a better approximation for the quantity of interest $\tilde{s}_h=\tilde{\ell}^O(u_h)$. 
The second line of research aims at the computation of certificates and guaranteed bounds for the quantity of interest, see for instance \cite{MR2114293,MR2116366,MR2186143,MR2186144,MR2399861,MR2899560,MR3633659,LADEVEZE2008,MR2981885}.
 Indeed, besides having an accurate approximation of the quantity of interest (either $\ell^O(u_h), \ell^O(\tilde{u}_h)$ or $\tilde{\ell}^O(u_h)$) in decision-making processes, it is important to be able to provide a guaranteed interval where the exact quantity of interest lies, that is, to guarantee that $s\in [\slb,\sub]$ where $\slb$ and $\sub$ should be fully computable, constant-free guaranteed upper and lower bounds. 
In this context, it is no longer important to directly assess the error in the original approximation of the quantity of interest $s-s_h$, but being able to compute a new improved approximation $\tilde{s}_h$ and providing a guaranteed bounding interval for the exact output, $[\slb,\sub]$, containing both $s$ and $\tilde{s}_h$.
It is also desirable that the new approximation and the bound gap $\sub-\slb$ converge faster than the original approximation.

The present work aims at addressing 
the computation of highly accurate approximations for the quantity of interest and 
providing certificates for the exact value of the quantity of interest. In particular, although a general framework for computing guaranteed bounds for quantities of interest is provided, 
accurate approximations for the quantity of interest and associated guaranteed bounds are obtained from hybridizable discontinuous Galerkin (HDG) approximations of the Poisson equation, where the superconvergence properties of the approximation are exploited to obtain optimally \Rem{convergent} approximations and bounds for the quantity of interest. Also, goal-oriented error indicators are provided to enhance the convergence of adaptive remeshing for non-smooth problems.

HDG methods have gained popularity in the last decade due to their reduced computational cost with respect to classical discontinuous Galerkin methods while retaining superconvergence properties \cite{HDGLab2020}. Also, a very attractive feature is that a simple post-process of the solution yields equilibrated $\mathcal{H}(\rm{div};\Omega)$ approximations of the fluxes. These fluxes are used to compute guaranteed bounds either for the energy norm or for quantities of interest \cite{phdthesis_Wong,MR3850360}. In the present work, the superconvergence properties of the high-order HDG method presented in \cite{MR2513831} are exploited to achieve optimal convergence when approximating and certifying quantities of interest.

The paper is organized as follows: In Section \ref{sec:model_problem}, we introduce the model problem and notations for the quantities of interest and adjoint problem. In Section \ref{sec:bounds}, a general framework to compute guaranteed bounds for quantities of interest by means of potential and equilibrated flux reconstructions is presented. In particular, Section \ref{subsec:bounds_equilibrated_projected}
presents an extension that allows both to compute bounds when non-polynomial data is present and to compute bounds using simplified zero-order equilibrated reconstructions. Section 
\ref{subsec:bounds_equilibrated_ZO}
particularizes the expression for the bounds to high-order projections of the flux reconstructions. Finally Section 
\ref{subsec:bounds_equilibrated_enhanc}
presents an exact representation for the quantity of interest allowing to enhance the bounds using lower bounds for the energy norm. In Section 
\ref{sec:HDG}, we particularize the results derived in Section 
\ref{sec:bounds}
to the HDG method, providing both an accurate alternative approximation for the quantity of interest and its associated guaranteed bounds. Section 
\ref{sec:numerical_examples}
shows the behavior of the proposed technique in two numerical examples, and we present some concluding remarks in Section 
\ref{sec:rem}. The proofs of the most significant results are presented in Appendices A, B, C and D.


%

\section{Model problem}
\label{sec:model_problem}

Consider the Poisson's equation in a polygonal/polyhedral domain $\Omega\subset\mathbb{R}^d$ for $d=2$ or $3$,
\begin{equation}\label{eq:strong}
\begin{array}{rcll}
    -\nabla\cdot(\nu \grad u)       &\!\!=\!\!& \source     & \text{ in } \Omega, \\
    u                           &\!\!=\!\!& \uD         & \text{ on } \GD, \\
    -\nu \grad  u \cdot \bn  &\!\!=\!\!& \trac       & \text{ on }
    \GN,
\end{array}
\end{equation}
where the boundary $\partial\Omega$ is divided into two disjoint parts $\GD$ and $\GN$ such that $\partial\Omega=\bar{\Gamma}_{\rm D} \cup \bar{\Gamma}_{\rm N} $, $\GD \cap \GN = \emptyset$ and $\GD$ is a non-empty set. The data are assumed to be sufficiently smooth, that is, $\source\in\mathcal{L}^2(\Omega)$, $\trac\in\mathcal{L}^2(\GN)$, $\uD\in\mathcal{C}(\GD)$ and $\nu\in\mathcal{L}^\infty(\Omega)$ is assumed to be strictly positive. 
Moreover, for simplicity,  $\nu$ is assumed to be piecewise constant on subdomains of $\Omega$.
%
 
The equivalent mixed formulation of \eqref{eq:strong} is 
\begin{equation}\label{eq:strongmixed}
\begin{array}{rcll}
\q &\!\!=\!\!&-\nu \grad  u  & \text{ in } \Omega,\\
\nabla \cdot \q &\!\!=\!\!& \source & \text{ in } \Omega,\\
u &\!\!=\!\!& \uD & \text{ on } \GD,\\
\q\cdot\bn &\!\!=\!\!& \trac & \text{ on } \GN.
\end{array}
\end{equation}

To introduce the weak form of \eqref{eq:strongmixed}, consider the test spaces $\test = \mathcal{H}^1(\Omega)$ and $\testq = \mathcal{H}(\rm{div};\Omega)=\{\bv\in\mathcal{L}^2(\Omega), \nabla\cdot\bv\in\mathcal{L}^2(\Omega)\}$, and the integral inner products 
\[
(\q,\bv)_\omega = \int_\omega \q\cdot\bv\dO\quad,\quad
(u,v)_\omega=\int_\omega uv\dO\quad\text{and}\quad
\langle u,v\rangle _{\gamma} = \int_\gamma uv\dG,
\]
$\omega$ being a domain in $\mathbb{R}^d$ and $\gamma$ being a domain in $\mathbb{R}^{d-1}$. The subscript $\omega$ is omitted when $\omega$ is the full domain $\Omega$. Recall that for any $\omega\subset\Omega$,  $\q\in\mathcal{H}(\rm{div};\Omega)$ and $w\in\test$ the following Green formula holds
\begin{equation}\label{eq:green}
(\q,\grad  w)_\omega + (\nabla\cdot\q, w)_\omega = \langle \q\cdot\bn ,w\rangle _{\partial\omega}.
\end{equation}
Then, the weak solution of \eqref{eq:strongmixed} is $(u,\q)\in\test\times\testq$ such that
\[
\begin{array}{ll}
\displaystyle
(\nu^{-1}\q,\bv)-(u,\nabla\cdot\bv) + \langle u,\bv\cdot\bn\rangle _{\GN} = - \langle \uD,\bv\cdot\bn\rangle _{\GD}& \forall \bv \in \testq,\\
-(\q,\grad w)  + \langle  \q\cdot\bn, w\rangle _{\GD} =  (\source,w) - \langle \trac,w\rangle _{\GN}& \forall w \in \test,
\end{array}
\]
or equivalently
\begin{equation}\label{eq:general_variational_formulation}
a(u,\q;w,\bv) = \ell(w,\bv) \quad \forall (w,\bv)\in \test\times\testq,
\end{equation}
for
\[
\begin{array}{l}
\displaystyle
a(u,\q;w,\bv) =(\nu^{-1}\q,\bv)-(u,\nabla\cdot\bv) + \langle u,\bv\cdot\bn\rangle _{\GN}-(\q,\grad w)  + \langle  \q\cdot\bn, w\rangle _{\GD},
\\
\displaystyle
\ell(w,\bv)=(\source,w) - \langle \uD,\bv\cdot\bn\rangle _{\GD}- \langle \trac,w\rangle _{\GN} .
\end{array}
\]

\begin{Remark}
\bigskip For any $(u,\q)\in \test\times \testq$ and $(w,\bv)\in \test\times \testq$ it holds that
\[
a(u,\q;w,\bv) =(\nu^{-1}\q,\bv) +(\bv,\grad u)-(\q,\grad w)  + \langle  \q\cdot\bn, w\rangle _{\GD}-\langle  \bv\cdot\bn, u\rangle _{\GD},
\]
and in particular
\begin{equation}\label{eq:a_wv_wv}
a(w,\bv;w,\bv) =(\nu^{-1}\bv,\bv)=\norm{\bv}^2,
\end{equation}
where $\norm{\cdot}$ denotes the energy norm in $\testq$.
\end{Remark}

We are interested in computing upper and lower bounds for linear functionals of the exact weak solution of \eqref{eq:strongmixed} of the form
\begin{equation}\label{eq:QoI}
s=\ell^O(u,\q)
=(\source^O,u)+\langle \uD^O, \q\cdot\bn\rangle _{\GD}+\langle \trac^O,u\rangle _{\GN},
\end{equation}
for $\source^O\in\mathcal{L}^2(\Omega)$, $\trac^O\in\mathcal{L}^2(\GN)$ and $\uD^O\in\mathcal{C}(\GD)$, namely, compute $\slb , \sub \in\mathbb{R}$ such that
\[
\slb  \leq s \leq \sub.
\]

To compute the bounds, we introduce the corresponding adjoint problem, which in strong form reads:
\begin{equation}\label{eq:strongmixedadjoint}
\begin{array}{rcll}
\qA &\!\!=\!\!&-\nu \grad  \xi  & \text{ in } \Omega,\\
\nabla \cdot \qA &\!\!=\!\!& \source^O & \text{ in } \Omega,\\
\xi &\!\!=\!\!& \uD^O & \text{ on } \GD,\\
\qA\cdot\bn &\!\!=\!\!& -\trac^O & \text{ on } \GN.
\end{array}
\end{equation}
\begin{Remark}
The weak form of the adjoint problem {is: find $(\xi,\qA)\in \test\times \testq$ such that}
\begin{equation}\label{eq:general_variational_formulation_adjoint}
a(w,\bv;\xi,-\qA)= \ell^O(w,\bv) \quad{\forall (w,\bv)\in \test\times\testq.}
\end{equation}
\end{Remark}

\section{Bounds for the Quantity of Interest from general non-orthogonal approximations}
\label{sec:bounds}

Upper and lower bounds for the quantity of interest can be computed given any equilibrated flux and potential reconstructions of the primal and adjoint problem, usually obtained from discrete approximations $(\uh,\qh)$ and $(\uhA,\qhA)$ of  \eqref{eq:general_variational_formulation} and \eqref{eq:general_variational_formulation_adjoint} respectively. The complexity of computing the reconstructions and
evaluating the bounds strongly depends on: 1) the properties of the discrete approximations $(\uh,\qh)$ and $(\uhA,\qhA)$, 2) the kind of data associated {with} the primal and adjoint problems and 3) the desired accuracy of the bounds.
This section presents three different approaches to compute bounds for the quantity of interest $s$. The first approach recovers the bounds by means of computing fully equilibrated fluxes, which in practice can only be used if the data are piecewise polynomial functions. In the second approach, the bounds are recovered \Rem{by} relaxing the equilibration conditions on the fluxes by means of introducing data oscillation errors\Rem{. Finally}, the third approach enhances the  bounds using a Helmholtz decomposition.

\subsection{Bounds from potential and equilibrated flux reconstructions}
\label{subsec:bounds_equilibrated}

Let $(\uh,\qh)$ and $(\uhA,\qhA)$ be two approximations of \eqref{eq:general_variational_formulation} and \eqref{eq:general_variational_formulation_adjoint} respectively. The pairs $(\uh,\qh)$ and $(\uhA,\qhA)$ are said to be potential and equilibrated flux reconstructions of the primal and adjoint problems if the following conditions hold:
\begin{equation}\label{eq:reconstructions}
\begin{tabular}{| llcl |}
\hline
Potential reconstructions:
&$\uh\in\test$ &\quad& $\uhA\in\test$\\[1ex]
&$\uh = \uD \text{ on } \GD$ && $\uhA = \uD^O \text{ on } \GD$
\\[1ex]\hline
Equilibrated flux reconstructions:
&$\qh\in\testq$ && $\qhA\in\testq$\\[1ex]
&$\nabla \cdot \qh = f \text{ in } \Omega$ && $\nabla \cdot \qhA = f^O\text{ in } \Omega$  \\[1ex]
&$\qh\cdot\bn = \trac \text{ on } \GN$ && $\qhA\cdot\bn = -\trac^O  
\text{ on } \GN$
\\\hline
\end{tabular}
\end{equation}

The next result shows that potential and equilibrated flux reconstructions allow computing constant-free bounds for the quantity of interest $s$.

\begin{Theorem}\label{th:property_QoI_01}
Let $(\uh,\qh)$ and $(\uhA,\qhA)$ be two potential and equilibrated flux reconstructions of the primal and adjoint problems satisfying \eqref{eq:reconstructions}. Then
\begin{equation}\label{eq:property_QoI_01}
\pm s 
 \geq
\pm \ell^O(\uh,\qh)
-
\dfrac{1}{2}{\norm{\qh+\nu\grad \uh}\,\norm{\qhA+\nu\grad \uhA}}
\pm \dfrac{1}{2}(\nu^{-1}(\qh+\nu\grad \uh),\qhA-\nu\grad \uhA)
\equiv \pm s_h^\mp,
\end{equation}
and therefore, the quantity of interest $s$ {is} bounded by
\begin{equation*}\label{eq:bounds_QoI_01}
\begin{array}{c}
\displaystyle
s>\slb  = \ell^O(\uh,\qh) 
+\dfrac{1}{2}(\nu^{-1}(\qh+\nu\grad \uh),\qhA-\nu\grad \uhA) 
-\dfrac{1}{2}{\norm{\qh+\nu\grad  \uh}\,\norm{\qhA+\nu\grad  \uhA}}
\\[2ex]
s<\sub  =
 \ell^O(\uh,\qh) 
+\dfrac{1}{2}(\nu^{-1}(\qh+\nu\grad  \uh),\qhA-\nu\grad  \uhA)
+\dfrac{1}{2}{\norm{\qh+\nu\grad  \uh}\,\norm{\qhA+\nu\grad  \uhA}}.
\end{array}
\end{equation*}
\end{Theorem}
The proof of this result is included in \ref{App:eq:property_QoI_01}. 

{
\begin{Remark} Equation \eqref{eq:property_QoI_01} should be interpreted as a shorthand expression for two equations where the $\pm$ and $\mp$ signs and superscripts are linked (each equation obtained by picking all the top/bottom signs/superscripts). Namely, equation \eqref{eq:property_QoI_01} represents the two equations
\begin{equation*}\label{eq:property_QoI_01_part}
\begin{array}{c}
\displaystyle
+ s 
 \geq
+ \ell^O(\uh,\qh)
-
\dfrac{1}{2}{\norm{\qh+\nu\grad \uh}\,\norm{\qhA+\nu\grad \uhA}}
+ \dfrac{1}{2}(\nu^{-1}(\qh+\nu\grad \uh),\qhA-\nu\grad \uhA)
\equiv + s_h^-,\\[2ex]
- s 
 \geq
- \ell^O(\uh,\qh)
-
\dfrac{1}{2}{\norm{\qh+\nu\grad \uh}\,\norm{\qhA+\nu\grad \uhA}}
- \dfrac{1}{2}(\nu^{-1}(\qh+\nu\grad \uh),\qhA-\nu\grad \uhA)
\equiv - s_h^+.
\end{array}
\end{equation*}
This notation is used throughout this paper.
\end{Remark}}

Once the upper and lower bounds for the quantity of interest $s$ are computed, one can compute the bound average 
\[
\sh = \frac{1}{2}(\sub + \slb) = 
\ell^O(\uh,\qh) + \frac{1}{2}(\nu^{-1}(\qh+\nu\grad \uh),\qhA-\nu\grad \uhA),
\]
and the bound gap
\begin{equation}\label{eq:gap}
\Delta_h = \sub - \slb = 
{\norm{\qh+\nu\grad  \uh}\,\norm{\qhA+\nu\grad  \uhA}}.
\end{equation}
The bound average $\sh$ is {seen as} an estimate of the output $s$. {Its} error with respect to $s$ can be easily bounded since
\begin{equation}\label{eq:bounds_QoI_01_mod}
|s- \sh|
\leq \dfrac{1}{2}\Delta_h.
\end{equation}

Mallik et al.  \cite{MR4008052} have recently presented a {result similar to} equation
\eqref{eq:bounds_QoI_01_mod}, but excluding the case of non-homogeneous Neumann boundary conditions. The derivation of the result is done using algebraic manipulations and reiterated use of {the} Cauchy-Schwartz inequality, instead of the reformulation of the output of interest as a constrained minimization problem, see \ref{App:eq:property_QoI_01}. The 
approach {introduced here} {is more general and enables the derivation of} the three improvements \Rem{
described} in the forthcoming sections and 
\Rem{
the extension of this approach} to other problems. 

\subsection{Bounds from potential and zero-order equilibrated flux reconstructions}
\label{subsec:bounds_equilibrated_ZO}

For non-polynomial data, it is not possible in general to find reconstructions 
satisfying \eqref{eq:reconstructions}, and therefore \eqref{eq:property_QoI_01} {cannot} be used to compute guaranteed bounds for the output. 
Fortunately, we can employ the technique described in 
\cite{MR1655854, LADEVEZE2008, MR2981885, MR2899560}
to recover bounds for the energy from projected equilibrated flux reconstructions
by means of introducing data oscillation errors \cite{MR3335498,MR3577961,MR4007990,MR3262938}.

Let $\triang$ be a collection of $d$-dimensional non-overlaping and non-degenerate simplices $\elem$ that partition $\Omega$, such that the intersection of a distinct pair of elements is either an empty set or their common node, edge or face (in three dimensions). Let $\triangE$ denote the set of all its facets, and define 
{
$\projk: \mathcal{L}^2(\elem) \to \mathbb{P}^{\hat{q}}(\elem)$ and $\projg: \mathcal{L}^2(\edge) \to \mathbb{P}^{\bar{q}}(\edge)$ to be the $\mathcal{L}^2(\elem)$
and $\mathcal{L}^2(\edge)$-orthogonal 
projection operators onto $\mathbb{P}^{\hat{q}}(\elem)$ and $\mathbb{P}^{\bar{q}}(\edge)$, respectively.} 
Finally, assume that $\triang$ is such that the data $\nu$ is constant in each element $\elem$, that is $\at{\nu}{\elem}=\nu_{\elem} \in \mathbb{R}$.

Then, the pairs $(\uh,\qhZ)$ and $(\uhA,\qhAZ)$ are said to be potential and zero-order equilibrated flux reconstructions of the primal and adjoint problems if the following conditions hold:
\begin{equation}\label{eq:reconstructions_ZO}
\begin{tabular}{| lll |}
\hline
Potential reconstructions:
&$\uh\in\test$ & $\uhA\in\test$\\[1ex]
&$\uh = \uD \text{ on } \GD$ & $\uhA = \uD^O \text{ on } \GD$
\\[1ex]\hline
\multicolumn{3}{| l |}{Zero-order equilibrated flux reconstructions:}\\[1ex]
&$\qhZ\in\testq$ & $\qhAZ\in\testq$\\[1ex]
\hspace{1.5cm}$\quad\forall\elem\in\triang$
&$(\nabla \cdot \qhZ,1)_\elem = (f,1)_\elem $ & $(\nabla \cdot \qhAZ,1)_\elem = (f^O,1)_\elem$\\[1ex]
\hspace{1.5cm}$\quad\forall\edge\in\triangE\cap\GN$
&$\langle \qhZ\cdot\bn,1\rangle_\edge = \langle \trac,1\rangle_\edge$ & $\langle \qhAZ\cdot\bn,1\rangle_\edge = \langle -\trac^O,1\rangle_\edge$
\\\hline
\end{tabular}
\end{equation}

{Note that the relaxation of the equilibrium conditions affect only fluxes, and \Rem{that} the conditions on the potentials are not weaker than in \eqref{eq:reconstructions}. Assuming that the conditions on the potentials $\uh$ and $\uhA$ are exact is not a strong restriction}
because any approximation can be easily modified on the Dirichlet boundary to exactly \Rem{
satisfy} the Dirichlet boundary conditions. This simplified approach can be considered here since the potential and flux reconstructions necessary to compute the bounds for $s$ are completely independent, as opposed to what occurs in other existing more involved approaches, see for instance \cite{MR2899560}.

If the bounds for the output are computed using zero-order equilibrated fluxes, the bounding property presented in \eqref{eq:property_QoI_01} is lost in general. The next result, proved in  \ref{App:eq:bounds_projected}, introduces a workaround to replace the exactly equilibrated fluxes reconstructions by its zero-order peers by means of introducing data oscillations errors. Indeed, constant-free bounds for the quantity of interest $s$ can be computed from potential and zero-order equilibrated flux reconstructions.

\begin{Theorem}\label{th:property_QoI_02}
Let $(\uh,\qhZ)$ and $(\uhA,\qhAZ)$ be two potential and zero-order equilibrated flux reconstructions of the primal and adjoint problems satisfying \eqref{eq:reconstructions_ZO} and $\kappa\in(0,+\infty)$ be an arbitrary {scaling} parameter. Then
\begin{equation}\label{eq:property_QoI_02}
\pm s 
 \geq
\pm(\source^O,\uh)
\pm\langle \trac^O,\uh\rangle _{\GN}
\pm (\source,\uhA)
\mp\langle \trac,\uhA\rangle _{\GN}
\mp(\nu\grad  \uh,\grad  \uhA)
-\dfrac{1}{4\kappa}\sum\limits_{\elem\in\triang} (\eta_K^{0\mp})^2,
\end{equation}
for
\begin{equation}\label{eq:etaK_ZO}
\begin{array}{rl}
\eta_K^{0\mp}
&=
{\norm{\pm (\qhAZ +\nu\grad \uhA) -\kappa (\qhZ+\nu\grad \uh)}}_K 
	+
	C_1 \nu_K^{-1/2}
	||\pm(\source^O-\nabla\cdot \qhAZ )-\kappa(\source-\nabla\cdot \qhZ)||_{\mathcal{L}^2(K)}	
\\[1ex]
&\displaystyle	
	+ \sum\limits_{\edge\in \GN \cap \partial\elem}
	\CTb \nu_K^{-1/2}||\mp(\trac^O+\qhAZ\cdot\bn)-\kappa (\trac-\qhZ\cdot\bn)||_{\mathcal{L}^2(\edge)},
\end{array}	
\end{equation}
{where $||\cdot||_{\mathcal{L}^2(\elem)}$ denotes the $\mathcal{L}^2(\elem)$ norm both in $\mathbb{R}$ and $\mathbb{R}^d$, $\norm{\cdot}_\elem $ is the restriction of the energy norm defined in \eqref{eq:a_wv_wv} to element $\elem$ and} the values for the constants $C_1$ and $\CTb$ are given in \ref{App:eq:bounds_projected}, equation \eqref{eq:C1C2}.
\end{Theorem}

\begin{Remark} The bounds provided in expression \eqref{eq:property_QoI_02} {coincide with} the bounds introduced in \eqref{eq:property_QoI_01} if
 $\qhZ=\qh$ and $\qhAZ=\qhA$ are exact equilibrated flux reconstructions and  one considers
$\kappa = \kappa_{\rm opt} = {\norm{\qhA +\nu\grad   \uhA }/\norm{\qh+\nu\grad  \uh}}$.
\end{Remark}

\begin{Remark}
The bounds given by \eqref{eq:property_QoI_02} are less accurate than the previously introduced in \eqref{eq:property_QoI_01} since they rely on the local Poincaré inequality, a trace inequality and reiterated applications of the Cauchy-Schwartz inequality. Therefore, if possible, the a posteriori error estimation technique should minimize the data oscillation errors included in 
$||\pm(\source^O-\nabla\cdot \qhAZ )-\kappa(\source-\nabla\cdot \qhZ)||_{\mathcal{L}^2(K)}$ and $||\mp(\trac^O+\qhAZ\cdot\bn)-\kappa (\trac-\qhZ\cdot\bn)||_{\mathcal{L}^2(\edge)}$.
\end{Remark}

\subsection{Bounds from potential and projected equilibrated flux reconstructions}
\label{subsec:bounds_equilibrated_projected}

In order to minimize the influence of the data oscillation errors \Rem{
and} to obtain computable expressions for the equilibrated flux reconstructions if the data for the problem are not piecewise polynomial fields, it is standard to introduce an intermediate step between the generally uncomputable exact equilibrated fluxes given by \eqref{eq:reconstructions} and the zero-order equilibrated fluxes given by \eqref{eq:reconstructions_ZO}. Indeed, $(\uh,\qhP)$ and $(\uhA,\qhAP)$ are said to be potential and projected equilibrated flux reconstructions of the primal and adjoint problems associated {with} the constant pair $(\hat{q},\bar{q})$ if the following conditions hold:

\begin{equation}\label{eq:reconstructions_Proj}
\begin{tabular}{| lll |}
\hline
Potential reconstructions:
&$\uh\in\test$ & $\uhA\in\test$\\[1ex]
&$\uh = \uD \text{ on } \GD$ & $\uhA = \uD^O \text{ on } \GD$
\\[1ex]\hline
Projected equilibrated flux reconstructions:
&$\qhP\in\testq$ & $\qhAP\in\testq$\\[1ex]
\hspace{5cm}$\quad\forall\elem\in\triang$
&$\nabla \cdot \qhP|_{\elem} = \projk \source$ & $\nabla \cdot \qhAP|_{\elem} = \projk \source^O $\\[1ex]
\hspace{5cm}$\quad\forall\edge\in\triangE\cap\GN$
&$\qhP\cdot\bn|_{\edge} =\projg \trac $ & $\qhAZ\cdot\bn|_{\edge}=-\projg\trac^O $
\\\hline
\end{tabular}
\end{equation}

In this case, bounds for the quantity of interest {are} obtained from \eqref{eq:property_QoI_02}
 where now the local elementary contributions {read}
\begin{equation}\label{eq:etaK_Proj}
\begin{array}{rl}
\eta_K^{\pi\mp}
&=
{\norm{\pm (\qhAP +\nu\grad \uhA) -\kappa (\qhP+\nu\grad \uh)}}_K 
	+
	C_1 \nu_K^{-1/2}
	||\pm(\source^O-\projk \source^O )-\kappa(\source-\projk \source)||_{\mathcal{L}^2(K)}	
\\[1ex]
&\displaystyle	
	+ \sum\limits_{\edge\in \GN \cap \partial\elem}
	\CTb \nu_K^{-1/2}||\mp(\trac^O-\projg\trac^O)-\kappa (\trac-\projg \trac) ||_{\mathcal{L}^2(\edge)},
\end{array}	
\end{equation}
and moreover, if $\hat{q}$ and $\bar{q}$ are greater or equal than $\text{degree}\{\uhA\}$  
then
\[
\pm(\source^O,\uh)
\pm\langle \trac^O,\uh\rangle _{\GN}
\pm (\source,\uhA)
\mp\langle \trac,\uhA\rangle _{\GN}
\mp(\nu\grad  \uh,\grad  \uhA)
=\pm\ell^O(\uh,\qhP)\mp (\qhP+\nu\grad  \uh,\grad  \uhA)
\]
yielding the alternative form of the bounds
\begin{equation}\label{eq:property_QoI_02_proj}
\pm s 
 \geq
\pm\ell^O(\uh,\qhP)\mp (\qhP+\nu\grad  \uh,\grad  \uhA)
-\dfrac{1}{4\kappa}\sum\limits_{\elem\in\triang} (\eta_K^{\pi\mp})^2.
\end{equation}

\begin{Remark} Finding the optimal value of $\kappa$ minimizing the bounds given in \eqref{eq:property_QoI_02} either for the expression or the local estimate $\eta_K$ given \eqref{eq:etaK_Proj} or \eqref{eq:etaK_ZO} is not trivial. Therefore, it is usual to use the value $\kappa = {\norm{\qhAP +\nu\grad   \uhA}/\norm{\qhP+\nu\grad  \uh}}$
or $\kappa = {\norm{\qhAZ +\nu\grad   \uhA}/\norm{\qhZ+\nu\grad  \uh}}$
 that optimizes the bounds assuming that no data oscillation errors are present.
\end{Remark}

\subsection{Exact representation for the quantity of interest - enhancement of the bounds using lower bounds for the energy}
\label{subsec:bounds_equilibrated_enhanc}

To improve the quality of the bounds given in the previous sections, the following result providing an exact representation for the quantity of interest can be used, see \ref{App:eq:property_QoI_02} for its proof.
\begin{Theorem}\label{th:property_QoI_03}
For any $(\uh,\qh)$ and $(\uhA,\qhA)$ in $\test\times\testq$ and $\kappa\in(0,+\infty)$ the following exact representation for the quantity of interest holds\begin{equation}\label{eq:exactQoI}
\begin{array}{ll}
\pm s
&=
\pm \hat{s}_h^\mp
+\dfrac{1}{4\kappa}{\norm{\bm{e}_{\bm{\varphi}}^\mp-\nu\grad  e_\phi^\mp}}^2 
+(u,\pm (\source^O - \nabla\cdot\qhA ) -\kappa(\source - \nabla\cdot\qh))
\\[1ex]
&
-\langle u,\mp (\trac^O + \qhA  \cdot\bn) -\kappa(\trac - \qh \cdot\bn)\cdot\bn\rangle _{\GN}
-\langle  \q\cdot\bn, \mp(\uD^O - \uhA) - \kappa (\uD - \uh)\rangle _{\GD},
\end{array}	
\end{equation}
where 
\[
\pm \hat{s}_h^\mp
=\pm\ell(\uhA,-\qhA )
+\kappa\ell(\uh,\qh) -\dfrac{1}{4\kappa}{\norm{\pm (\qhA +\nu\grad   \uhA )-\kappa (\qh-\nu\grad  \uh)}}^2,
%
\]
and
$e_\phi^\mp 
=\mp(\xi - \uhA) - \kappa (u - \uh)$ and 
$
\bm{e}_{\bm{\varphi}}^\mp 
= \pm (\qA - \qhA ) -\kappa(\q - \qh).
$
\end{Theorem}
\begin{Remark}
In the case where $(\uh,\qh)$ and $(\uhA,\qhA)$ in $\test\times\testq$ are potential and equilibrated flux reconstructions of the primal and adjoint problems, $\tilde{s}_h^{\mp}$ coincides with $s_h^{\mp}$, expressed in two different forms in equations \eqref{eq:property_QoI_01_mod_mod} and \eqref{eq:app:smph_alternative}.
\end{Remark}

Many a posteriori error estimation techniques can be derived from this exact representation of the quantity of interest. For instance, it is possible to devise error estimators incorporating possible errors in the Dirichlet boundary conditions, error estimators incorporating the data oscillation errors outside $\hat{s}_h^\mp$ in contrast to the strategy described in Section \ref{subsec:bounds_equilibrated_ZO}, or error estimators incorporating the term ${\norm{\bm{e}_{\bm{\varphi}}^\mp-\nu\grad  e_\phi^\mp}}^2$ in the final bounds.

Here, this expression {is only} used to introduce two enhancements of the bounds. The first error estimation technique derived from Theorem \ref{th:property_QoI_03} is summarized in Remark \ref{Rem:QoI_ub_lb}. This technique mimics the standard expression used in a posteriori error estimation to compute bounds for quantities of interest for standard Galerkin orthogonal finite element approximations. That is, it allows obtaining bounds for the quantity of interest by means of computing upper and lower bounds for the energy norm.

\begin{Remark}\label{Rem:QoI_ub_lb}
Let $\uh$ and $\uhA$ be two potential reconstructions of $u$ and $\xi$ respectively and consider
$\qh = \q =-\nu \grad  u$ and $\qhA = \qA=-\nu \grad  \xi$. Noting that in this case $\tilde{s}_h^{\mp} = s_h^{\mp}$, so that $\tilde{s}_h^{\mp}$ can be rewritten as shown in equation \eqref{eq:app:smph_alternative}, the exact representation for the quantity of interest \eqref{eq:exactQoI} yields after some rearrangements to
\begin{equation}\label{eq:exact_QoI_standard_FE}
\begin{array}{ll}
\pm s
& =
\pm(\source^O,\uh)
\pm\langle \trac^O,\uh\rangle _{\GN}
\pm (\source,\uhA)
\mp\langle \trac,\uhA\rangle _{\GN}
\mp(\nu\grad  \uh,\grad  \uhA)
\\[1ex]
&
-\dfrac{1}{4\kappa}{\norm{\nu\grad( \uA -   \uhA \mp\kappa( u-  \uh))}}^2
+\dfrac{1}{4\kappa}{\norm{\nu\grad  (\xi - \uhA \pm \kappa (u - \uh))}}^2, 
\end{array}
\end{equation}
and therefore bounds for the quantity of interest may be recovered computing upper and lower bounds for the energy norm of the adequate combined primal/adjoint problems as
\begin{equation}\label{eq:QoI_ub_lb}
\begin{array}{ll}
\pm s
& \geq
\pm(\source^O,\uh)
\pm\langle \trac^O,\uh\rangle _{\GN}
\pm (\source,\uhA)
\mp\langle \trac,\uhA\rangle _{\GN}
\mp(\nu\grad  \uh,\grad  \uhA)
\\[1ex]
&
-\dfrac{1}{4\kappa}{\norm{\nu\grad( \uA -   \uhA \mp\kappa( u-  \uh))}}^2_{_\textup{UB}}
+\dfrac{1}{4\kappa}{\norm{\nu\grad  (\xi - \uhA \pm \kappa (u - \uh))}}^2_{_\textup{LB}}. 
\end{array}
\end{equation}
In fact, expanding the norms appearing in equation \eqref{eq:exact_QoI_standard_FE} allows obtaining the following exact expression for the quantity of interest
\begin{equation}\label{eq:exactQoI_nodual}
s
=
(\source^O,\uh)
+\langle \trac^O,\uh\rangle _{\GN}
+(\source,\uhA)
-\langle \trac,\uhA\rangle _{\GN}
-(\nu\grad  \uh,\grad  \uhA)
+(\nu\grad( \uA -   \uhA), \grad( u-  \uh))
\end{equation}
from where equation \eqref{eq:exact_QoI_standard_FE}
{is} be recovered back using the standard parallelogram identity applied to the last scalar product.
\end{Remark}


The second technique devised from Theorem \ref{th:property_QoI_03} incorporates the error in the term
${\norm{\bm{e}_{\bm{\varphi}}^\mp-\nu\grad  e_\phi^\mp}}^2$ in the final expression of the bounds. For simplicity of presentation, this technique is only described assuming that no data oscillation errors are present, that is, assuming it is possible to compute $(\uh,\qh)$ and $(\uhA,\qhA)$ being potential and equilibrated flux reconstructions of the primal and adjoint problems satisfying \eqref{eq:reconstructions}. In this case, the quantity of interest
{is} rewritten using equation \eqref{eq:exactQoI} as 
\begin{equation}\label{eq:exactQoI_nodataoscillations}
\begin{array}{ll}
\pm s
&=
\pm s^\mp_h
+\dfrac{1}{4\kappa}{\norm{\bm{e}_{\bm{\varphi}}^\mp-\nu\grad  e_\phi^\mp}}^2
\geq 
\pm s^\mp_h
+\dfrac{1}{4\kappa}{\norm{\bm{e}_{\bm{\varphi}}^\mp-\nu\grad  e_\phi^\mp}}^2_{_\textup{LB}},
\end{array}	
\end{equation}
and therefore, the bounds can be improved by introducing a lower bound of the energy norm of 
$\bm{e}_{\bm{\varphi}}^\mp-\nu\grad  e_\phi^\mp$.
These lower bounds {are} incorporated using the result detailed in  \ref{App:LB}. Indeed, the following representation holds
\begin{equation}\label{eq:exactQoI_nodataoscillations_LB}
\begin{array}{ll}
\pm s
&=
\pm s^\mp_h
+\hspace{-.4cm}
\sup\limits_{\scriptsize\begin{array}{c}w^\mp\in H^1_0(\Omega)\\[0ex]\psi^\mp\in [H^1(\Omega)]^{2d-3}\end{array}}
\hspace{-.4cm}
\dfrac{1}{4\kappa}\dfrac{(\ell_\times^\mp(w^\mp,\psi^\mp))^2}{
{\norm{\nu\grad w^\mp + \curl \psi^\mp}}^2}
\geq
\pm s^\mp_h
+\dfrac{1}{4\kappa}
\dfrac{(\ell_\times^\mp(w^\mp,\psi^\mp))^2}{
{\norm{\nu\grad w^\mp + \curl \psi^\mp}}^2},
\end{array}	
\end{equation}
for any $w^\mp\in H^1_0(\Omega)$, $\psi^\mp\in [H^1(\Omega)]^{2d-3}$, where
\[
\ell_\times^\mp(w^\mp,\psi^\mp)=
\mp(\nu^{-1}(\qhA +\nu\grad  \uhA),\nu\grad w^\mp+\curl \psi^\mp)
- \kappa( \nu^{-1}(\qh+\nu\grad \uh),\nu\grad w^\mp-\curl \psi^\mp)
\]
and $\grad\times$ denotes the standard curl operator, see \cite{MR851383}.

\section{Bounds for the Quantity of Interest using the Hybridizable Discontinuous Galerkin Method}
\label{sec:HDG}

This section details how to compute bounds for a quantity of interest using the HDG method introduced in \cite{MR2513831} as a means to {obtain} the approximations $(\uh,\qh)$ and $(\uhA,\qhA)$ in $\test\times\testq$ of the primal and adjoint problems. For simplicity, {only} the construction of the potential and equilibrated flux reconstructions  for the primal problem are described. The constructions for the adjoint problem are analogous.

\subsection{Notations and HDG approximation}

To introduce the HDG approximation of \eqref{eq:general_variational_formulation}, some notations have to be introduced, see \cite{MR2513831}.

Let $\triang$ be a disjoint partition of $\Omega$, see Section \ref{subsec:bounds_equilibrated_ZO}, and consider the set of all its facets $\triangE=\triangE^o \cup \triangE^\partial$, where $\triangE^\partial$ consists of the facets lying on the boundary $\partial\Omega$, and  $\triangE^o$ are the remaining interior facets.
Also denote by $\partial\triang$ the mesh skeleton $\{\partial \elem: \elem \in\triang\}$. 


Given two elements $\elem^+$ and $\elem^-$ of $\triang$ sharing a common facet $e=\partial\elem^+\cap \partial\elem^-\in\triangE^o$, let $\bn^+$ and $\bn^-$ be the outward unit normals to $\elem^+$ and $\elem^-$, respectively, and let $(\q^\pm,u^\pm)$ be the traces of $(\q,u)$ on $\edge$ from the interior of $\elem^\pm$, that is $\q^\pm = \at{\q}{\elem^\pm}$ and $u^\pm=\at{u}{\elem^\pm}$. Then, we define the mean values $\mean{\cdot}$ and jumps $\jump{\cdot}$ as follows. For $e=\partial\elem^+\cap \partial\elem^-\in\triangE^o$, we set
\[
\begin{array}{ll}
\mean{\q} = (\q^+ + \q^-)/2 & \mean{u} = (u^+ + u^-)/2 \\
\jump{\q\cdot\bn} = \q^+\cdot\bn^+ + \q^-\cdot\bn^- & \jump{u\bn} = u^+\bn^+ + u^-\bn^-,
\end{array}
\]
whereas for $\edge\in\triangE^\partial$, the set of boundary facets on which $\q$ and $u$ are single valued, we set
\[
\mean{\q} = \q \quad,\quad \mean{u} = u \quad,\quad
\jump{\q\cdot\bn} = \q\cdot\bn \quad,\quad \jump{u\bn} = u\bn,
\]
where $\bn$ is the outward normal to $\partial\Omega$. Note that the jump in $u$ is a vector, but the jump in $\q$ is a scalar which only involves the
normal component of $\q$. Furthermore, the jump will be zero for a continuous function.

The discontinuous finite dimensional spaces $\testh$ and $\testqh$ are defined by
\[
\begin{array}{l}
\testh=\{w\in \mathcal{L}^2(\Omega) :  \at{w}{\elem }\in \mathbb{P}^p(\elem )  \,\,\forall \elem \in\triang \},\\[1ex]
\testqh=\{\bv\in [\mathcal{L}^2(\Omega)]^d : \at{\bv}{\elem }\in [\mathbb{P}^p(\elem )]^d \,\,\forall \elem \in\triang \},\\[1ex]
\testeh = \{\mu \in \mathcal{L}^2(\triangE) :  \at{\mu}{\edge}\in \mathbb{P}^p(\edge)  \,\,\forall \edge\in\triangE \},\\[1ex]
\testeh(\uD) = \{\mu \in \testeh :  \at{\mu}{\edge} =\proj_\edge^{p} \uD   \,\,\forall \edge\in\triangE^\partial\cap\GD \},
\end{array}
\]
where $\mathbb{P}^p(D)$ denotes the set of polynomials of degree at most $p\geq 0$ on $D$ and $\proj_\edge^{p}$ denotes the $\mathcal{L}^2$ projection defined in Section \ref{subsec:bounds_equilibrated_ZO}.

Finally, let
\[
(w,v)_{\triang} = \sum\limits_{\elem\in\triang} (w,v)_{\elem} 
\quad,\quad
\langle \zeta,\rho\rangle _{\partial\triang} = \sum\limits_{\elem\in\triang} \langle \zeta,\rho\rangle _{\partial\elem} 
\quad,\quad
\langle \mu,\nu\rangle _{\triangE} = \sum\limits_{\edge\in\triangE} \langle \mu,\nu\rangle _{\edge} 
\]
for scalar or vector functions $w, v$ defined on $\triang$, $\zeta, \rho$ defined on $\partial\triang$ and $\mu, \nu$ on $\triangE$.

The HDG method seeks an approximation $(\uhHDG,\qhHDG)\in \testh\times \testqh$ to the exact solution  $(u,\q)\in \test\times \testq$  such that for all $\elem\in\triang$
\begin{equation}\label{eq:HDG_weak01}
\begin{array}{ll}
(\nu^{-1}\qhHDG,\bv)_{\elem}
-(\uhHDG,\nabla\cdot\bv)_{\elem}
+\langle \uhHDGT,\bv\cdot\bn\rangle _{\partial\elem} = 0 & \forall \bv\in[\mathbb{P}^p(\elem )]^d
\\[2ex]
-(\qhHDG,\nabla w)_{\elem}
+\langle \qhHDGT\cdot\bn,w\rangle _{\partial\elem} 
= (\source,w)_{\elem}& \forall w\in \mathbb{P}^p(\elem ),
\end{array}
\end{equation}
where the numerical traces are defined as
\begin{equation}\label{eq:num_traces}
\begin{array}{ll}
 \displaystyle
\uhHDGT =
\frac{\tau^+}{\tau^++\tau^-}\uhHDG^+
+\frac{\tau^-}{\tau^++\tau^-}\uhHDG^-
+\frac{1}{\tau^++\tau^-}\jump{\qhHDG\cdot\bn}  & \text{ on } \triangE^o,
\\[1ex]
\displaystyle
\qhHDGT=
\frac{\tau^-}{\tau^++\tau^-}\qhHDG^+
+\frac{\tau^+}{\tau^++\tau^-}\qhHDG^-
+\frac{\tau^+\tau^-}{\tau^++\tau^-}\jump{\uhHDG\bn}  & \text{ on } \triangE^o,
\\[2ex]
\displaystyle
\uhHDGT =\proj_\edge^{p} \uD
\qquad,\qquad
\qhHDGT \cdot \bn=
\qhHDG\cdot \bn
+\tau(\uhHDG - \proj_\edge^{p}\uD)
 & \text{ on } \triangE^\partial\cap\GD,
\\[1ex]
\displaystyle
\qhHDGT \cdot \bn=\proj_\edge^{p}\trac \hspace{-.12cm}
\qquad,\qquad
\uhHDGT = \uhHDG
+\frac{1}{\tau}(\qhHDG\cdot\bn - \proj_\edge^{p}\trac)
 & \text{ on } \triangE^\partial\cap\GN,
\end{array}
\end{equation}
and $\tau$ is the {strictly positive stabilization parameter which plays a crucial role on the stability, accuracy and convergence properties of the HDG method, see for instance \cite{MR2556564,MR2513831}. 
The stabilization function $\tau$ is defined for each element $\elem\in\triang$ so that $\tau^+$ and $\tau^-$ denote its restriction to elements $\elem^+$ and $\elem^-$ respectively, namely $\tau^\pm=\at{\tau}{\elem^\pm}$. Note that for each facet $e=\partial\elem^+\cap \partial\elem^-\in\triangE^o$, in general, $\at{\tau^+}{e}\neq \at{\tau^-}{e}$ so that the stabilization parameter is double-valued on $\triangE^o$.
}

As shown in \cite{MR2513831}, the distinctive feature of the HDG method is that both $\uhHDG$ and $\qhHDG$ converge with the optimal order $p + 1$ in the $\mathcal{L}^2$-norm. Moreover, it is shown that  $\uhHDG$ and  $\uhHDGT$ superconverge with order $p+2$ to some $\mathcal{L}^2$-like projections of the exact variable $u$. As a consequence, a post-processing of the approximate solution provides an approximation of the potential converging with order $p + 2$.
We can see from \eqref{eq:num_traces} that the HDG method belongs to a family of DG methods whose numerical traces are of the form
\[
\begin{array}{ll}
\displaystyle
\uhHDGT = \mean{\uhHDG} - \bC_{12} \jump{\uhHDG\bn} + C_{22} \jump{\qhHDG\cdot\bn}  & \text{ on } \triangE^o,
\\[1ex]
\displaystyle
\qhHDGT=\mean{\qhHDG} + \bC_{12} \jump{\qhHDG\cdot\bn} + C_{11} \jump{\uhHDG\bn}
 & \text{ on } \triangE^o,
\end{array}
\]
where the penalization parameters are such that $|\bC_{12}|$ is finite, $C_{11} > 0$, and $C_{22}\geq 0$. This family of DG methods were studied first in \cite{MR1655854} and more thoroughly in \cite{MR2429868, MR2448694}, wherein it was shown that if one chooses $C_{22} \sim 1/C_{11}$, then both $\uhHDG$ and $\qhHDG$ converge in $\mathcal{L}^2$-norm with the optimal order $p+1$. Since the HDG method satisfies this condition for any value of $\tau$ such that $\tau^+ = \tau^- > 0$, the method possesses the optimal and super-convergence properties as mentioned above. Note that, 
for some other DG methods such as the LDG method \cite{MR1655854}  with $C_{22} = 0$, the approximate flux $\qhHDG$ converges with order $p$ in  $\mathcal{L}^2$-norm, which is suboptimal. 
 Furthermore, one can show that the HDG method is consistent, adjoint consistent, and locally and globally conservative by following the analysis given in \cite{MR1951054,MR1885715}.

\begin{Remark} The weak problem given by \eqref{eq:HDG_weak01} and 
\eqref{eq:num_traces} is equivalent to the following alternative weak formulation: find $(\uhHDG,\qhHDG)\in \testh\times \testqh$ such that
\begin{equation}\label{eq:HDG_weak02}
\aHDG(\uhHDG,\qhHDG;w,\bv) = \ellHDG(w,\bv) \qquad \forall (w,\bv)\in \testh\times \testqh
\end{equation}
where
\[
\begin{array}{l}
\aHDG(u,\q;w,\bv) =
{(\nu^{-1}\q,\bv)_{\triang }
-(u,\nabla\cdot\bv)_{\triang }
-(\q,\nabla w)_{\triang }}
+\langle \widehat{u},\bv\cdot\bn\rangle _{\partial\triang \backslash \GD}
+\langle \widehat{\q}\cdot\bn,w\rangle _{\partial\triang \backslash\GN} 
 \\[2ex] 
\ellHDG(w,\bv)  =
 (\source,w)_{\triang }
 -\langle \uD,\bv\cdot\bn\rangle _{\GD}
 -\langle \trac,w\rangle _{\GN}.
 \end{array}
\]
Moreover, if $(u,\q)$ is the solution of \eqref{eq:general_variational_formulation} and $(w,\bv)\in\test\times\testq$ then,
\[
\aHDG(u,\q;w,\bv) = a(u,\q;w,\bv)+\frac{1}{\tau}\langle 
\trac - \proj_\edge^{p}\trac,\bv\cdot\bn\rangle _{\GN}
+\tau\langle \uD - \proj_\edge^{p}\uD,w\rangle _{\GD} 
\] 
and $\ellHDG(w,\bv)  = \ell(w,\bv)$, and since $\testh\not\subset \test , \testqh\not\subset \testq$,
the approximation $(\uhHDG,\qhHDG)\in \testh\times \testqh$ can be seen as a non-conforming approximation of the exact solution $(u,\q)\in\test\times\testq$ such that
\[
\aHDG(u,\q;w,\bv) = \ellHDG(w,\bv)
-\frac{1}{\tau}\langle 
\trac - \proj_\edge^{p}\trac,\bv\cdot\bn\rangle _{\GN}
-\tau\langle \uD - \proj_\edge^{p}\uD,w\rangle _{\GD} 
 \qquad \forall (w,\bv)\in \test\times \testq.
\]
\end{Remark}

\begin{Remark} The numerical traces $\uhHDGT$ and $\qhHDGT$ defined in \eqref{eq:num_traces} are single-valued functions for each edge $\edge\in\triangE^o$ and verify $\at{\qhHDGT \cdot \bn}{\GN}=\proj_\edge^{p}\trac$, $
\at{\uhHDGT}{\GD}=\proj_\edge^{p} \uD$. Moreover, from equation \eqref{eq:HDG_weak01} it holds that {for all} $\elem\in\triang$
\begin{equation}\label{eq:qhHDGT_ZOeq}
\langle \qhHDGT\cdot\bn,1\rangle _{\partial\elem} 
= (\source,1)_{\elem}.
\end{equation}
\end{Remark}

\subsection{HDG projected equilibrated flux reconstruction}
\label{subsec:flux_rec}

As described in \cite{MR2513831}, thanks to the single-valuedness of the normal component of the numerical trace $\qhHDGT$ and using \eqref{eq:qhHDGT_ZOeq} it is possible to recover a projected equilibrated flux reconstruction $\qhP\in\testq$ using an element-by-element procedure and converging in an optimal fashion. Indeed, let $RT^p(\elem)=[\mathbb{P}^p(\elem)]^d +  \bm{x}\,\mathbb{P}^p(\elem)$ be the Raviart-Thomas finite element space of order $p$, see \cite{MR0483555,MR592160,MR1115205}, and let $\qhP\in\testq$ be the post-processed flux defined in \cite{MR2513831}, namely for each element $\elem\in\triang$, $\at{\qhP}\elem\in RT^p$ is such that
\begin{equation}\label{eq:qhP}
\begin{array}{ll}
\langle (\qhP-\qhHDGT)\cdot\bn,\mu\rangle_\edge = 0
&\forall \mu\in  \mathbb{P}^p(\edge)\quad,\quad
\forall \edge\in\partial\elem,
\\[1ex]
(\qhP-\qhHDG,\bv)_\elem = 0
&\forall \bv \in  [\mathbb{P}^{p-1}(\elem)]^d \quad,\quad \text{ if } p\geq 1.
\end{array}
\end{equation}
Then, the projected conditions \eqref{eq:reconstructions_Proj} are satisfied for $\hat{q}=\bar{q}=p$. {This is proven by letting} $w\in \mathbb{P}^p(\elem)$ and $\mu\in  \mathbb{P}^p(\edge)$. Since $\grad w\in [\mathbb{P}^{p-1}(\elem)]^d$ and $\at{w}{e}\in \mathbb{P}^p(\edge)$, using equation \eqref{eq:green} with $\omega=K$, and equations \eqref{eq:qhP}, \eqref{eq:HDG_weak01} and \eqref{eq:num_traces} it holds that
\begin{equation*}
\begin{array}{ll}
(\nabla \cdot \qhP,w)_\elem
&
=\langle\qhP\cdot \bn,w\rangle_{\partial\elem}-(\qhP ,\grad w)_\elem
\\[1ex]
&
=\langle\qhHDGT\cdot \bn,w\rangle_{\partial\elem}-(\qhHDG ,\grad w)_\elem
=(\source,w)_{\elem}
=(\proj_K^{p} \source,w)_{\elem}
\end{array}
\end{equation*}
and 
\begin{equation*}
\langle \qhP\cdot\bn,\mu\rangle_\edge 
=\langle \qhHDGT\cdot\bn,\mu\rangle_\edge 
=\langle \proj_\edge^{p}\trac,\mu\rangle_\edge
\qquad \text{for $e\in\partial\elem\cap\GN$},
\end{equation*}
which concludes the proof using that $\nabla \cdot \qhP\in \mathbb{P}^p(\elem)$ and $\qhHDGT\cdot\bn\in\mathbb{P}^p(\edge)$.

\subsection{HDG potential reconstruction}
\label{subsec:pot_rec}

A potential reconstruction $\uh$ {is} computed taking into account the single-valuedness of the numerical trace $\uhHDGT$ or {alternatively,} by simply averaging $\uhHDG$. However, the equilibrated projected flux reconstruction $\qhP$ converges with order $p+1$, {and therefore} optimal convergence for the quantity of interest is only achieved if the potential reconstruction $\uh$ superconverges with order $p+2$. Luckily, the {post-processed scalar variable $\uhHDG^\ast\in \mathcal{W}_h^{p+1}$ introduced in Section 4.2 of  \cite{MR2513831} can be used to achieve this desired superconvergence, namely}
\begin{equation}\label{eq:uhHDGast}
\begin{array}{ll}
(\nabla \uhHDG^\ast,\nabla w)_\elem
=(\nabla \cdot \qhP ,w)_\elem-(\qhP\cdot \bn,w)_{\partial\elem} 
\qquad \forall w\in  \mathcal{W}_h^{p+1},
\\[1ex]
(\uhHDG^\ast,1)_{\elem} = (\uhHDG,w)_{\elem}.
\end{array}
\end{equation}
Then the \emph{continuous} potential reconstruction $\uh\in\mathcal{W}_h^{p+1}\cap\test$ is recovered using a simple averaging of $\uhHDG^\ast$ at the element interfaces and exactly enforcing the Dirichlet boundary conditions. {Recall that the condition regarding the values of $\uh$ on $\GD$ are exact, namely $\uh = \uD$ on $\GD$ \eqref{eq:reconstructions_Proj}. Therefore, on the edges $e\in\triangE^\partial\cap\GD$ for which $\at{\uD}{e}\in \mathbb{P}^{p+1}(\edge)$, the nodal values of $\uh$ lying on $\GD$ are modified to match $\uD$. Otherwise, local extension operators are used to exactly enforce the boundary conditions.
}


\subsection{Local optimization of the bounds}\label{subsec:localopt}

The quality of the bounds for the quantity of interest {is} measured using the bound gap introduced in \eqref{eq:gap}. Therefore, the optimal reconstructions are the ones minimizing ${\norm{\qh+\nu\grad  \uh}}$. Thanks to the single-valuedness of the numerical traces $\uhHDGT$ and $\qhHDGT$ one could recover $\uh\in\test$ and $\qh\in\test$ verifying \eqref{eq:reconstructions_Proj} by first averaging $\uhHDGT$ at the mesh vertices and then using a constrained local optimization procedure in each element. This strategy, however, does not provide optimal convergence for the quantity of interest {because} it does not recover a superconvergent potential reconstruction $\uh$. 

However, once the flux and potential reconstructions are obtained using the strategies described in Subsections \ref{subsec:flux_rec} and \ref{subsec:pot_rec}, an extra local minimization procedure can be performed in each element to improve the bounds. Indeed, let $\uh\in\test$ and  $\qhP\in\testq$ be the reconstructions defined in the aforementioned subsections.
Then for each element of the mesh, the improved value for the reconstructions is computed as: find $\at{\uh^\ast}{\elem}\in \mathbb{P}^{p+1}(\elem)$ 
and  $\at{(\qhP)^\ast}{\elem}\in[\mathbb{P}^{p+1}(\elem)]^d$ minimizing 
${\norm{(\qhP)^\ast+\nu\grad  \uh^\ast}}_K$ such that
\begin{equation*}
\begin{array}{ll}
\nabla \cdot (\qhP)^\ast = \projk \source & \text{ in } \elem\\
(\qhP)^\ast\cdot\bn =\qhP\cdot\bn & \text{ on } \partial\elem\\
\uh^\ast =\uh & \text{ on } \partial\elem.
\end{array}
\end{equation*}
It is worth noting that this improvement is only relevant for large values of $p$ where the degrees of freedom are not concentrated on the boundaries. Also, the local interpolation degree of $\uh^\ast$ and $\at{(\qhP)^\ast}{\elem}$ {could} be increased but no gain on the global convergence {would be} obtained.

{A summary of the procedure devised above}
to determine the bounds for $s$ from HDG approximations of the primal and adjoint {problems is shown} in Figure \ref{box:bounds}.

\begin{figure}[h!]
\noindent\framebox{ \begin{minipage}[t]{.93\textwidth}
\begin{itemize}
\item[0.-] Compute the HDG approximations of the primal and adjoint problems $(\uhHDG,\qhHDG)$ and $(\uhHDGA,\qhHDGA)\in \testh\times \testqh$ such that $\forall (w,\bv)\in \testh\times \testqh$
\begin{equation*}
\aHDG(\uhHDG,\qhHDG;w,\bv) = \ellHDG(w,\bv)    \text{ and }  
\aHDG(\uhHDGA,\qhHDGA;w,\bv)= \ellHDG^O(w,\bv)
\end{equation*}
where $\ellHDG^O(w,\bv)  =
 (\source^O,w)_{\triang }
 -\langle \uD^O,\bv\cdot\bn\rangle _{\GD}
 +\langle \trac^O,w\rangle _{\GN}.$

\item[1.-] Compute the potential and projected equilibrated flux reconstructions 
$\uh, \uhA\in\test$ and $\qhP,\qhAP\in\testq$ such that, $\at{\qhP}{\elem}, \qhAP\at{\!}{\elem}\in RT^p$ verify
\begin{equation*}
\begin{array}{l}
\langle (\qhP-\qhHDGT)\cdot\bn,\mu\rangle_\edge = 0
   \text{ and }  
\langle (\qhAP-\qhHDGTA)\cdot\bn,\mu\rangle_\edge = 0
\quad\forall \mu\in  \mathbb{P}^p(\edge)\quad,\quad
\forall \edge\in\partial\elem
\\[1ex]
(\qhP-\qhHDG,\bv)_\elem = 0
   \text{ and }  
(\qhAP-\qhHDGA,\bv)_\elem = 0
\quad\forall \bv \in  [\mathbb{P}^{p-1}(\elem)]^d \quad,\quad \text{ if } p\geq 1,
\end{array}
\end{equation*}
and $\uh$ and $\uhA$ are continuous averages (exactly verifying the Dirichlet boundary conditions) of $\uhHDG^\ast$ and $\uhHDGA^\ast\in \mathcal{W}_h^{p+1}$ satisfying
\begin{equation*}\label{eq:uhHDGast}
\begin{array}{ll}
(\nabla \uhHDG^\ast,\nabla w)_\elem
=(\nabla \cdot \qhP ,w)_\elem-(\qhP\cdot \bn,w)_{\partial\elem} 
\qquad \forall w\in  \mathcal{W}_h^{p+1},
\\[1ex]
(\nabla \uhHDGA^\ast,\nabla w)_\elem
=(\nabla \cdot \qhAP ,w)_\elem-(\qhAP\cdot \bn,w)_{\partial\elem} 
\qquad \forall w\in  \mathcal{W}_h^{p+1},
\\[1ex]
(\uhHDG^\ast,1)_{\elem} = (\uhHDG,w)_{\elem}
   \text{ and }  
(\uhHDGA^\ast,1)_{\elem} = (\uhHDGA,w)_{\elem}.
\end{array}
\end{equation*}

\item[2.-] For each element of the mesh compute $\kappa = {\norm{\qhAP +\nu\grad\uhA}/\norm{\qhP+\nu\grad\uh}}$ and 
\begin{equation*}
\hspace{-.7cm}
\begin{array}{rl}
\eta_K^{\pi-}
&=
{\norm{\qhAP +\nu\grad \uhA -\kappa (\qhP+\nu\grad \uh)}}_K 
	+
	C_1 \nu_K^{-1/2}
	||\source^O-\projk \source^O -\kappa(\source-\projk \source)||_{\mathcal{L}^2(K)}	
\\[1ex]
&\displaystyle	
	+ \sum\limits_{\edge\in \GN \cap \partial\elem}
	\CTb \nu_K^{-1/2}||\trac^O-\projg\trac^O+\kappa (\trac-\projg \trac) ||_{\mathcal{L}^2(\edge)}
\\[3ex]	
\eta_K^{\pi+}
&=
{\norm{\qhAP +\nu\grad \uhA+\kappa (\qhP+\nu\grad \uh)}}_K 
	+
	C_1 \nu_K^{-1/2}
	||\source^O-\projk \source^O +\kappa(\source-\projk \source)||_{\mathcal{L}^2(K)}	
\\[1ex]
&\displaystyle	
	+ \sum\limits_{\edge\in \GN \cap \partial\elem}
	\CTb \nu_K^{-1/2}||\trac^O-\projg\trac^O-\kappa (\trac-\projg \trac) ||_{\mathcal{L}^2(\edge)}		
\end{array}	
\end{equation*}

\item[3.-]  Compute the approximation of $s$ 
\[
\sh =(\source^O,\uh)
+\langle \trac^O,\uh\rangle _{\GN}
+ (\source,\uhA)
-\langle \trac,\uhA\rangle _{\GN}
-(\nu\grad  \uh,\grad  \uhA)
\]
and the bounds for the quantity of interest $\sub$ and $\slb$
\begin{equation*}
\slb 
=
\sh
-\dfrac{1}{4\kappa}\sum\limits_{\elem\in\triang} (\eta_K^{\pi-})^2
\quad\text{ and } \quad 
\sub 
=
\sh
+\dfrac{1}{4\kappa}\sum\limits_{\elem\in\triang} (\eta_K^{\pi+})^2.
\end{equation*}

\end{itemize}
\end{minipage}
}\caption{Bounds for the quantity of interest from the HDG approximations}\label{box:bounds}
\end{figure}

\section{Numerical examples}
\label{sec:numerical_examples}

The behavior of the bounding procedure described above is analyzed in two numerical examples. 
Four estimates of $s$ are considered: the upper and lower bounds ($\sub$ and $\slb$ respectively), their average $\sh = (\sub + \slb)/2$ and the quantity of interest given by the HDG finite element approximation, denoted by $s_h=\ell^O(\uhHDG,\qhHDG)$. {The stabilization parameter is set to $\tau=1$ in all the cases.}

A measure of the accuracy of the bounds is the half bound gap $\Delta_h/2 = (\sub - \slb)/2$ since it is an upper bound for the error between the approximation $\sh$ and the exact output, see equation \eqref{eq:bounds_QoI_01_mod}. 
The bound gap also provides local error information which {is} used as an indicator for mesh adaptivity. Indeed, the bound gap associated {with} the bounding strategy described in Figure \ref{box:bounds} {is split} using the local elemental contributions 
\[
\Delta_h = \sub - \slb
=
\dfrac{1}{4\kappa}\sum\limits_{\elem\in\triang} ((\eta_K^{\pi+})^2+(\eta_K^{\pi-})^2)
=
\sum\limits_{\elem\in\triang} \Delta_h^K.
\]
The elemental contributions $\Delta_h^K$ {are} informative mesh adaptivity indicators for controlling the error in the quantity of interest. Note that these indicators take into account the error in both the primal and adjoint problems and also the data oscillation errors, and therefore, the mesh is refined both in the areas most contributing to the error and in the areas where the data cannot be properly represented using its projection.

Two remeshing strategies are considered, see for instance \cite{MR2186143,MR2914423}. In the first strategy, 
given a target bound gap $\Delta_\text{tol}$, a uniform error distribution assumption is used and, at each level of refinement, the elements with $\Delta_h^K\geq(\Delta_\text{tol})/\numel$ are refined where $\numel=|\triang|$ denotes the number of triangles of the mesh. The second strategy refines the elements according to a bulk criterion, that is, given a prescribed scalar parameter $\Theta\in(0,1]$, selects a subset $M$ of $\triang$ such that 
$\Theta (\sum_{\elem\in\triang} \Delta_h^K) \leq \sum_{\elem\in M} \Delta_h^K$.

Finally, in the numerical experiments we monitor the convergence of the estimates via the computational order of convergence calculated as follows. {We denote by} $e(\numel)$ and $e(\tilde{\numel})$ any error-like quantity for two consecutive triangulations with $\numel$ and $\tilde{\numel}$ number of triangles. Then, the computational ratio of convergence is given by
\[
-2 \dfrac{\log(e(\numel)/e(\tilde{\numel}))}{\log(\numel/\tilde{\numel})}.
\]

\subsection{Example 1 - Smooth solution}
First, we investigate the order of convergence of the bounds for smooth solutions. 
Consider the Poisson equation in the square plate $\Omega=(0,1)^2$
with homogeneous Dirichlet boundary conditions and empty Neumann boundary, namely
$\partial\Omega=\GD, \nu=1$ and $\uD=0$ in equation \eqref{eq:strong}. The right-hand side $\source$ is chosen such that the exact solution is given by
\[
u(x,y) = \sin(\pi x)\sin(\pi y).
\]
Two quantities of interest are considered. The first one, $s_1$,  is an average of the solution over the whole domain, and the second one, $s_2$, is a weighted average of the normal flux in the Dirichlet boundary. These quantities of interest are given by equation \eqref{eq:QoI} for 
\begin{itemize}
\item Data for $s_1$: $\uD^O=\trac^O=0$ and $\source^O(x,y) = 1$, where $s_1=4/\pi^2$ for
\[
\xi(x,y)=\dfrac{x(1-x)}{2}-\dfrac{4}{\pi^3} \sum\limits_{\text{odd } k} 
\dfrac{\sin(k\pi x)(\sinh(k\pi y)+\sinh(k\pi (1-y)))}{k^3\sinh(k\pi)}.
\]
\item Data for $s_2$: $\source^O=\trac^O=0$ and $\uD^O=\frac{\pi}{2}\sin(\pi y)$ on $x=1$ and $\uD^O=0$ elsewhere, where $\xi(x,y) = \pi \sin(\pi y)\sinh(\pi x)/(2\sinh(\pi))$ and $s_2 = \pi^2/4$.
\end{itemize}
The HDG approximations of both the primal and adjoint problems associated {with} $s_2$ have an optimal convergence and both ${\norm{\qh+\nu\grad  \uh}}$ and ${\norm{\qhA+\nu\grad  \uhA}}$ superconverge with order $p+1$, and therefore the bound gap is expected to converge with order $2(p+1)$ or $\mathcal{O}(\numel^{-(p+1)})$. However, the adjoint solution associated {with} $s_1$ verifies $\xi\in\mathcal{H}^3(\Omega)$ and therefore we expect that ${\norm{\qhA+\nu\grad  \uhA}}$ converges with order $2$ for $p\geq 1$, see \cite{MR1813251,MR3725770}, yielding an expected convergence of the bound gap of order $p+3$ or $\mathcal{O}(\numel^{-(p+3)/2})$.

\begin{figure}[ht!]
\psfrag{a}[cc][cc][0.85]{Number of elements}
\psfrag{b}[cc][cc][0.85]{Half bound gap}
\psfrag{c}[Bl][Bl][0.85]{$p=1$}
\psfrag{d}[Bl][Bl][0.85]{$p=2$}
\psfrag{e}[Bl][Bl][0.85]{$p=3$}
\psfrag{f}[Bl][Bl][0.85]{$p=4$}
\psfrag{g}[Bl][Bl][0.85]{$\numel^{-2}$}
\psfrag{h}[Bl][Bl][.85]{$\numel^{-5/2}$}
\psfrag{i}[Bl][Bl][0.85]{$\numel^{-3}$}
\psfrag{j}[Bl][Bl][0.85]{$\numel^{-7/2}$}
\psfrag{k}[cc][cc][0.85]{$10^2$}
\psfrag{m}[cc][cc][0.85]{$10^3$}
\psfrag{n}[cc][cc][0.85]{$10^4$}
\psfrag{o}[cc][cc][0.85]{$10^{-4}$}
\psfrag{p}[cc][cc][0.85]{$10^{-6}$}
\psfrag{q}[cc][cc][0.85]{$10^{-8}$}
\centering
\includegraphics[width=.8\textwidth]{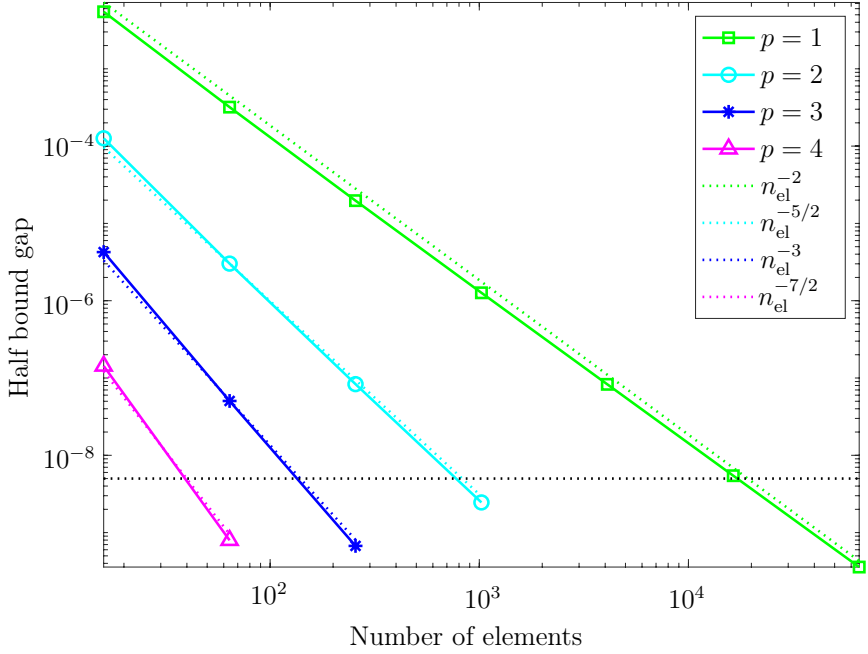}
\caption{Example 1: $s_1$  - Convergence of the half bound gap for a uniform mesh refinement (optimal convergence of order $p+3$ or $\mathcal{O}(\numel^{-(p+3)/2})$).}
\label{fig:SMOOTH_Uniforced_halfgap}
\end{figure}

The numerical results for the first quantity of interest $s_1$ for a uniform mesh refinement are shown in Table \ref{table:SMOOTH_Uniforced} and Figure \ref{fig:SMOOTH_Uniforced_halfgap}, where the stopping criteria is set to achieve $\Delta_h<10^{-8}$. The initial structured mesh consists of $16$ triangles and at each refinement, every triangle is divided into four similar triangles. For all the values of $p$, the optimal order of convergence $p+3$ predicted by the theory is achieved. It is also worth noting that for high order polynomials, the required precision is achieved with very coarse meshes.
Since the number of degrees of freedom of the global system of HDG computations is $n_\text{edge} \approx 3(p+1)\numel/2$ and taking into account that the manipulation of the mesh takes up a significant amount of computational effort both in the HDG computation and in the a posteriori error estimation procedure, working with high-order polynomials seems to be advantageous. For this problem and with this particular quantity of interest, using adaptive mesh refinement strategies does not provide significantly more accurate bounds, since the error is uniformly distributed both for the primal and adjoint problems.
\begin{table}[hbt!]
\centering
\begin{tabular}{cc|cc|cc} 
\hline
$\numel$ & $n_\text{edge}$ & $\sh \pm \Delta_h/2$  & order & $|s-s_h|$ & $|s-\sh|$    \\
\hline
\multicolumn{5}{c}{$p=1$}\\
\hline
      16 &       56 & 0.406021554922 $\pm$ 5.47e-03 & -- & 1.90e-03 & 7.37e-04 \\ 
      64 &      208 & 0.405317075580 $\pm$ 3.19e-04 & 4.10 & 3.64e-04 & 3.23e-05 \\ 
     256 &      800 & 0.405286596697 $\pm$ 1.97e-05 & 4.02 & 5.01e-05 & 1.86e-06 \\ 
    1024 &     3136 & 0.405284843586 $\pm$ 1.27e-06 & 3.96 & 6.52e-06 & 1.09e-07 \\ 
    4096 &    12416 & 0.405284741107 $\pm$ 8.28e-08 & 3.93 & 8.32e-07 & 6.54e-09 \\ 
   16384 &    49408 & 0.405284734968 $\pm$ 5.45e-09 & 3.93 & 1.05e-07 & 3.99e-10 \\ 
   65536 &   197120 & 0.405284734592 $\pm$ 3.58e-10 & 3.93 & 1.20e-08 & 2.24e-11 \\ 
\hline    
\multicolumn{5}{c}{$p=2$}\\
\hline
      16 &       84 & 0.405275669432 $\pm$ 1.26e-04 & -- & 6.64e-05 & 9.07e-06 \\ 
      64 &      312 & 0.405284783569 $\pm$ 3.02e-06 & 5.38 & 1.10e-06 & 4.90e-08 \\ 
     256 &     1200 & 0.405284735937 $\pm$ 8.33e-08 & 5.18 & 2.14e-08 & 1.37e-09 \\ 
    1024 &     4704 & 0.405284734592 $\pm$ 2.46e-09 & 5.08 & 4.86e-10 & 2.26e-11 \\ 
\hline   
\multicolumn{5}{c}{$p=3$}\\
\hline
      16 &      112 & 0.405284626142 $\pm$ 4.25e-06 & -- & 8.77e-08 & 1.08e-07 \\ 
      64 &      416 & 0.405284735218 $\pm$ 5.04e-08 & 6.40 & 7.91e-09 & 1.05e-08 \\ 
     256 &     1600 & 0.405284734574 $\pm$ 6.73e-10 & 6.23 & 3.37e-11 & 4.53e-12 \\ 
\hline    
\multicolumn{5}{c}{$p=4$}\\
\hline
      16 &      140 & 0.405284734710 $\pm$ 1.43e-07 & -- & 4.17e-08 & 1.41e-10 \\ 
      64 &      520 & 0.405284734520 $\pm$ 7.95e-10 & 7.49 & 1.71e-10 & 4.96e-11 \\ 
\hline
\end{tabular}
\caption{Example 1: $s_1$  - Uniform mesh refinement: effect of the polynomial degree $p$.}
\label{table:SMOOTH_Uniforced}
\end{table}

To compute the bounds for the second quantity of interest, it is worth noting that in this case, a simple averaging of the post-processed HDG approximation $\uhHDGA^\ast\in \mathcal{W}_h^{p+1}$ does not yield a potential reconstruction since it does not exactly verify the Dirichlet boundary conditions $\uhA = \uD^O=\frac{\pi}{2}\sin(\pi y)$ on the right edge ($x=1$). In this case, even though more elaborate extensions operators could be used, see for instance \cite{Vejchodsk2004,MR2513831}, 
since the bounding procedure is valid for any potential reconstruction $\uhA$, the exact Dirichlet boundary conditions are enforced via an easy modification of $\uhA$ in a small band around $x=1$. Specifically, $\uhA$ is obtained as
\begin{itemize}
\item[1.] the post-processed scalar variable $\uhHDGA^\ast$ is averaged to obtain a continuous reconstruction $\uhA$
\item[2.] the maximum value $x_{band}\in[0,1)$ such that the straight line $x=x_{band}$ does not intersect any element interior is computed

\item[3.] introducing the following extension of the Dirichlet boundary conditions 
\[
\uDe^O = \frac{\pi}{2}\sin(\pi y)\dfrac{x-x_{band}}{1-x_{band}} \text{ for }
x\in [x_{band},1] \text{ and } \uDe^O = 0 \text{ otherwise},
\]
and its global nodal interpolant 
$\mathcal{I}^{p+1}_h\left(\uDe^O\right)$,
the value of $\uhA$ is modified on the edge $x=1$, 
$\uhA|_{\{x=1\}} = \mathcal{I}^{p+1}_h\left(\uDe^O\right)|_{\{x=1\}}$

\item[4.] for each element inside the band $[x_{band},1]\times[0,1]$, {the final value of $\uhA|_{\elem}$ is set adding the interpolation error}
\[
\uhA|_{\elem} + \uDe^O|_{\elem}  -\mathcal{I}^{p+1}_h\left(\uDe^O\right)|_{\elem}.
\]
\end{itemize}

\begin{figure}[ht!]
\centering
\begin{minipage}[c]{.35\textwidth}
\psfrag{a}[cc][cc][0.85]{$x=x_{band}$}
\includegraphics[width=\textwidth]{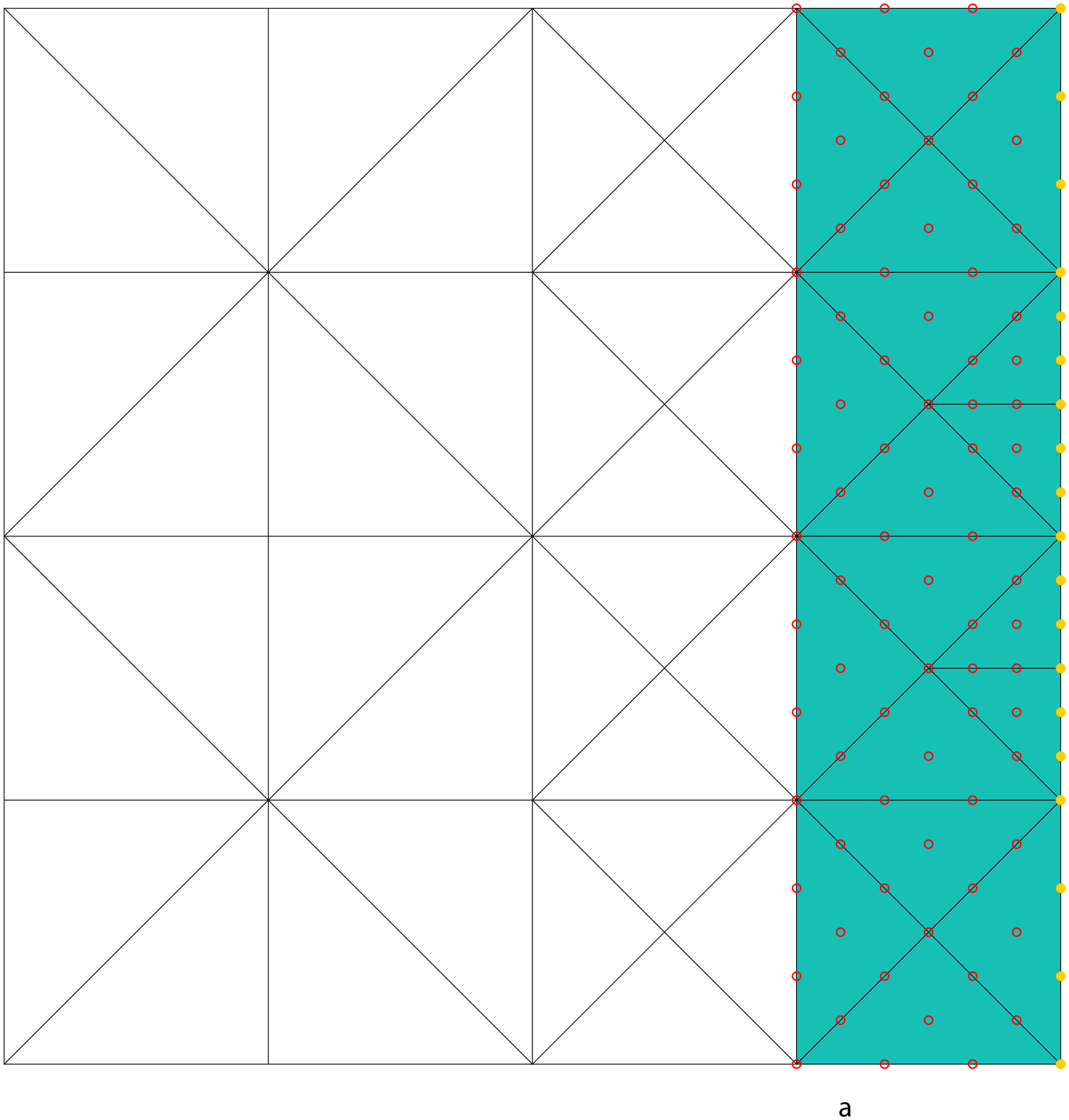}
\end{minipage}
\begin{minipage}[c]{.45\textwidth}
\includegraphics[width=\textwidth]{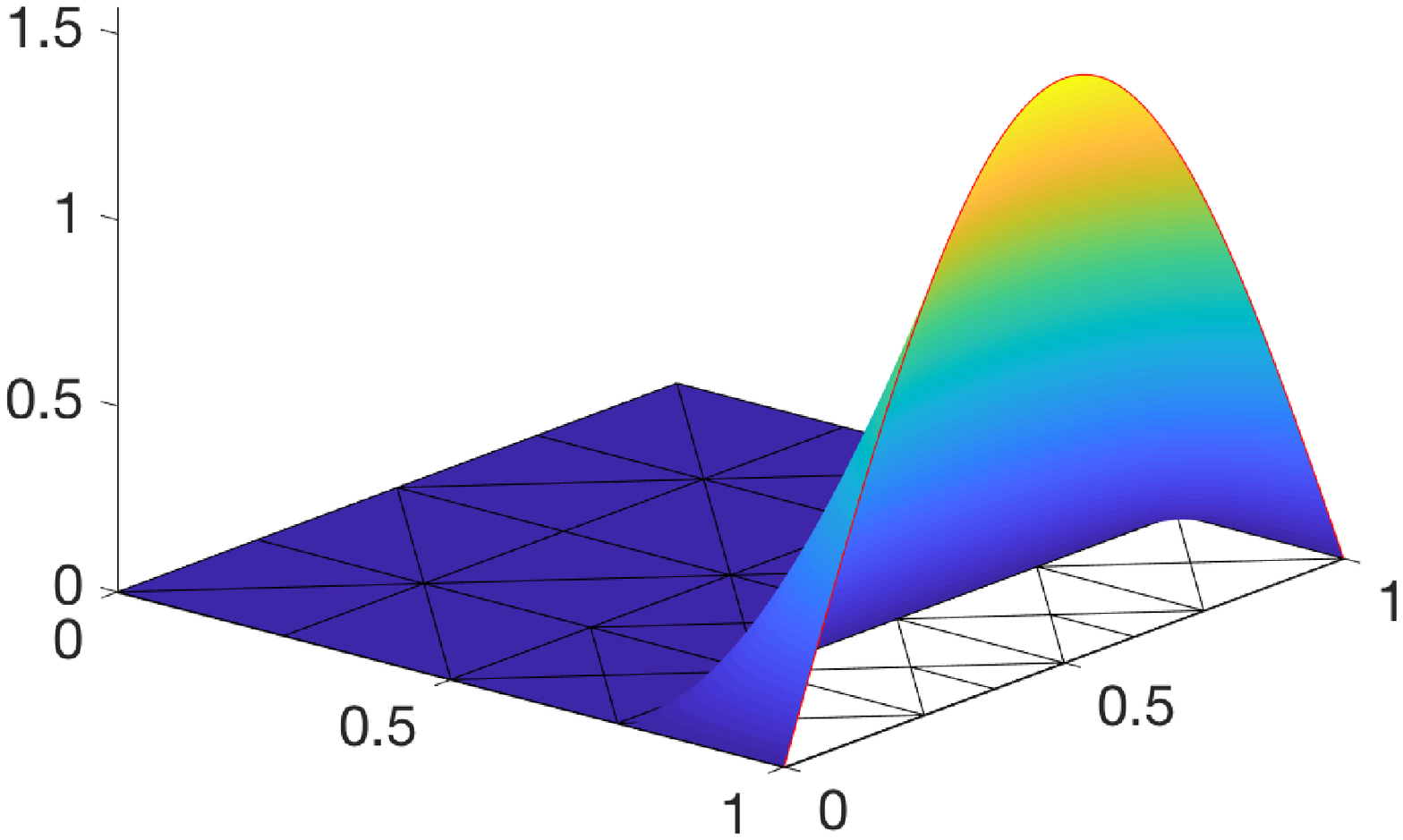}
\end{minipage}
\caption{Example 1: $s_2$  - Exact enforcement of the Dirichlet boundary conditions: band for $p=3$ (left) and 
$\uDe^O$  (right).
}
\label{fig:SMOOTH_CCD_exact_ccd1}
\end{figure}

\begin{figure}[ht!]
\centering
\includegraphics[width=.32\textwidth]{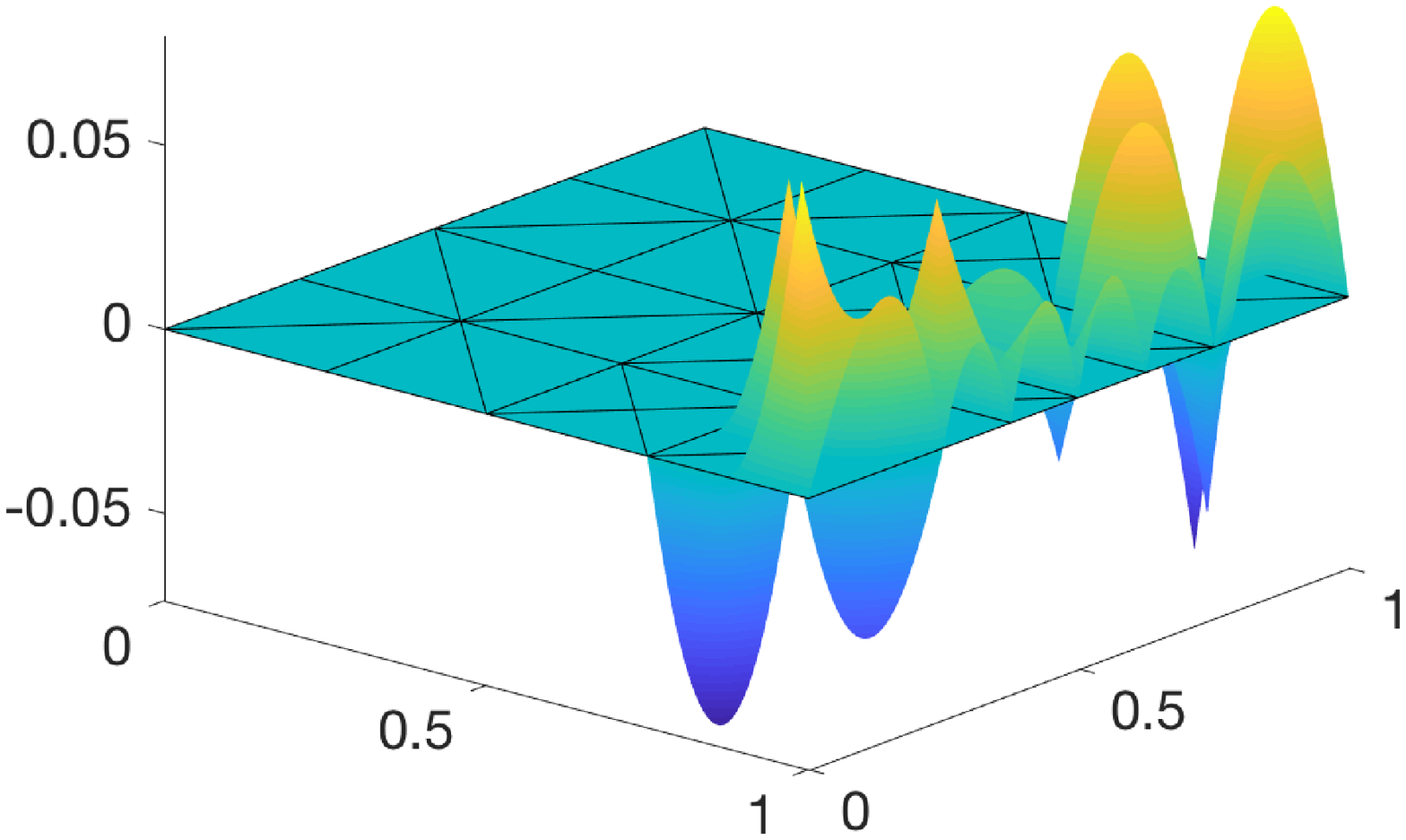}
\includegraphics[width=.32\textwidth]{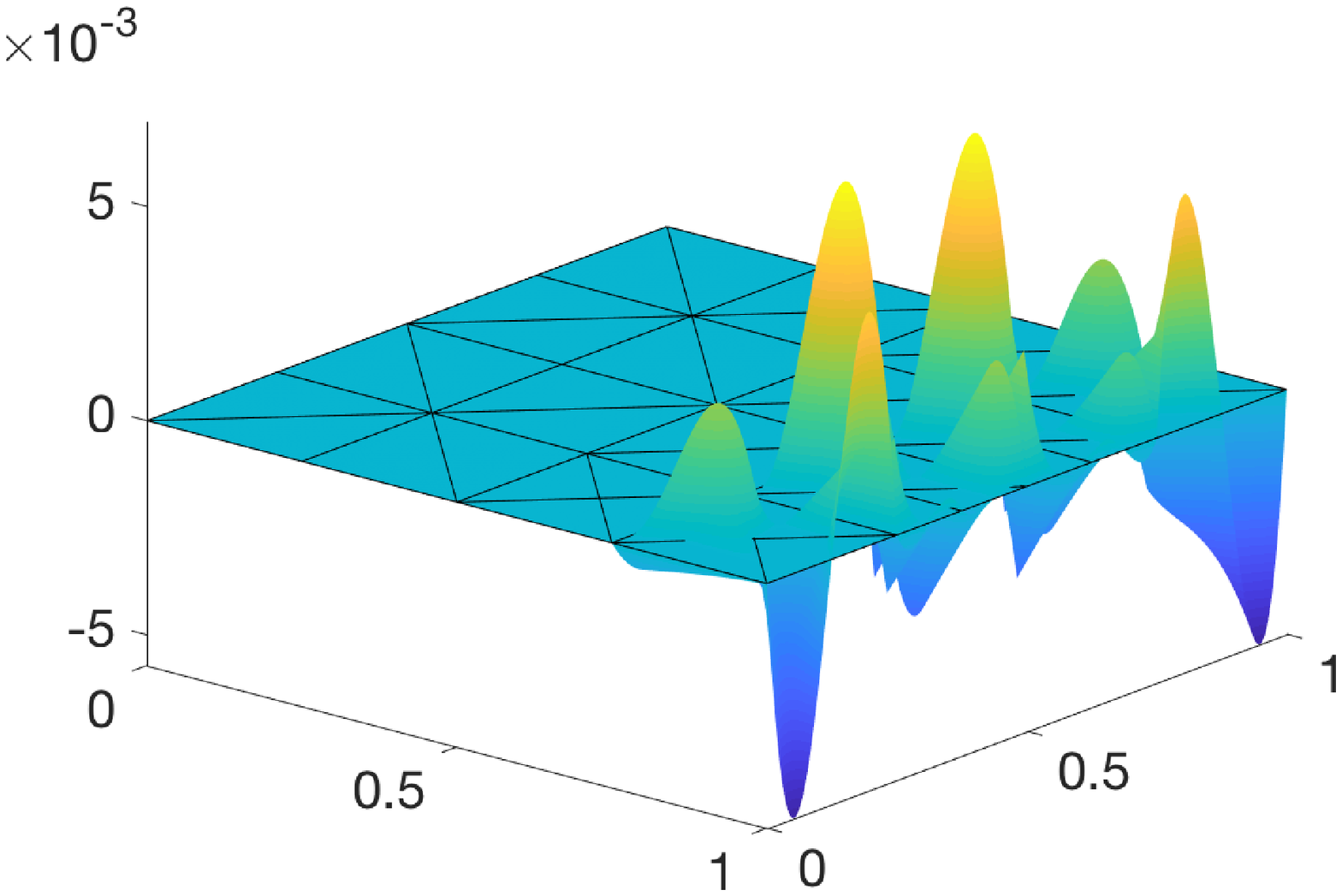}
\includegraphics[width=.32\textwidth]{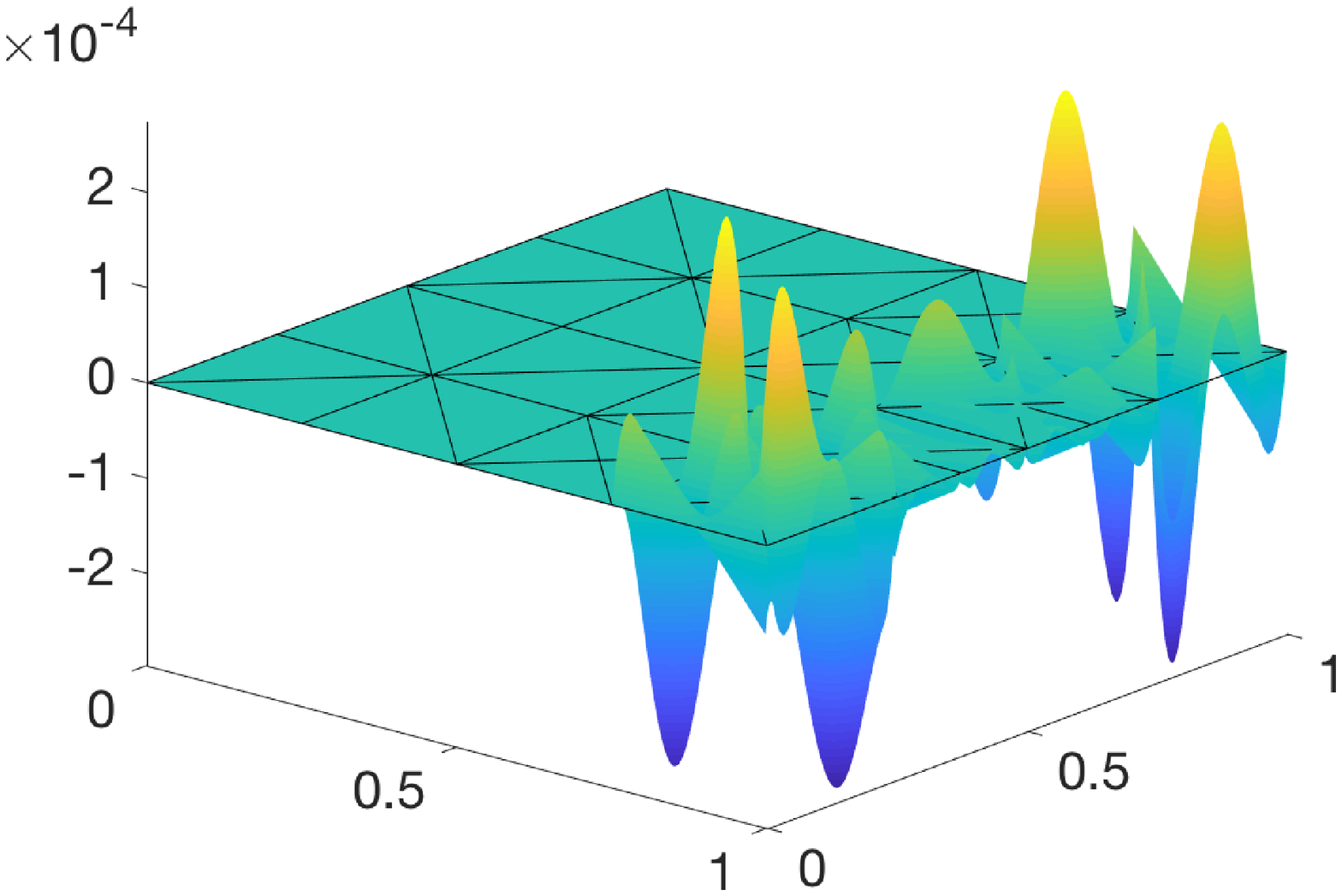}
\caption{Example 1: $s_2$  - {Plots of $\uDe^O-\mathcal{I}^{p+1}_h\left(\uDe^O\right)$
for $p=1,2,3$ (exact enforcement of the Dirichlet boundary conditions)}.
}
\label{fig:SMOOTH_CCD_exact_ccd2}
\end{figure}

Figure \ref{fig:SMOOTH_CCD_exact_ccd1} shows the band where the solution is modified and the shape of $\uDe^O$ for a particular mesh while Figure  \ref{fig:SMOOTH_CCD_exact_ccd2} shows the magnitude of the modifications given by the functions $\uDe^O-\mathcal{I}^{p+1}_h\left(\uDe^O\right)$.
It can be seen that the proposed procedure only introduces relevant modifications to the adjoint approximation for small values of $p$ and coarse meshes. In these cases, more involved strategies could be considered if no adaptive procedures alleviating the influence of the boundary conditions are available.

\begin{figure}[ht!]
\centering
\psfrag{a}[cc][cc][0.85]{Number of elements}
\psfrag{b}[cc][cc][0.85]{Half bound gap}
\psfrag{c}[Bl][Bl][0.75]{uniform $p=1$}
\psfrag{d}[Bl][Bl][0.75]{uniform $p=2$}
\psfrag{e}[Bl][Bl][0.75]{uniform $p=3$}
\psfrag{f}[Bl][Bl][0.75]{uniform $p=4$}
\psfrag{g}[Bl][Bl][0.75]{$\numel^{-2}$}
\psfrag{h}[Bl][Bl][.85]{$\numel^{-3}$}
\psfrag{i}[Bl][Bl][0.75]{$\numel^{-4}$}
\psfrag{j}[Bl][Bl][0.75]{$\numel^{-5}$}
\psfrag{k}[Bl][Bl][0.75]{adapt $p=1$ (opt)}
\psfrag{m}[Bl][Bl][0.75]{adapt $p=2$ (opt)}
\psfrag{n}[Bl][Bl][0.75]{adapt $p=3$ (opt)}
\psfrag{o}[Bl][Bl][0.75]{adapt $p=4$ (opt)}
\psfrag{p}[Bl][Bl][0.75]{adapt $p=1$}
\psfrag{q}[Bl][Bl][0.75]{adapt $p=2$}
\psfrag{r}[Bl][Bl][0.75]{adapt $p=3$}
\psfrag{s}[Bl][Bl][0.75]{adapt $p=4$}
\includegraphics[width=.7\textwidth]{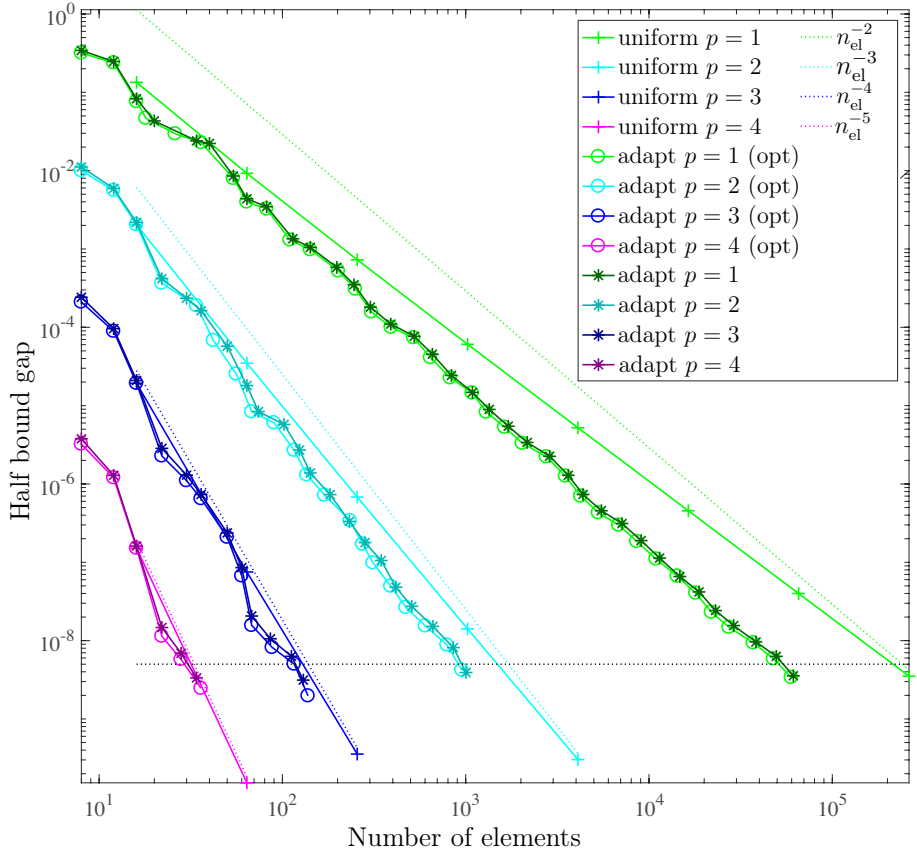}
\caption{Example 1: $s_2$  - Convergence of the half bound gap both a uniform and adaptive mesh refinements.}
\label{fig:SMOOTH_CCD_halfgap}
\end{figure}

\begin{figure}[ht!]
\begin{center}
\includegraphics[width=.24\textwidth]{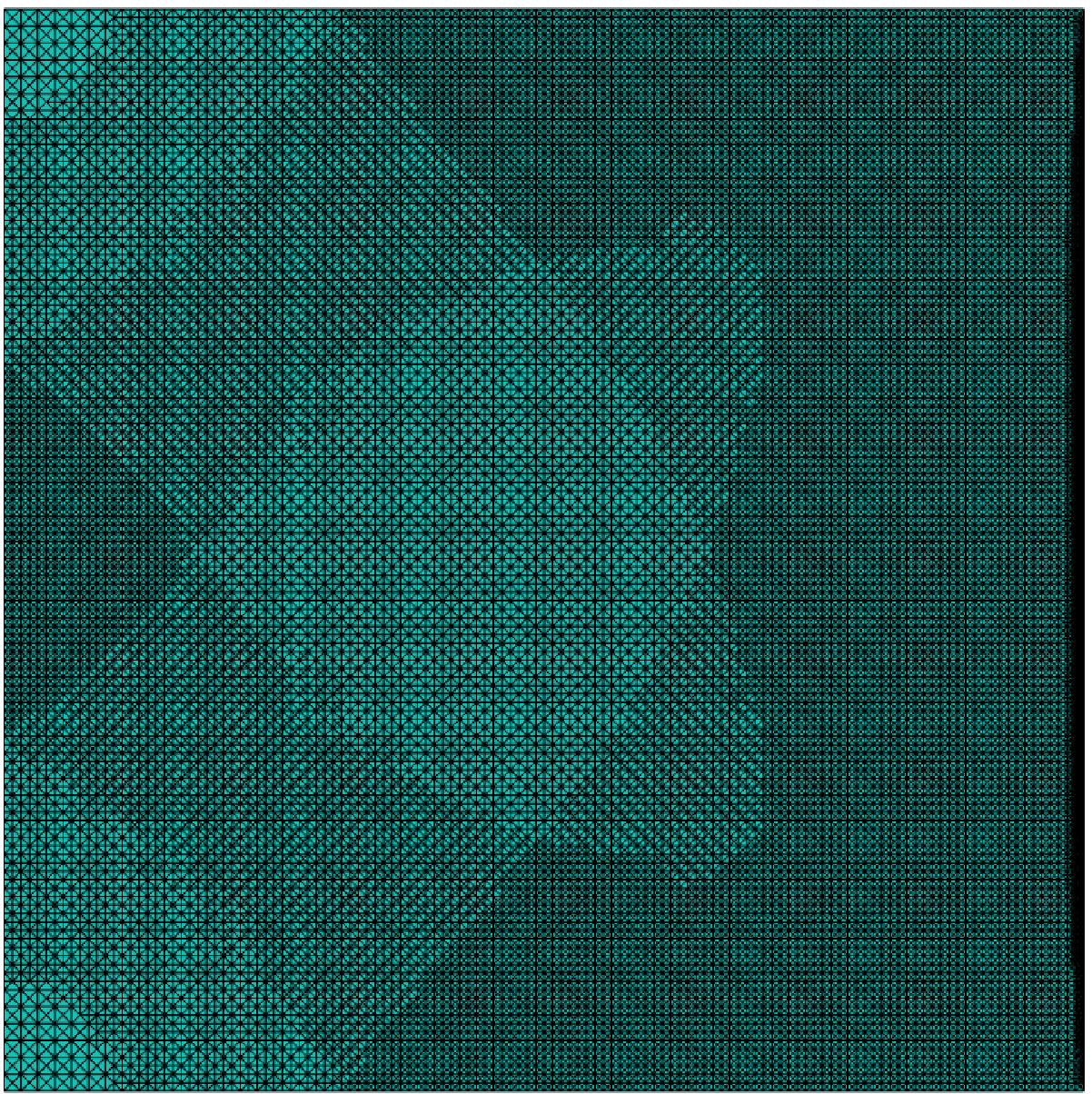}
\includegraphics[width=.24\textwidth]{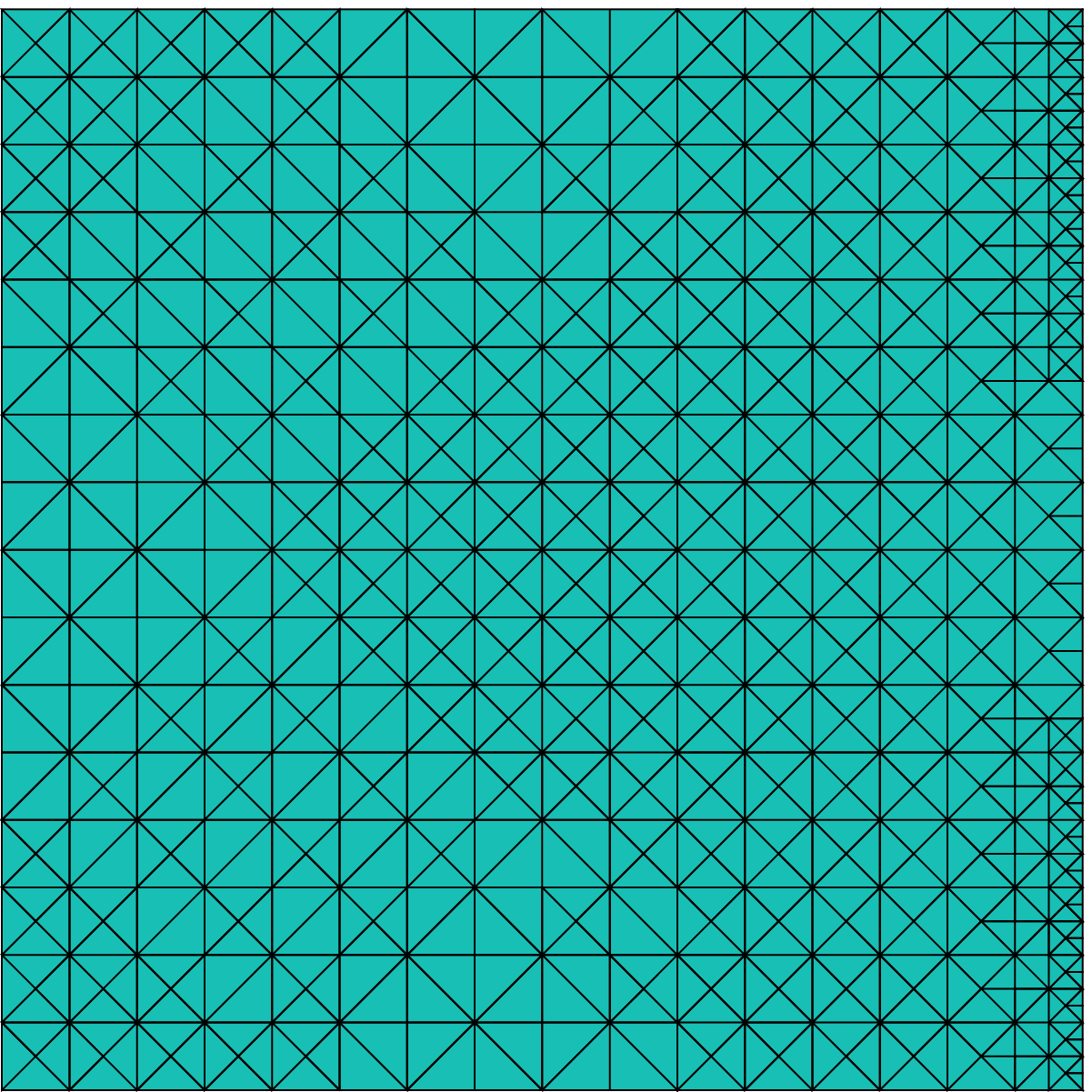}
\includegraphics[width=.24\textwidth]{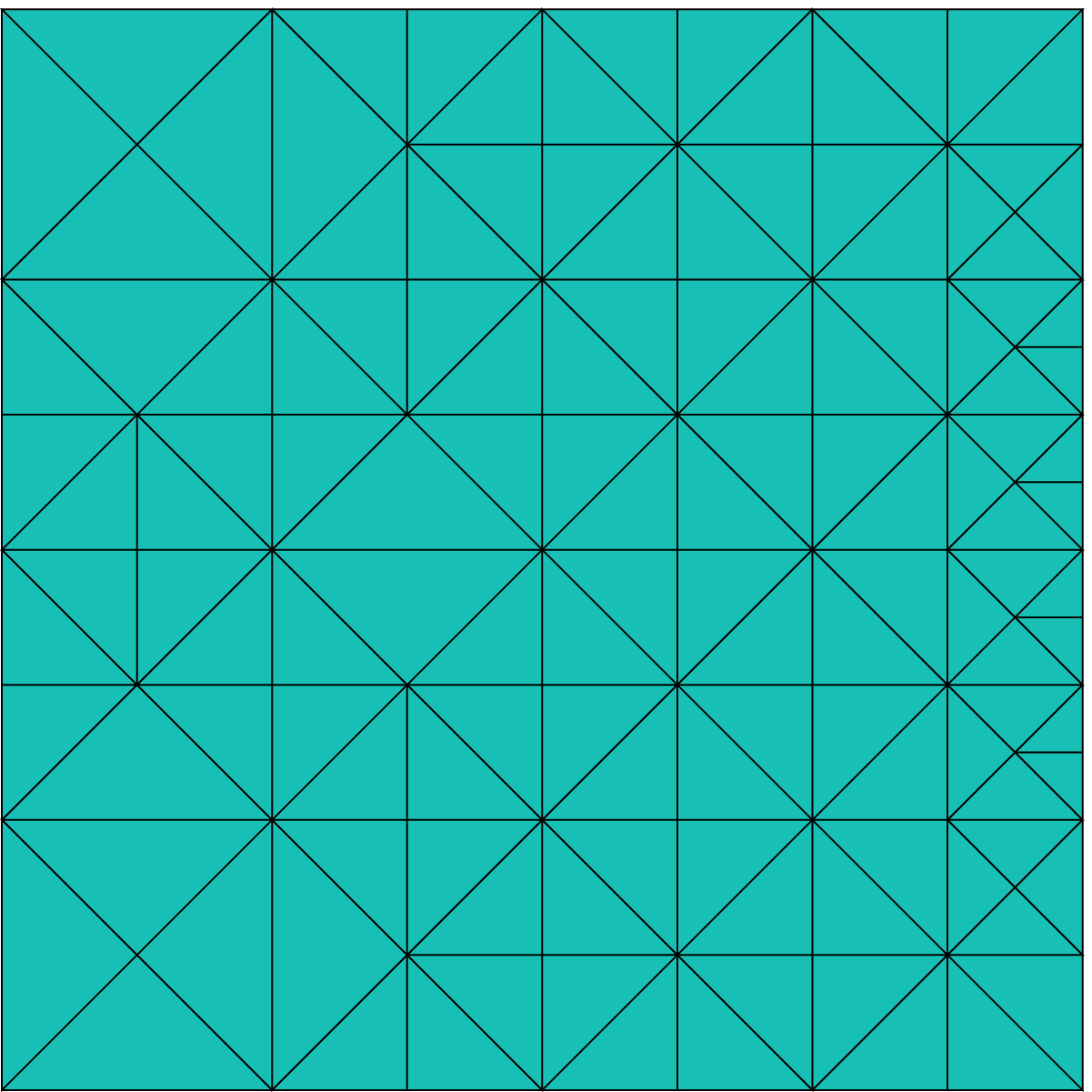}
\includegraphics[width=.24\textwidth]{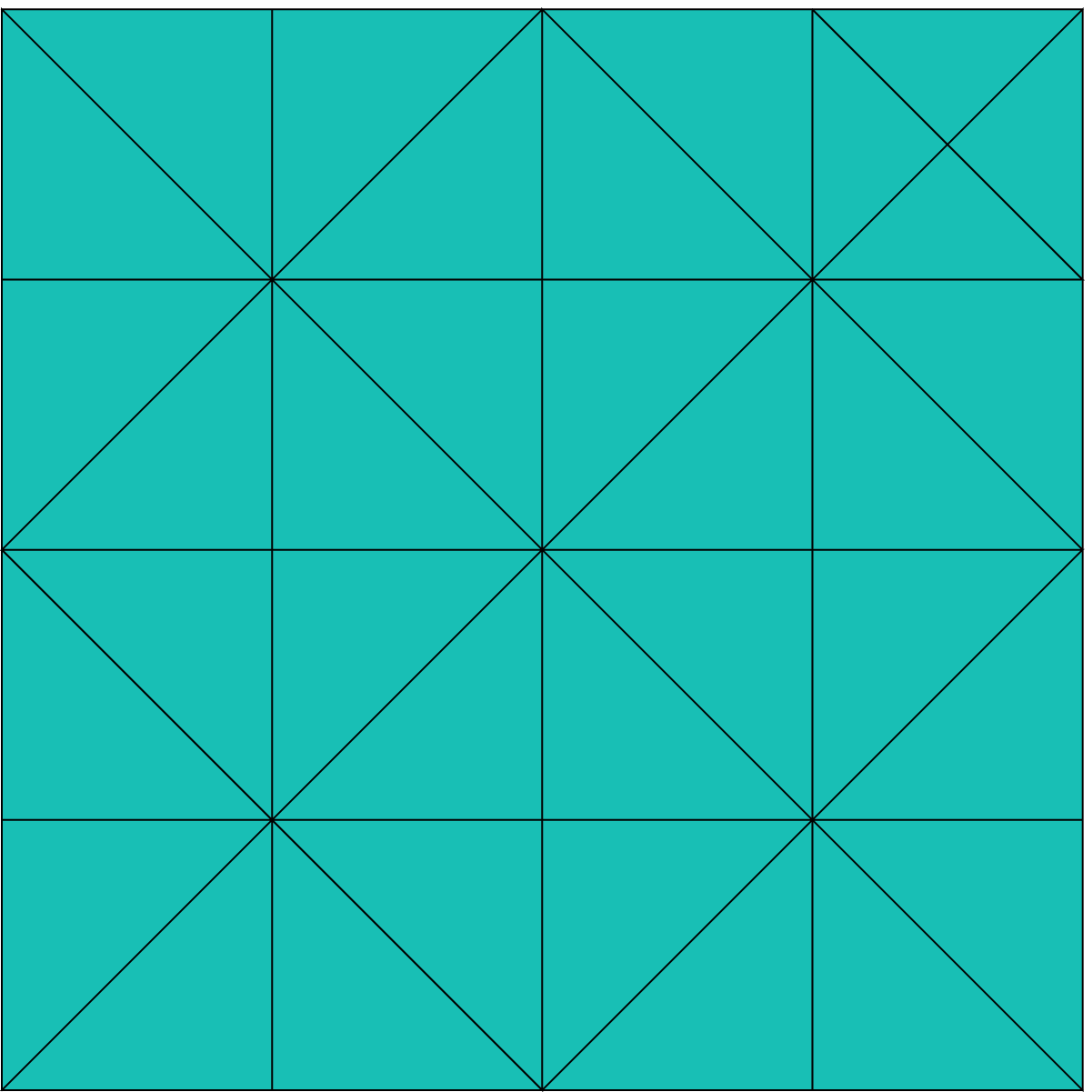}\\
\includegraphics[width=.24\textwidth]{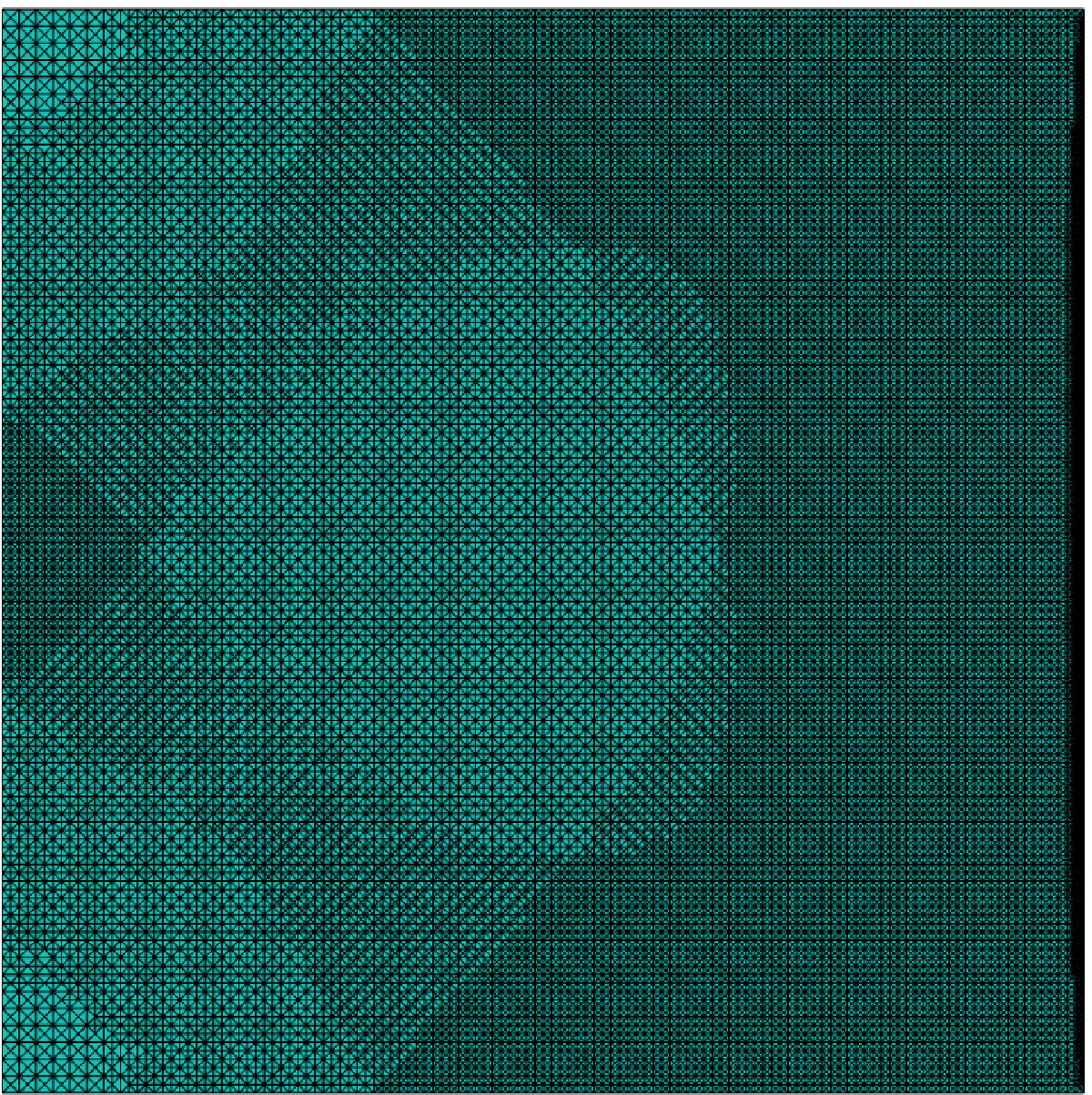}
\includegraphics[width=.24\textwidth]{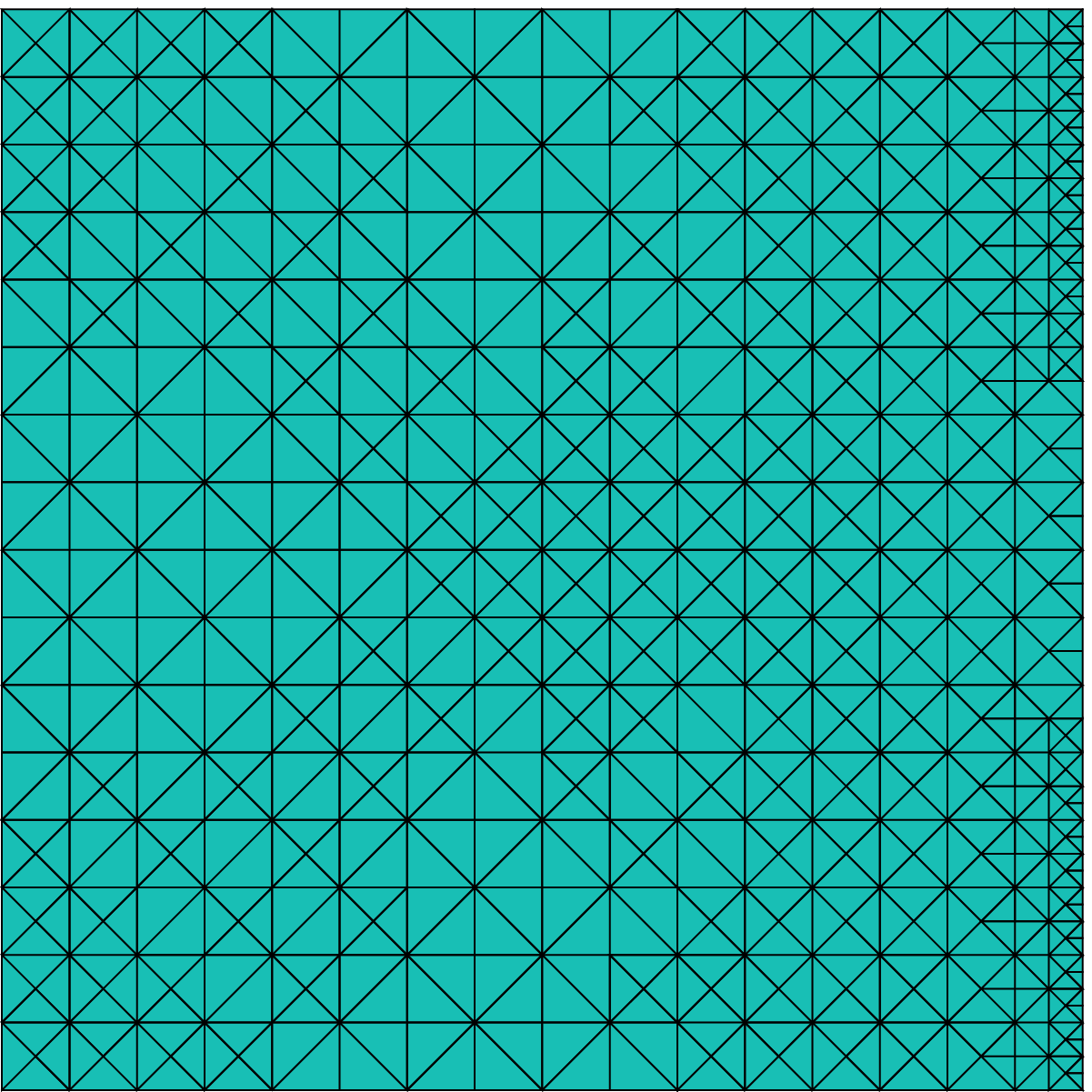}
\includegraphics[width=.24\textwidth]{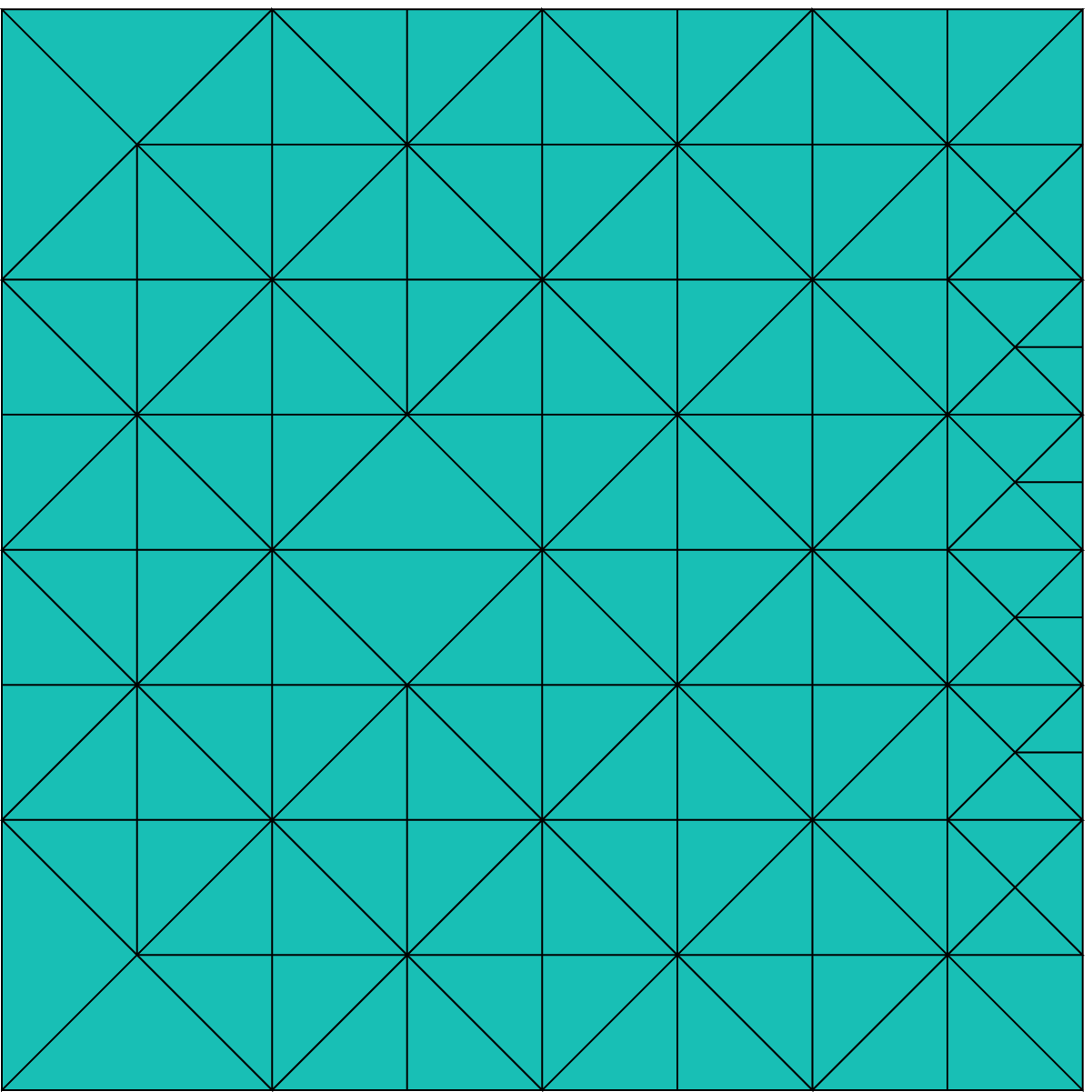}
\includegraphics[width=.24\textwidth]{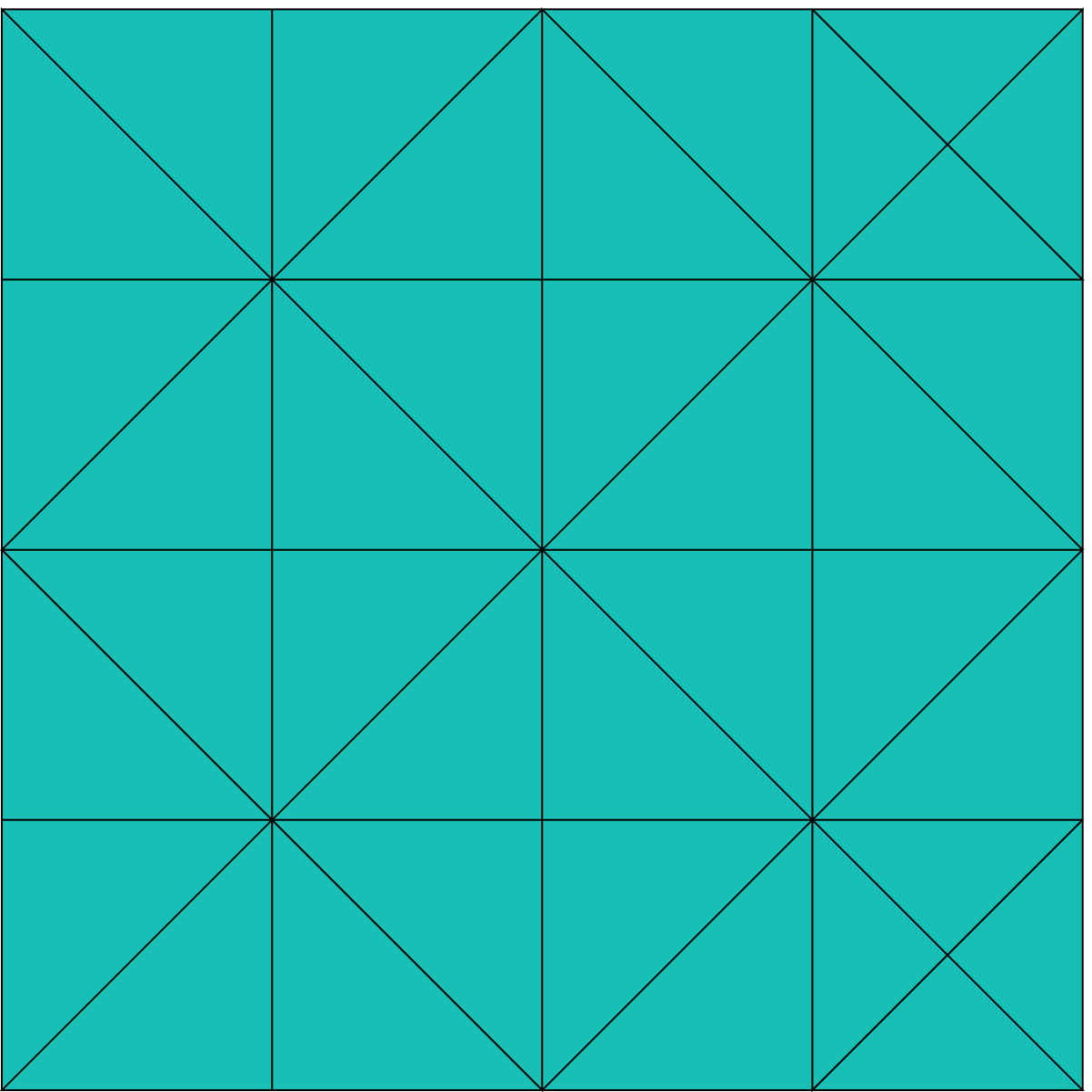}\\
\end{center}
\caption{Example 1: $s_2$  - Final meshes of the adaptive procedure for $p=1, 2, 3$ and $4$ from left to right with $\numel = 61310, 1004, 130$ and $34$ respectively (top) and final meshes with extra optimization procedure with $\numel = 59762, 952, 138$ and $36$ respectively (bottom).}
\label{fig:SMOOTH_CCD__final_meshes} 
\end{figure}

Figure \ref{fig:SMOOTH_CCD_halfgap} shows the convergence of the half bound gap for the first quantity of interest $s_1$  and for a final tolerance limit $\Delta_h = 10^{-8}$. The convergence is shown both for a uniform mesh refinement and the adaptive strategy following the bulk criterion for $\Theta=0.5$. For the adaptive procedure, both the bounds associated {with} the reconstructions shown in Sections \ref{subsec:pot_rec} and \ref{subsec:flux_rec} and the bounds obtained adding the extra local optimization procedure detailed in Subsection \ref{subsec:localopt} are shown. As can be seen both in this figure and in Table \ref{table:SMOOTH_CCD_finalmesh}, the extra local optimization procedure provides an improvement of the value for the half bound gap that becomes more relevant as $p$ increases. 
\begin{table}[hbt!]
\centering
\begin{tabular}{cc|c|c} 
\hline
& $\numel$ & $\sh \pm \Delta_h/2$  & $|s-\sh|$    \\
\hline
\multicolumn{4}{c}{$p=1$}\\
\hline
&   61310 & 2.467401099996 $\pm$ 3.56e-09 & 2.76e-10 \\ 
optimized &   59762 & 2.467401100185 $\pm$ 3.40e-09 & 8.71e-11 \\    
\hline
\multicolumn{4}{c}{$p=2$}\\
\hline   
&    1004 & 2.467401100039 $\pm$ 3.92e-09 & 2.33e-10 \\ 
optimized &     952 & 2.467401100554 $\pm$ 4.22e-09 & 2.81e-10 \\ 
\hline
\multicolumn{4}{c}{$p=3$}\\
\hline     
&     130 & 2.467401100099 $\pm$ 3.13e-09 & 1.73e-10 \\ 
optimized &     138 & 2.467401100343 $\pm$ 1.98e-09 & 7.02e-11 \\ 
\hline
\multicolumn{4}{c}{$p=4$}\\
\hline     
&      34 & 2.467401100022 $\pm$ 3.34e-09 & 2.50e-10 \\ 
optimized &      36 & 2.467401100173 $\pm$ 2.46e-09 & 9.88e-11 \\ \hline
\end{tabular}
\caption{Example 1: $s_2$  - Bounds for the final meshes of the adaptive procedure.}
\label{table:SMOOTH_CCD_finalmesh}
\end{table}
Also note that for $p=1$ and $p=2$ optimal convergence is only reached when adaptive procedures are used, due to the simple procedure used to exactly impose the Dirichlet boundary conditions in the adjoint problem. If no adaptive procedures are available, more involved techniques could be used to achieve optimal convergence, see \cite{Vejchodsk2004,MR2513831}.
Finally, Figure \ref{fig:SMOOTH_CCD__final_meshes} shows the final meshes obtained in the adaptive procedures. As can be seen, using the extra local optimization procedure does not significantly introduce changes in the final meshes while providing slightly better results with a small extra computational cost.

\subsection{Example 2 - L-shaped domain}

Consider the Poisson equation, $\nu=1$, in the L-shaped domain $\Omega=[-1,1]^2\backslash (0,1)\times (-1,0)$ with right-hand side $\source=1$ and all  homogeneous Dirichlet boundary conditions, that is, $\partial\Omega=\GD$ and $\uD=0$. The exact solution is unknown, but it's energy norm is ${\norm{\nu \grad u}}^2 = 0.2140758036140825$, see \cite{MR2667298}. The solution has a typical corner singularity at the origin and a theoretical convergence rate of the error in the energy norm is $\mathcal{O}(h^{2/3})=\mathcal{O}(\numel^{-1/3})$.

Two quantities of interest are considered. The first quantity of interest is associated {with} $\uD^O=\trac^O=0$ and $\source^O(x,y) = \source(x,y) = 1$. In this case, the primal and adjoint problems coincide yielding $s_1={\norm{\nu \grad u}}^2$. The second quantity of interest, $s_2$ is taken from \cite{MR2470943,MR2899560} and is associated {with} the data $\uD^O=\trac^O=0$ and 
\[
\source^O(x,y) = -\dfrac{3(2y-1)}{10^{-4} + ( (-2x+0.5)^2 + (2y-1)^2 )^{2.5} )}.
\]

\begin{figure}[ht]
\begin{center}
\begin{minipage}{.7\textwidth}
\hspace{-1cm}
\includegraphics[width=\textwidth]{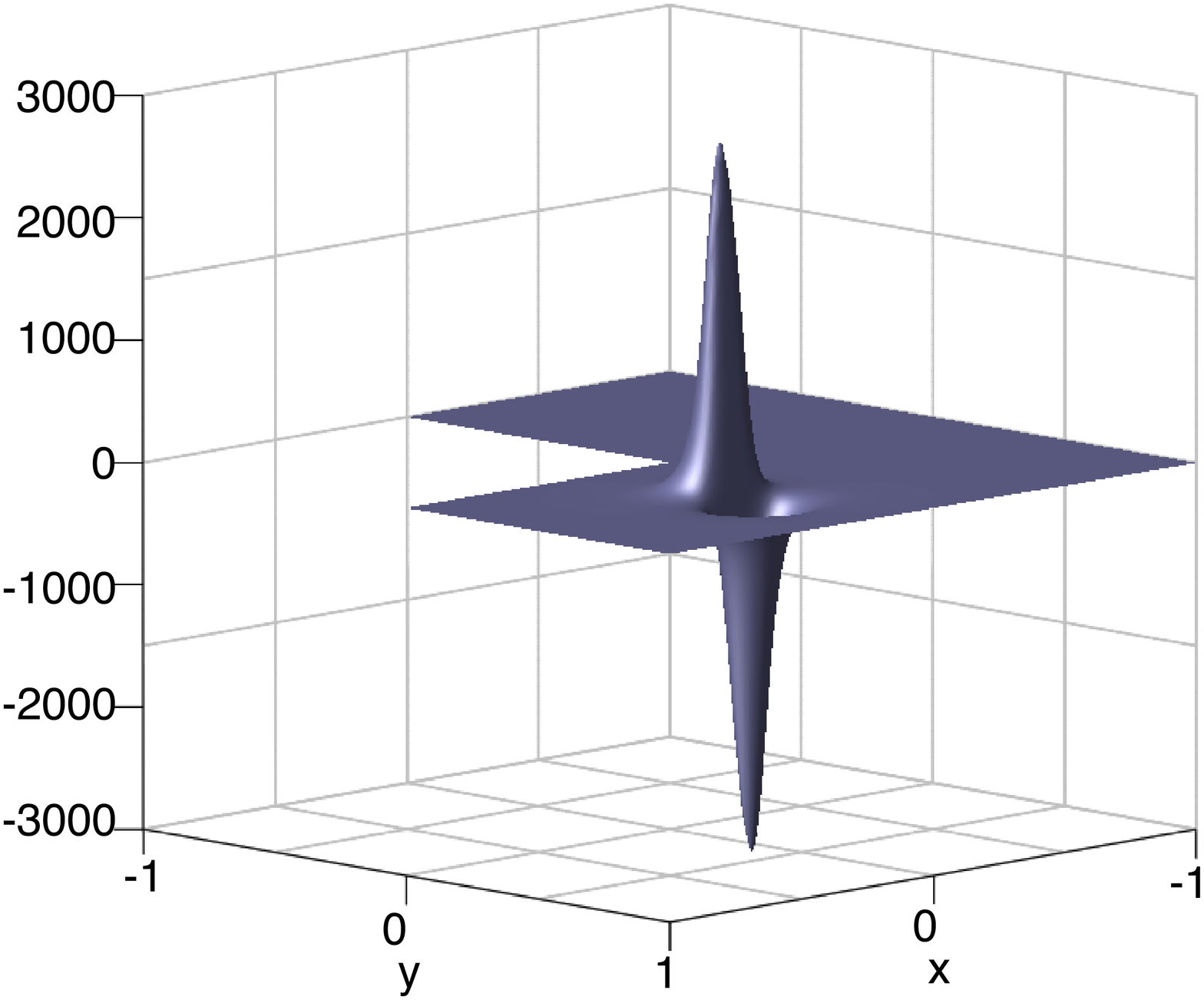}
\hspace{-1cm}
\end{minipage}
\begin{minipage}{.29\textwidth}
\hspace{-1cm}
\includegraphics[width=\textwidth]{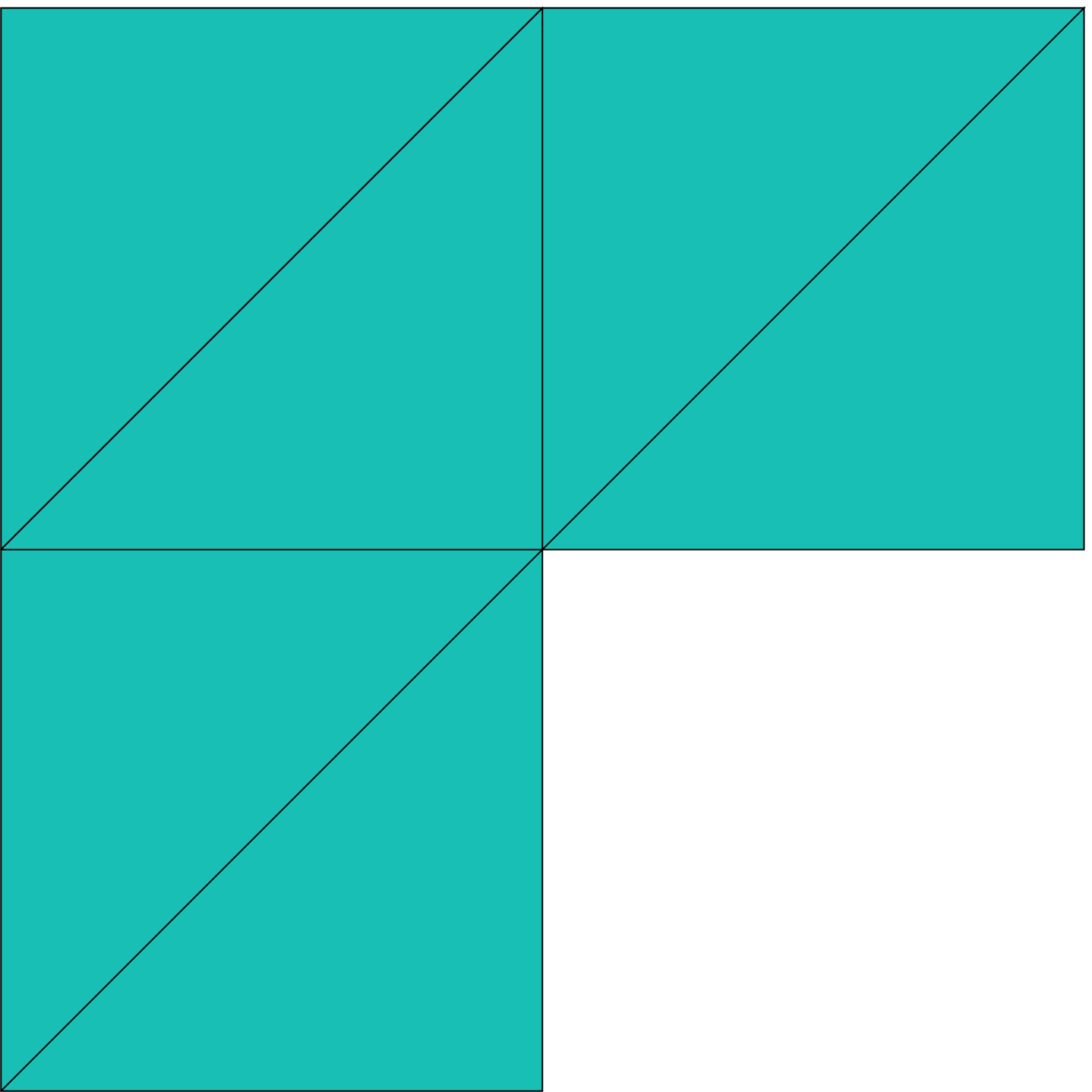}
\end{minipage}
\end{center}
\caption{Example 2: Source term $\source^O(x,y)$ associated {with} the second quantity of interest (left) and initial mesh.}
\label{fig:Lshaped_output_initial_mesh}
\end{figure}
Figure \ref{fig:Lshaped_output_initial_mesh} shows the source term of the adjoint problem associated {with} $s_2$ and the initial mesh for all the computations.

\begin{figure}[ht!]
\psfrag{a}[cc][cc][0.85]{Number of elements}
\psfrag{b}[cc][cc][0.85]{Half bound gap}
\psfrag{c}[Bl][Bl][0.75]{uniform $p=1$}
\psfrag{d}[Bl][Bl][0.75]{uniform $p=2$}
\psfrag{e}[Bl][Bl][0.75]{uniform $p=3$}
\psfrag{f}[Bl][Bl][0.75]{adapt $\Delta_\text{tol}$ , $p=1$}
\psfrag{g}[Bl][Bl][0.75]{adapt $\Delta_\text{tol}$ , $p=2$}
\psfrag{h}[Bl][Bl][0.75]{adapt $\Delta_\text{tol}$ , $p=3$}
\psfrag{i}[Bl][Bl][0.75]{adapt $\Theta = 0.5$ , $p=1$}
\psfrag{j}[Bl][Bl][0.75]{adapt  $\Theta = 0.5$ , $p=2$}
\psfrag{k}[Bl][Bl][0.75]{adapt  $\Theta = 0.5$ , $p=3$}
\psfrag{l}[Bl][Bl][0.75]{adapt  $\Theta = 0.25$ , $p=1$}
\psfrag{m}[Bl][Bl][0.75]{adapt  $\Theta = 0.25$ , $p=2$}
\psfrag{n}[Bl][Bl][0.75]{adapt  $\Theta = 0.25$ , $p=3$}
\psfrag{o}[Bl][Bl][.85]{$\numel^{-2/3}$}
\psfrag{p}[Bl][Bl][0.75]{$\numel^{-2}$}
\psfrag{q}[Bl][Bl][0.75]{$\numel^{-4}$}
\psfrag{r}[Bl][Bl][0.75]{$\numel^{-6}$}
\centering
\includegraphics[width=.8\textwidth]{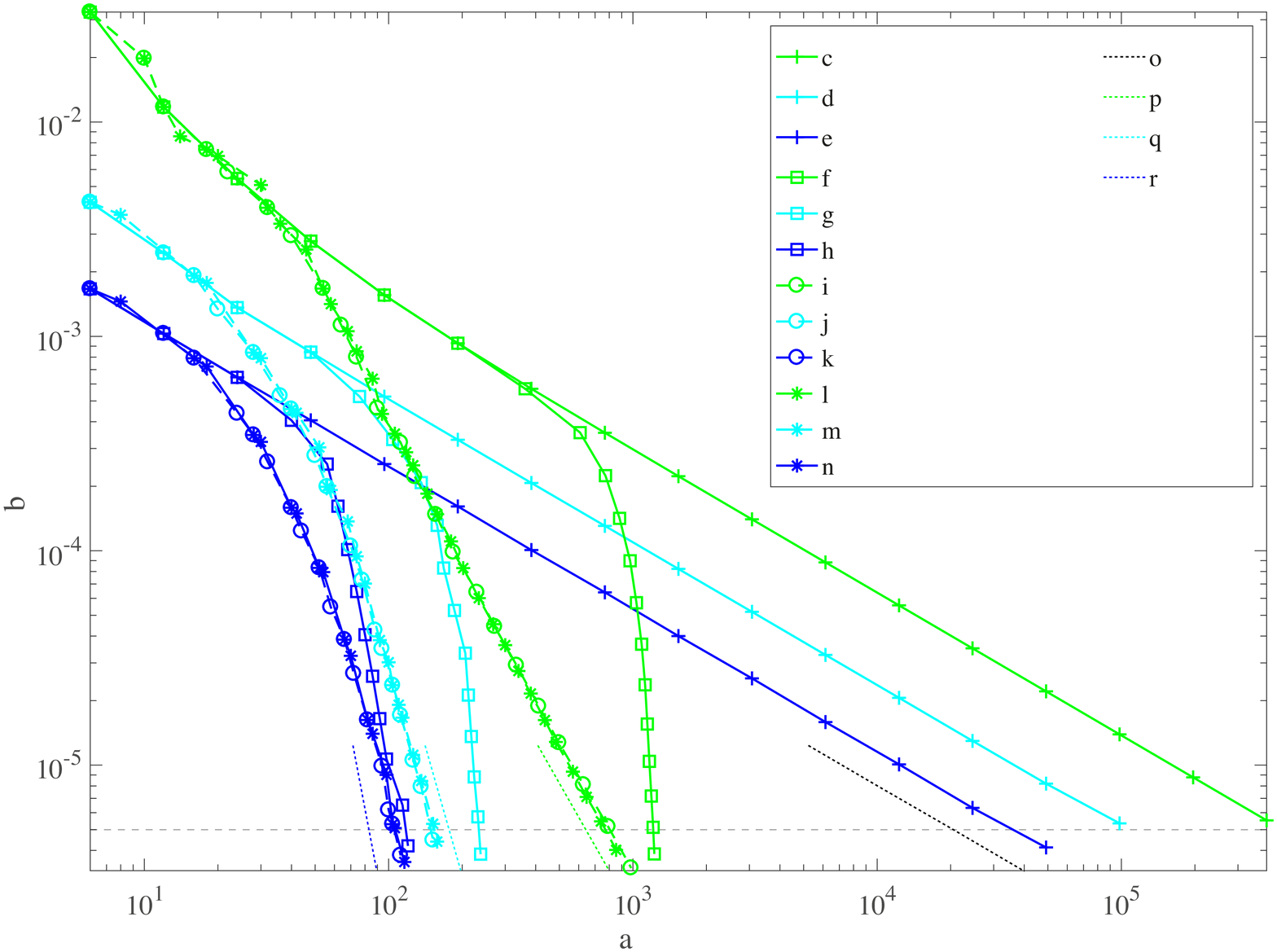}
\caption{Example 2: $s_1={\norm{\nu \grad u}}^2$  - Convergence of the half bound gap for both uniform and adaptive mesh refinements.}
\label{fig:Lshaped_halfgap}
\end{figure}

The behavior of the proposed strategy is first shown for the energy quantity of interest, $s_1={\norm{\nu\grad u}}^2$, using both a uniform mesh refinement (where in each step each triangle is bisected splitting its longest edge) and three different criteria for the adaptive procedure. The three adaptive procedures are all associated {with} a final bound gap of $\Delta_\text{tol}=10^{-5}$ (or an equivalence target for the half bound gap of $0.5\cdot 10^{-5}$): the first adaptive strategy assumes a uniform error distribution while the two others use a bulk criterion with $\Theta = 0.5$ and $\Theta=0.25$. 
Figure \ref{fig:Lshaped_halfgap} shows the convergence of the half bound gap obtained from the HDG approximations of order $p=1, 2$ and $3$.
As can be seen, using a uniform mesh refinement the expected convergence rate of $\mathcal{O}(\numel^{-2/3})$ is achieved, since the HDG method in this case convergences as $\mathcal{O}(\numel^{-1/3})$ regardless of the value of $p$, see \cite{MR1813251}. 
The adaptive strategies using both bulk criterions asymptotically converge as $\mathcal{O}(\numel^{-2p})$. On the other hand, the adaptive strategy based on a uniform error distribution assumption reaches the same accuracy with a similar number of elements, but with a very different convergence behavior. In the initial steps of the adaptive procedure, the meshes are uniformly refined resulting in a slow convergence, and once the adaptive strategy starts refining the elements around the singularity, convergence is reached in few iterations. 

\begin{figure}[ht]
\begin{center}
\begin{minipage}{.24\textwidth}
\includegraphics[width=\textwidth]{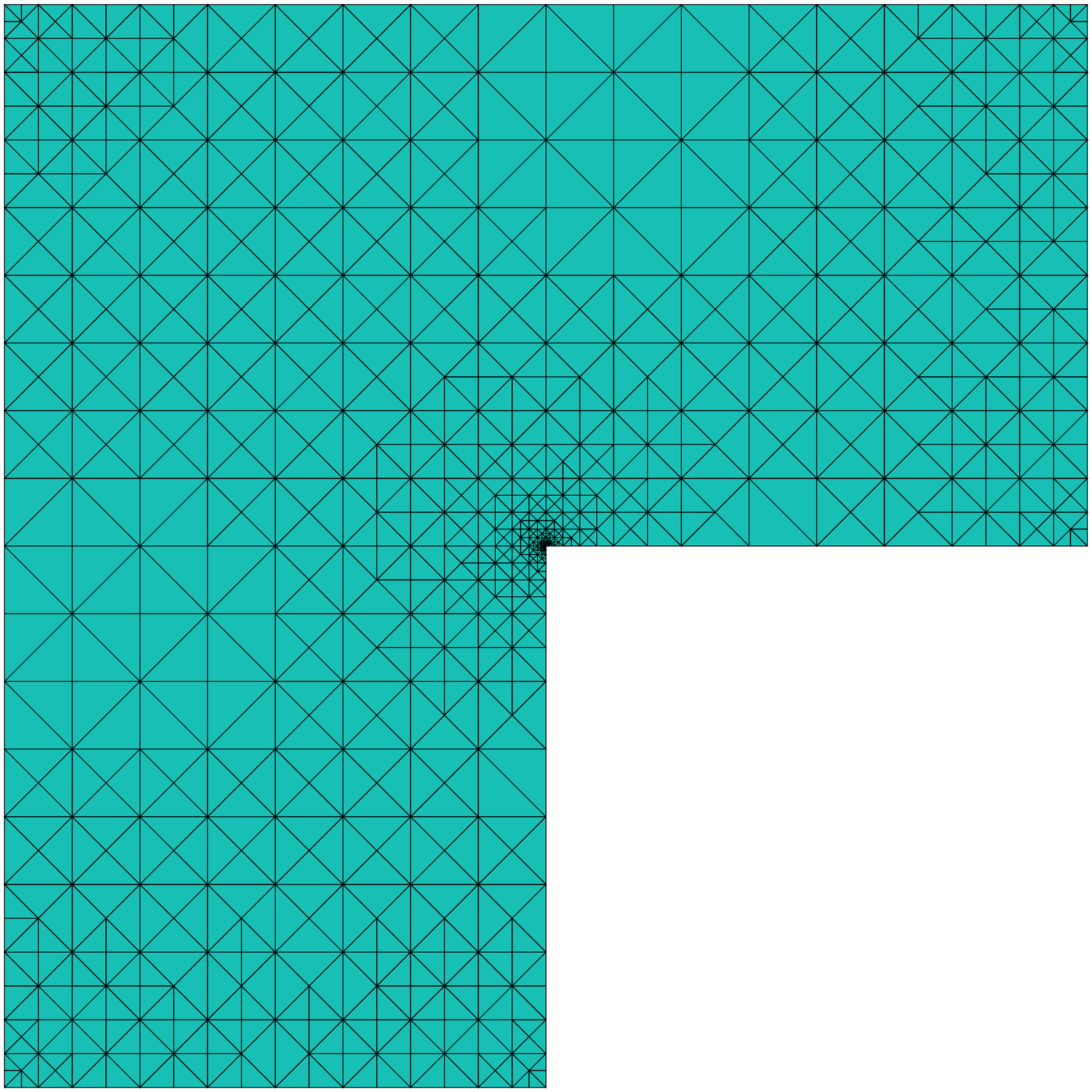}\\
\includegraphics[width=\textwidth]{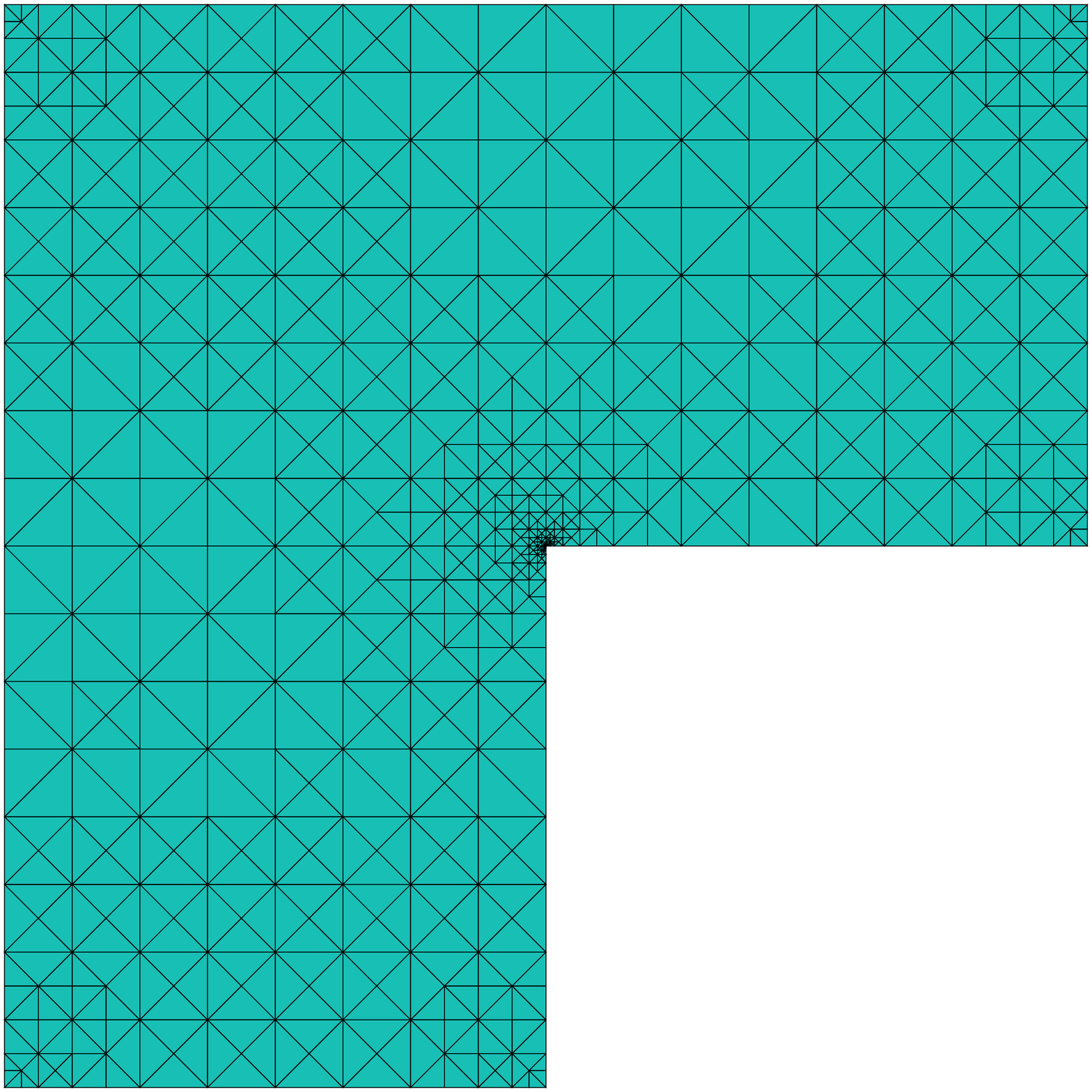}\\
\includegraphics[width=\textwidth]{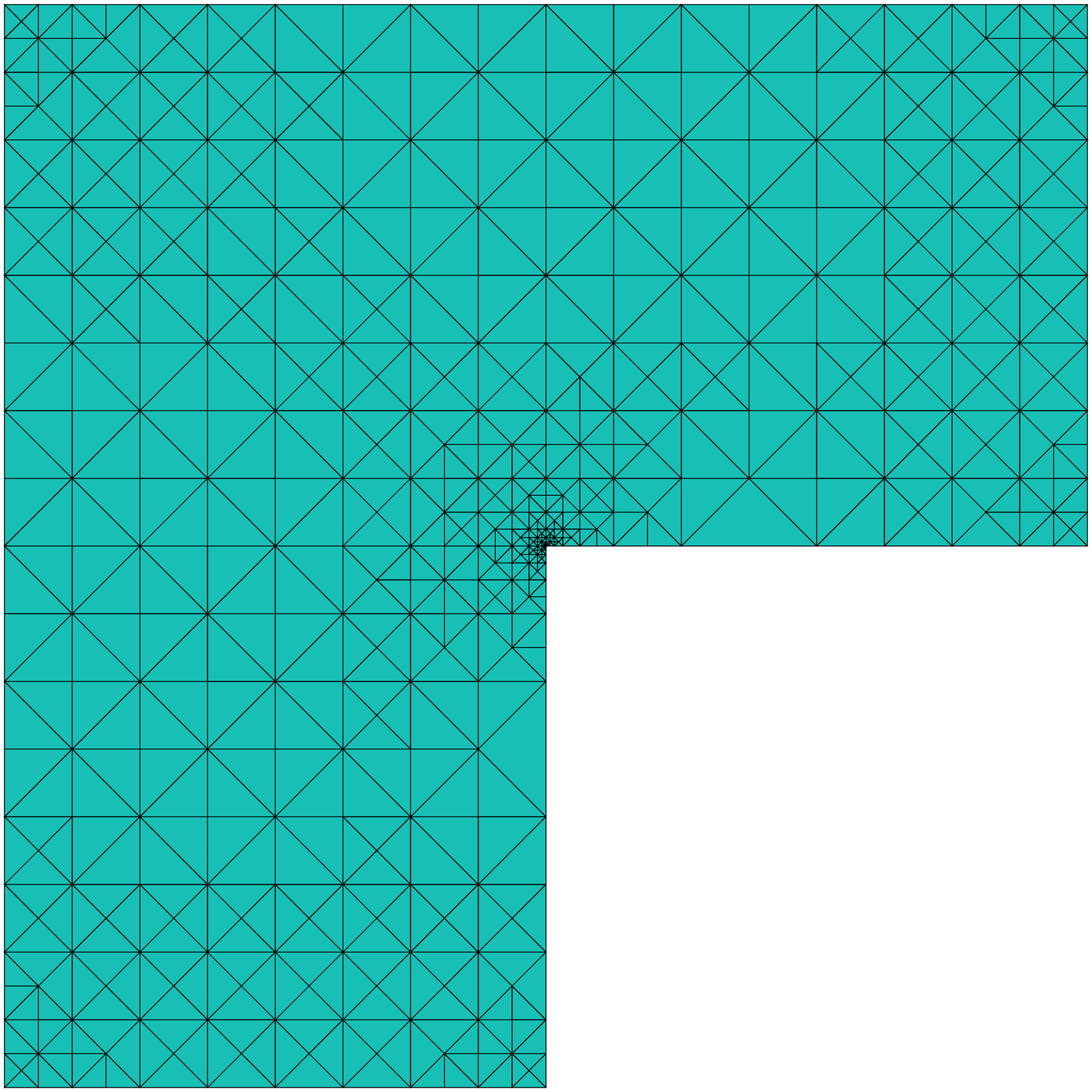}\\
\end{minipage}
\begin{minipage}{.24\textwidth}
\includegraphics[width=\textwidth]{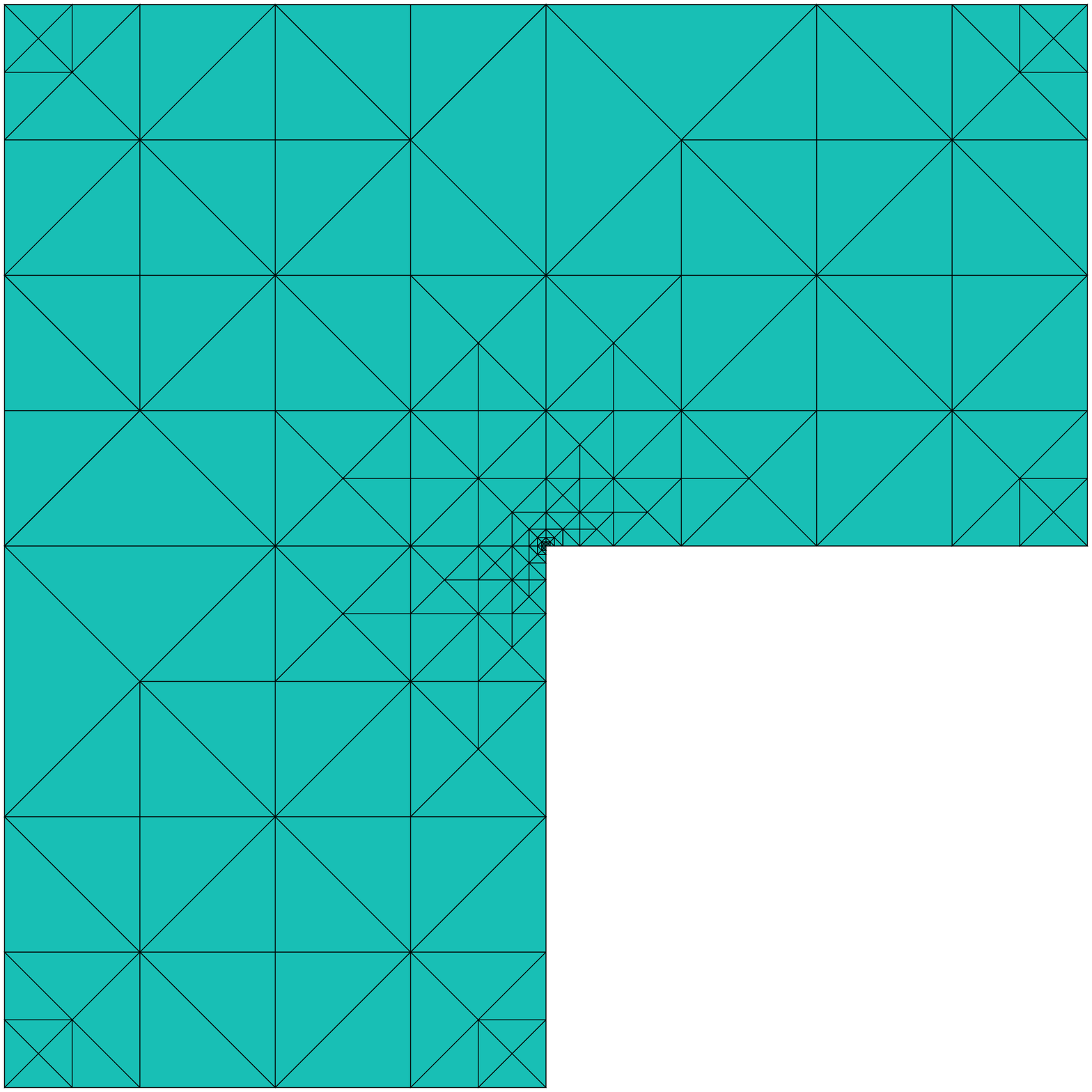}\\
\includegraphics[width=\textwidth]{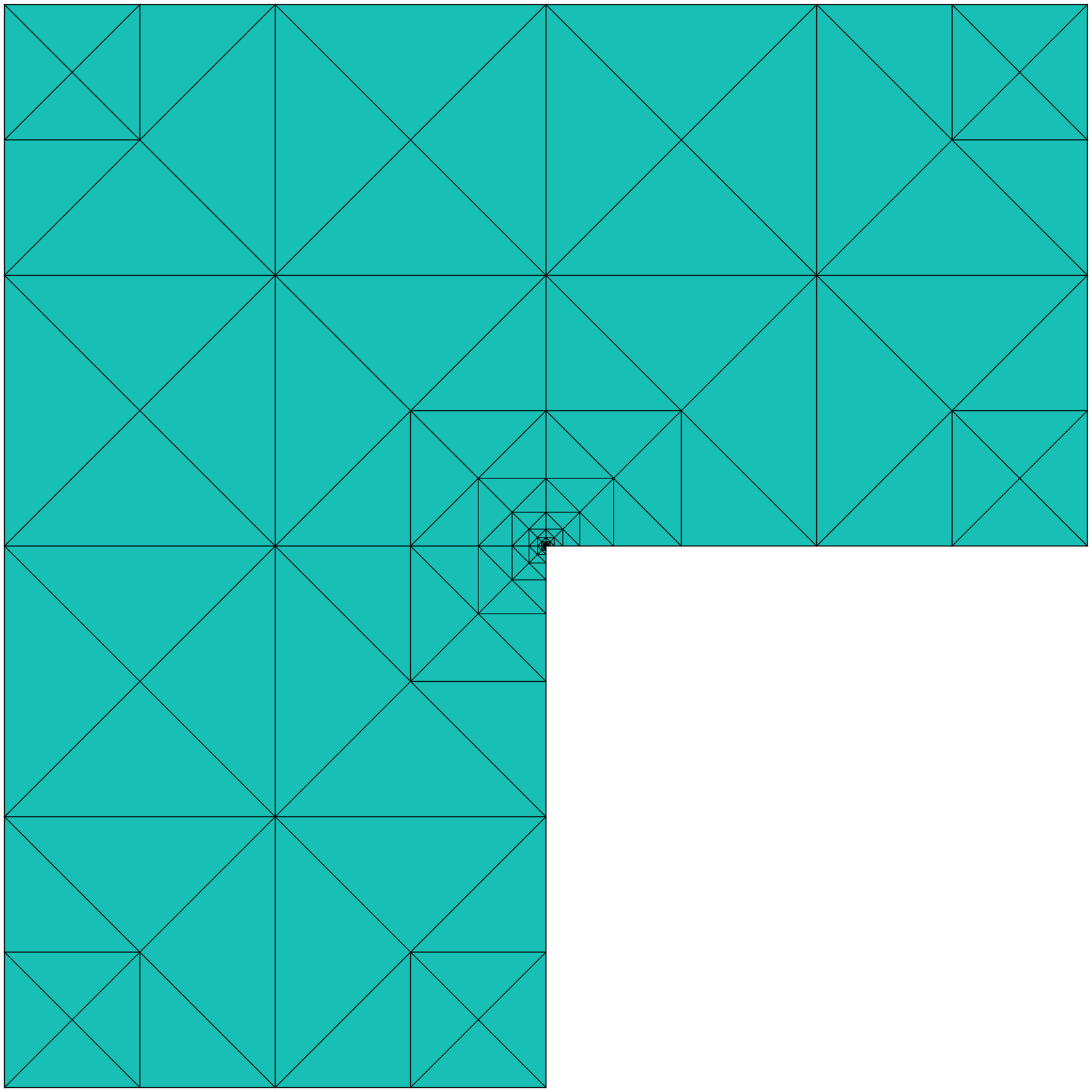}\\
\includegraphics[width=\textwidth]{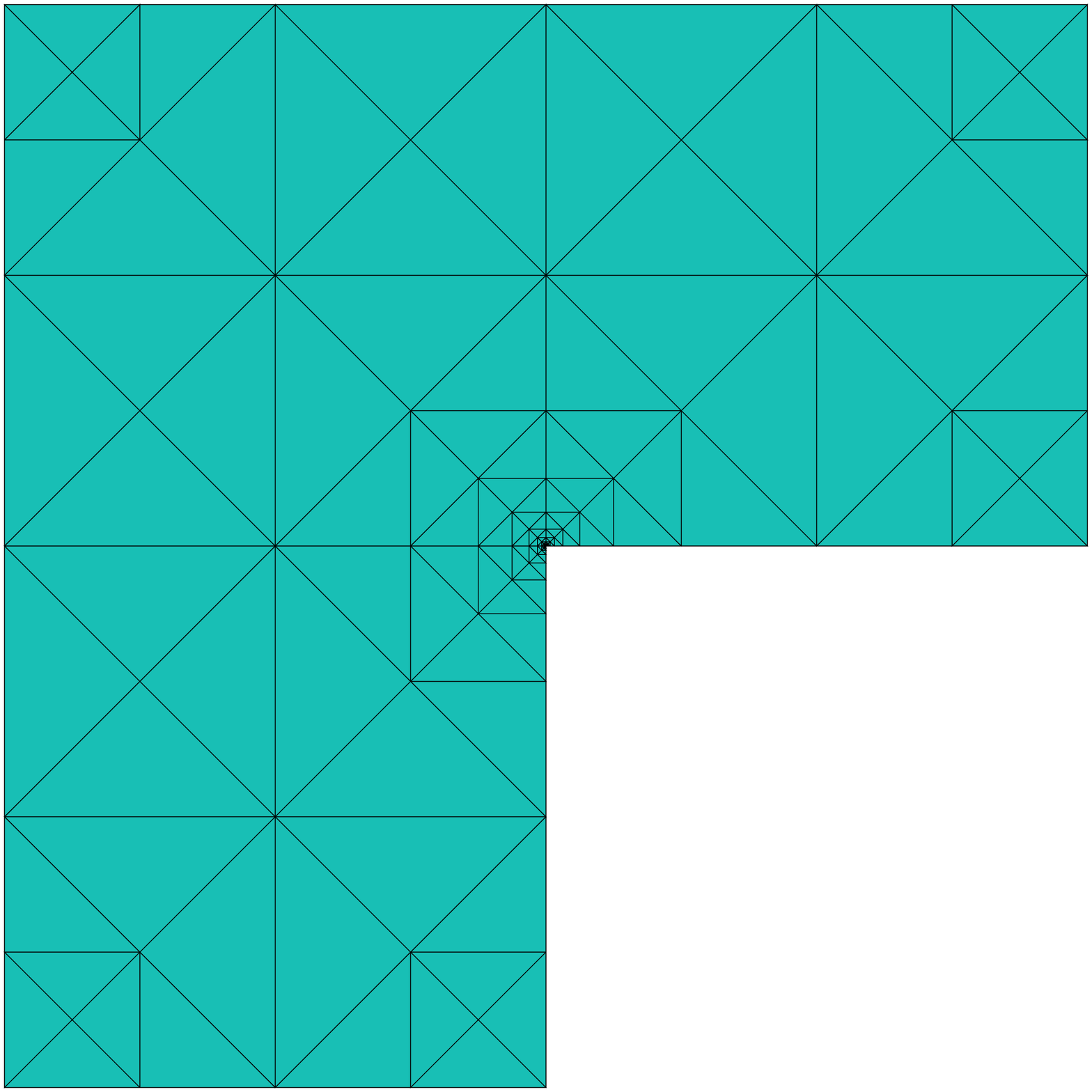}\\
\end{minipage}
\begin{minipage}{.24\textwidth}
\includegraphics[width=\textwidth]{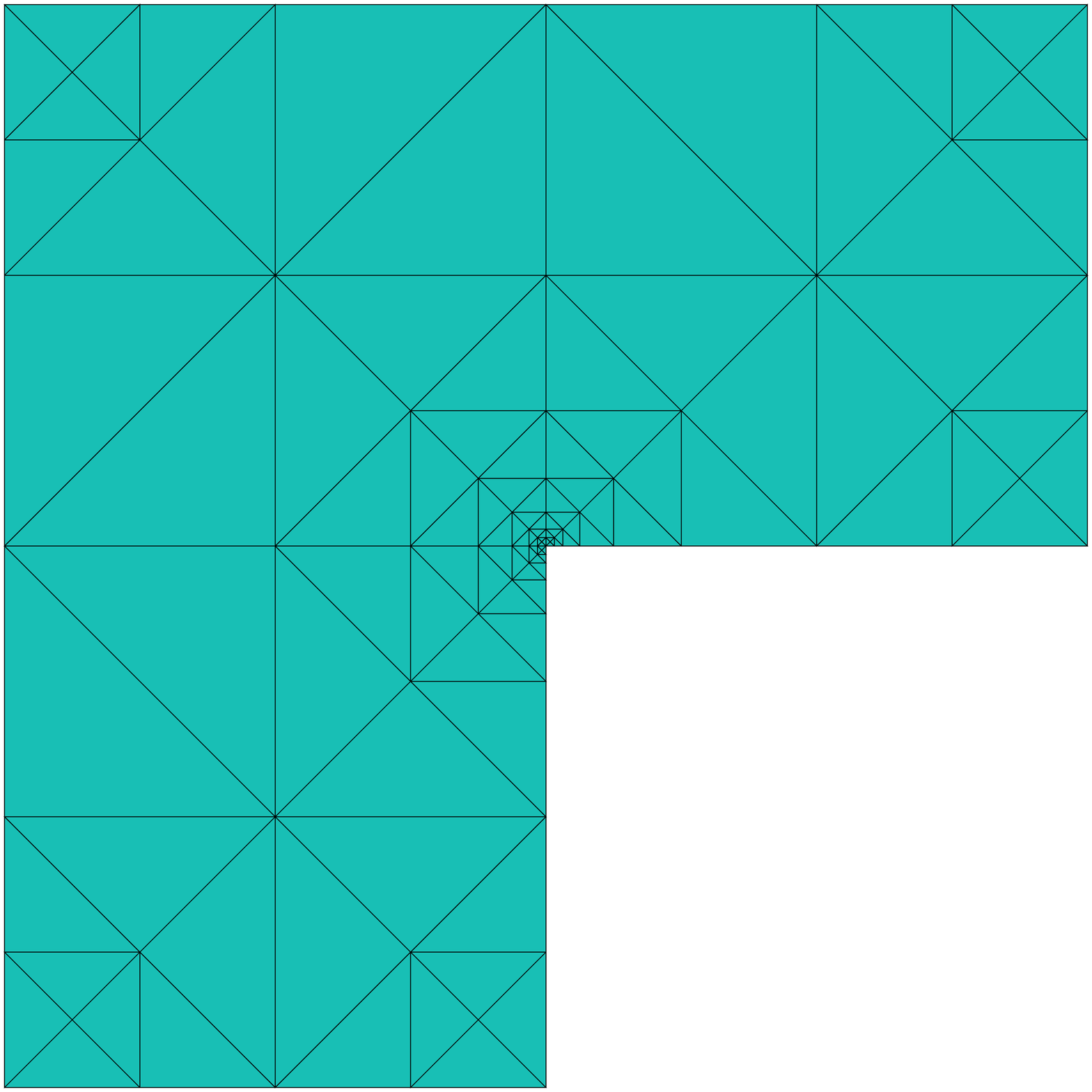}\\
\includegraphics[width=\textwidth]{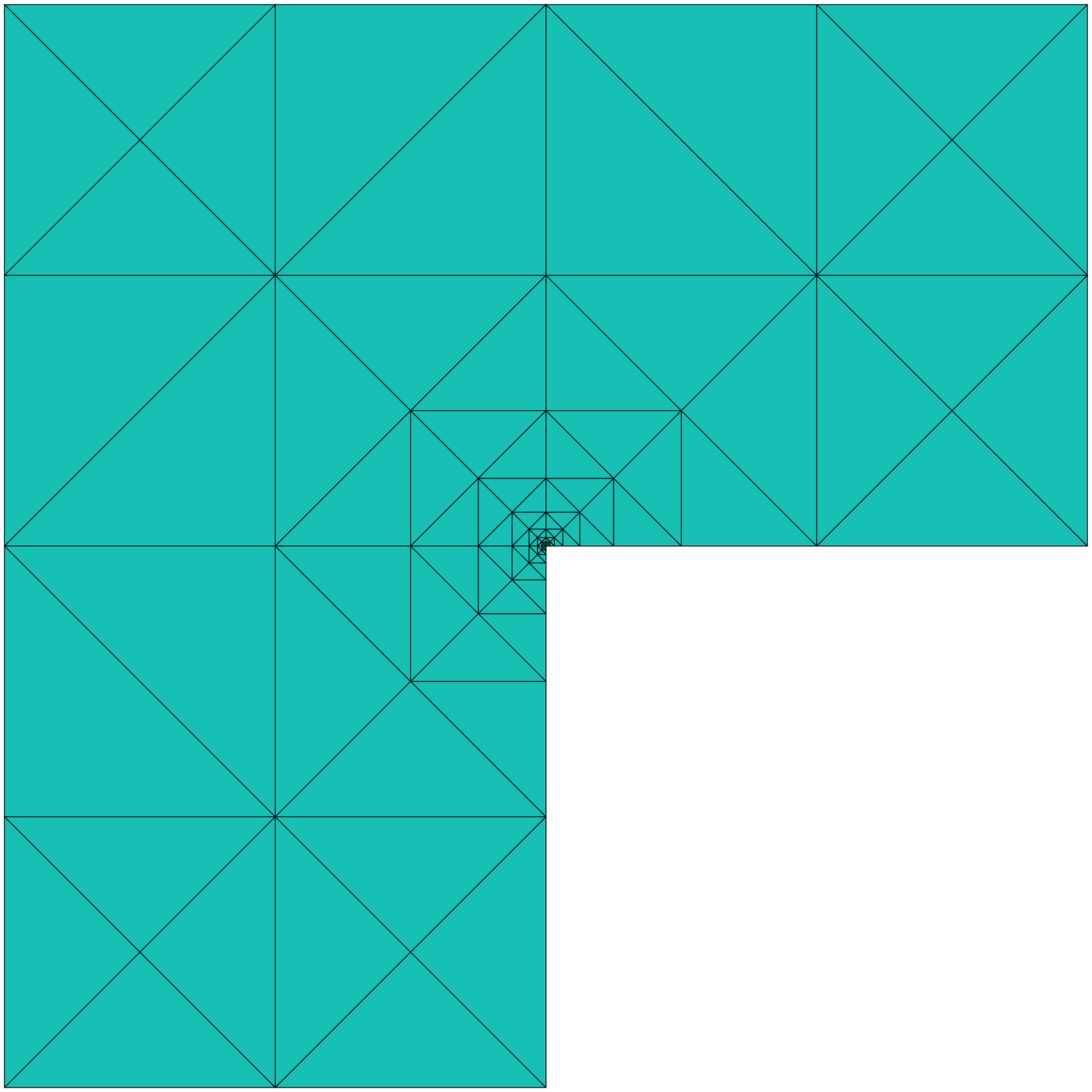}\\
\includegraphics[width=\textwidth]{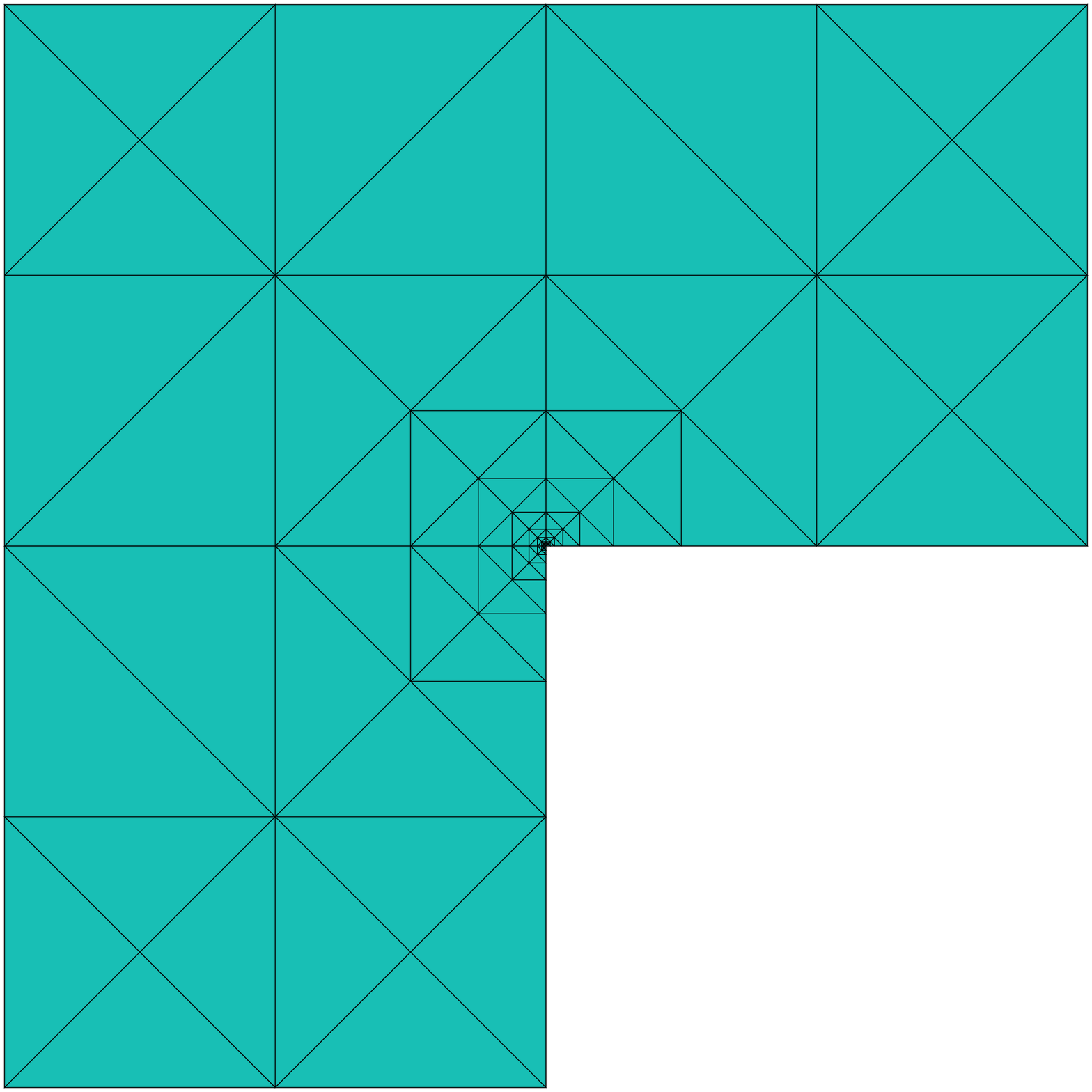}\\
\end{minipage}
\end{center}
\caption{Example 2: $s_1={\norm{\nu \grad u}}^2$  - Final meshes of the adaptive procedure ($\Delta_\textrm{tol}=10^{-5}$ top, $\Theta=0.5$ middle $\Theta=0.25$ bottom) for $p=1, 2, 3$ from left to right.}
\label{fig:Lshaped_final_meshes} 
\end{figure}

The final meshes of the adaptive procedures are shown in Figure \ref{fig:Lshaped_final_meshes}. As can be seen, all the adaptive strategies provide similar final meshes (although the intermediate meshes vary significantly in the first steps of the adaptive procedures when using a uniform error distribution strategy than when using a bulk criterion). 
Also, since the adaptive strategies converge like $\mathcal{O}(\numel^{-2p})$, there is a clear difference between the final meshes associated {with} $p=1$ and $p>1$. 

Table \ref{table:Lshaped} summarizes the results associated {with} the initial mesh, the intermediate iterations associated {with} a half bound gap lower than $0.5\cdot10^{-3}, 0.5\cdot10^{-4}$ and for the final mesh where $\Delta_\text{tol}/2<0.5\cdot10^{-5}$. The results for $\Theta=0.25$ are omitted since they are similar to the ones associated {with} $\Theta=0.5$.
\begin{table}[hbt!]
\centering
{\small
\begin{tabular}{c|cc|cc|c|cc} 
\hline
&& $\numel$ & $\slb$  & $\sub$ & $\sh \pm \Delta_h/2$ & $|s-\sh|$ 
& $|s-s_h|$ 
\\\hline
\hline        
\parbox[t]{2mm}{\multirow{16}{*}{\rotatebox[origin=c]{90}{$p=1$}}} 	
& \parbox[t]{2mm}{\multirow{4}{*}{\rotatebox[origin=c]{90}{uniform}}} 	
&        6 & 0.1740651 & 0.2392014 & 0.2066332 $\pm$ 3.26e-02 & 7.44e-03 & 2.62e-03 \\ &
&      768 & 0.2136240 & 0.2143344 & 0.2139792 $\pm$ 3.55e-04 & 9.66e-05 & 2.64e-04 \\ &
&    24576 & 0.2140310 & 0.2141012 & 0.2140661 $\pm$ 3.50e-05 & 9.69e-06 & 2.86e-05 \\ &
&   393216 & 0.2140687 & 0.2140798 & 0.2140743 $\pm$ 5.52e-06 & 1.53e-06 & 4.54e-06 \\ 
\cline{2-8} &
\parbox[t]{2mm}{\multirow{4}{*}{\rotatebox[origin=c]{90}{$\Delta_\textrm{tol}=10^{-5}$}}} 	
&        6 & 0.1740651 & 0.2392014 & 0.2066332 $\pm$ 3.26e-02 & 7.44e-03  & 2.62e-03 \\ &
&      608 & 0.2136232 & 0.2143347 & 0.2139789 $\pm$ 3.56e-04 & 9.69e-05  & 2.53e-04 \\ &
&     1086 & 0.2140289 & 0.2141023 & 0.2140656 $\pm$ 3.66e-05 & 1.02e-05  & 7.50e-06 \\ &
&     1224 & 0.2140708 & 0.2140786 & 0.2140747 $\pm$ 3.85e-06 & 1.14e-06  & 1.83e-05 \\ 
\cline{2-8} &
\parbox[t]{2mm}{\multirow{4}{*}{\rotatebox[origin=c]{90}{$\Theta=0.5$}}} 	
&        6 & 0.1740651 & 0.2392014 & 0.2066332 $\pm$ 3.26e-02 & 7.44e-03 & 2.62e-03 \\ &
&       90 & 0.2134465 & 0.2143698 & 0.2139081 $\pm$ 4.62e-04 & 1.68e-04 & 9.74e-04 \\ &
&      272 & 0.2140174 & 0.2141062 & 0.2140618 $\pm$ 4.43e-05 & 1.40e-05 & 1.51e-04 \\ &
&      984 & 0.2140714 & 0.2140781 & 0.2140748 $\pm$ 3.31e-06 & 1.04e-06 & 2.48e-05 \\ 
\hline
\hline        
\parbox[t]{2mm}{\multirow{16}{*}{\rotatebox[origin=c]{90}{$p=2$}}} 	
&\parbox[t]{2mm}{\multirow{4}{*}{\rotatebox[origin=c]{90}{uniform}}} 	
&        6 & 0.2084763 & 0.2169298 & 0.2127031 $\pm$ 4.23e-03 & 1.37e-03 & 2.52e-03 \\ &
&      192 & 0.2136675 & 0.2143265 & 0.2139970 $\pm$ 3.29e-04 & 7.88e-05 & 2.89e-04 \\ &
&     6144 & 0.2140352 & 0.2141007 & 0.2140680 $\pm$ 3.27e-05 & 7.84e-06 & 2.90e-05 \\ &
&    98304 & 0.2140694 & 0.2140802 & 0.2140748 $\pm$ 5.35e-06 & 1.03e-06 & 4.57e-06 \\ 
\cline{2-8} &
\parbox[t]{2mm}{\multirow{4}{*}{\rotatebox[origin=c]{90}{$\Delta_\textrm{tol}=10^{-5}$}}} 	
&        6 & 0.2084763 & 0.2169298 & 0.2127031 $\pm$ 4.23e-03 & 1.37e-03 & 2.52e-03 \\ &
&      104 & 0.2136664 & 0.2143269 & 0.2139966 $\pm$ 3.30e-04 & 7.92e-05 & 2.87e-04 \\ &
&      206 & 0.2140344 & 0.2141011 & 0.2140677 $\pm$ 3.33e-05 & 8.07e-06 & 2.85e-05 \\ &
&      238 & 0.2140709 & 0.2140787 & 0.2140748 $\pm$ 3.84e-06 & 1.03e-06 & 2.48e-06 \\ 
\cline{2-8} &
\parbox[t]{2mm}{\multirow{4}{*}{\rotatebox[origin=c]{90}{$\Theta=0.5$}}} 	
&        6 & 0.2084763 & 0.2169298 & 0.2127031 $\pm$ 4.23e-03 & 1.37e-03 & 2.52e-03 \\ &
&       40 & 0.2135426 & 0.2144586 & 0.2140006 $\pm$ 4.58e-04 & 7.52e-05 & 3.61e-04 \\ &
&       88 & 0.2140251 & 0.2141104 & 0.2140678 $\pm$ 4.26e-05 & 8.03e-06 & 3.39e-05 \\ &
&      152 & 0.2140700 & 0.2140790 & 0.2140745 $\pm$ 4.47e-06 & 1.30e-06 & 6.65e-07 \\  
\hline
\hline        
\parbox[t]{2mm}{\multirow{16}{*}{\rotatebox[origin=c]{90}{$p=3$}}} 	
&\parbox[t]{2mm}{\multirow{4}{*}{\rotatebox[origin=c]{90}{uniform}}} 	
&        6 & 0.2120143 & 0.2153474 & 0.2136809 $\pm$ 1.67e-03 & 3.95e-04 & 1.49e-03 \\ &
&       48 & 0.2135839 & 0.2143962 & 0.2139901 $\pm$ 4.06e-04 & 8.57e-05 & 3.72e-04 \\ &
&     1536 & 0.2140275 & 0.2141076 & 0.2140676 $\pm$ 4.00e-05 & 8.23e-06 & 3.72e-05 \\ &
&    49152 & 0.2140709 & 0.2140793 & 0.2140751 $\pm$ 4.13e-06 & 7.19e-07 & 3.70e-06 \\ 
\cline{2-8} &
\parbox[t]{2mm}{\multirow{4}{*}{\rotatebox[origin=c]{90}{$\Delta_\textrm{tol}=10^{-5}$}}} 	
&        6 & 0.2120143 & 0.2153474 & 0.2136809 $\pm$ 1.67e-03 & 3.95e-04 & 1.49e-03 \\ &
&       40 & 0.2135838 & 0.2143962 & 0.2139900 $\pm$ 4.06e-04 & 8.58e-05 & 3.72e-04 \\ &
&       80 & 0.2140266 & 0.2141079 & 0.2140672 $\pm$ 4.06e-05 & 8.59e-06 & 3.76e-05 \\ &
&      120 & 0.2140706 & 0.2140791 & 0.2140749 $\pm$ 4.20e-06 & 9.09e-07 & 3.86e-06 \\ 
\cline{2-8} &
\parbox[t]{2mm}{\multirow{4}{*}{\rotatebox[origin=c]{90}{$\Theta=0.5$}}} 	
&        6 & 0.2120143 & 0.2153474 & 0.2136809 $\pm$ 1.67e-03 & 3.95e-04 & 1.49e-03 \\ &
&       24 & 0.2135489 & 0.2144235 & 0.2139862 $\pm$ 4.37e-04 & 8.96e-05 & 3.99e-04 \\ &
&       66 & 0.2140296 & 0.2141068 & 0.2140682 $\pm$ 3.85e-05 & 7.60e-06 & 3.71e-05 \\ &
&      112 & 0.2140713 & 0.2140790 & 0.2140752 $\pm$ 3.79e-06 & 6.38e-07 & 3.95e-06 \\ 
\hline
\end{tabular}
}
\caption{Example 2: $s_1={\norm{\nu \grad u}}^2$  - Bounds for both uniform and adaptive mesh refinements for $p=1, 2$ and $3$.}
\label{table:Lshaped}
\end{table}
It is worth noting that as expected, the half bound gap $\Delta_h/2$ provides indeed an upper bound for the error in the quantity of interest associated {with} $\sh$, namely, $s-\sh$. In fact, even though the bounding procedure is devised to minimize the bound gap and not to produce accurate upper bounds for $s-\sh$, the effectivities measuring the quality of the half bound gap as an upper bound of $s-\sh$ are quite good in most cases. 

\begin{figure}[ht!]
\psfrag{a}[cc][cc][0.85]{Number of elements}
\psfrag{b}[cc][cc][0.85]{Half bound gap}
\psfrag{c}[Bl][Bl][0.75]{uniform $p=1$}
\psfrag{d}[Bl][Bl][0.75]{uniform $p=2$}
\psfrag{e}[Bl][Bl][0.75]{uniform $p=3$}
\psfrag{f}[Bl][Bl][0.75]{adapt $\Delta_\text{tol}$ , $p=1$}
\psfrag{g}[Bl][Bl][0.75]{adapt $\Delta_\text{tol}$ , $p=2$}
\psfrag{h}[Bl][Bl][0.75]{adapt $\Delta_\text{tol}$ , $p=3$}
\psfrag{i}[Bl][Bl][0.75]{adapt $\Theta = 0.5$ , $p=1$}
\psfrag{j}[Bl][Bl][0.75]{adapt  $\Theta = 0.5$ , $p=2$}
\psfrag{k}[Bl][Bl][0.75]{adapt  $\Theta = 0.5$ , $p=3$}
\psfrag{p}[Bl][Bl][0.75]{$\numel^{-2}$}
\psfrag{q}[Bl][Bl][0.75]{$\numel^{-4}$}
\psfrag{r}[Bl][Bl][0.75]{$\numel^{-6}$}
\centering
\includegraphics[width=.8\textwidth]{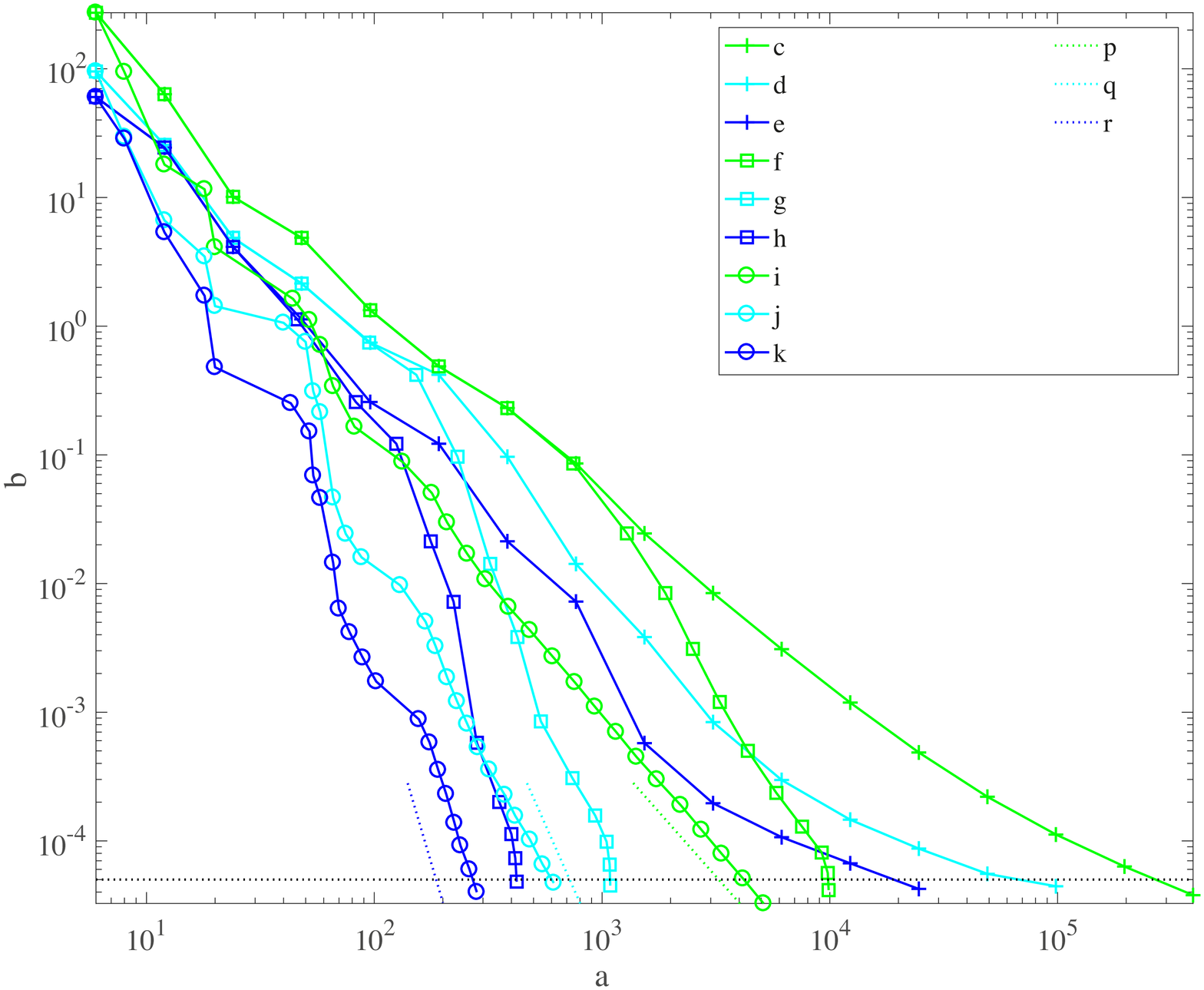}
\caption{Example 2: $s_2$ - Convergence of the half bound gap for both uniform and adaptive mesh refinements.}
\label{fig:Lshaped_QoI_halfgap}
\end{figure}

The results associated {with} the second quantity of interest are shown in Figure \ref{fig:Lshaped_QoI_halfgap} for $\Delta_\text{tol}=10^{-4}$. Two adaptive strategies are used: the uniform error distribution assumption and the bulk criterion for $\Theta=0.5$.
It can be seen that in the first iterations of both the uniform and adaptive refinements, the estimators are governed by the large data oscillation errors associated {with} the adjoint problem yielding to pessimistic bounds. However, since the data oscillation errors are of high order, after few iterations the half bound gaps converge as expected.
It is again clear that using higher order elements is advantageous, because for about the same accuracy high order elements result in meshes with fewer triangles and less global edge degrees of freedom. Also, the order of convergence of the adaptive procedures for larger values of $p$ makes a difference in the computational effort required to achieve a desired prescribed tolerance.

\begin{figure}[ht!]
\begin{center}
\includegraphics[width=.3\textwidth]{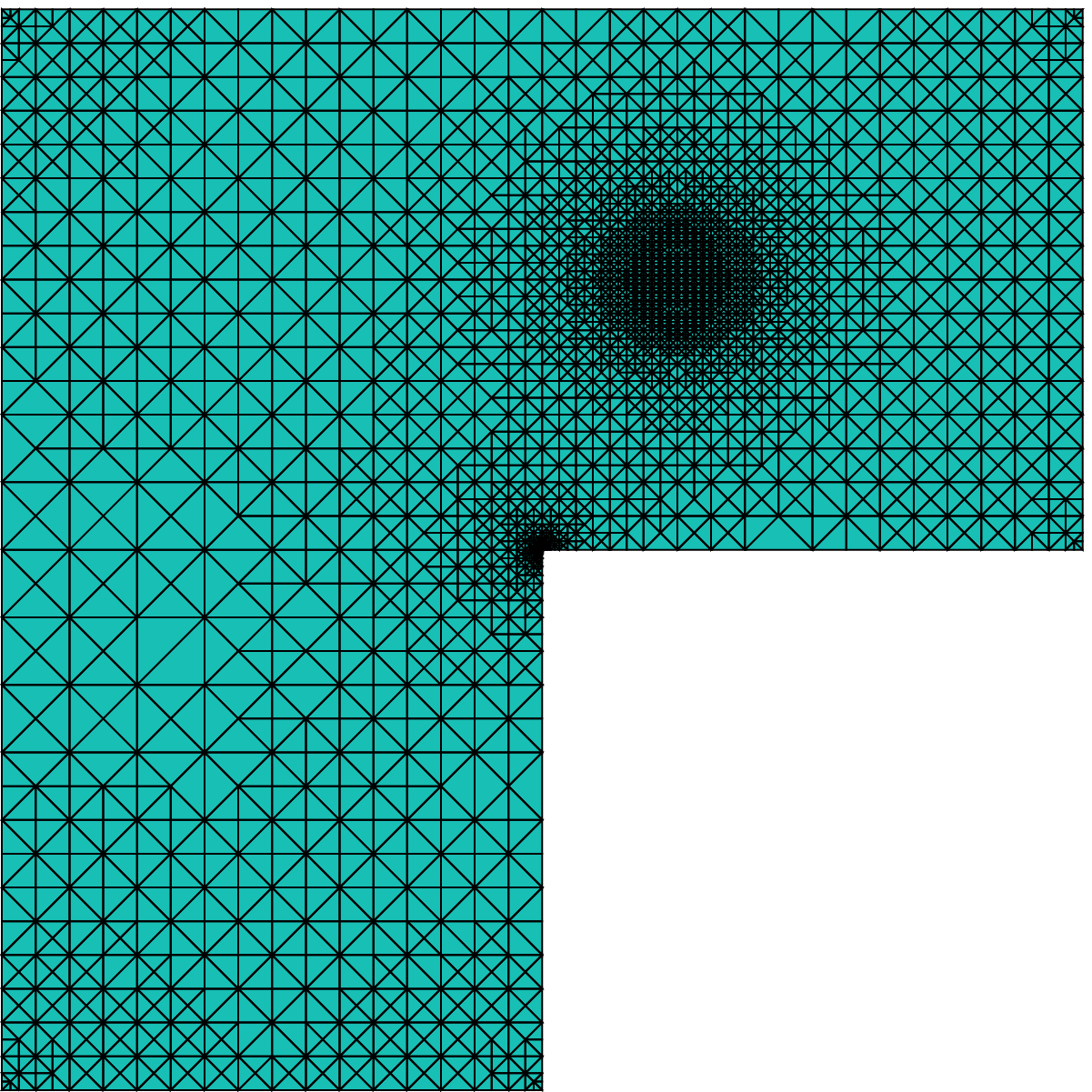}
\includegraphics[width=.3\textwidth]{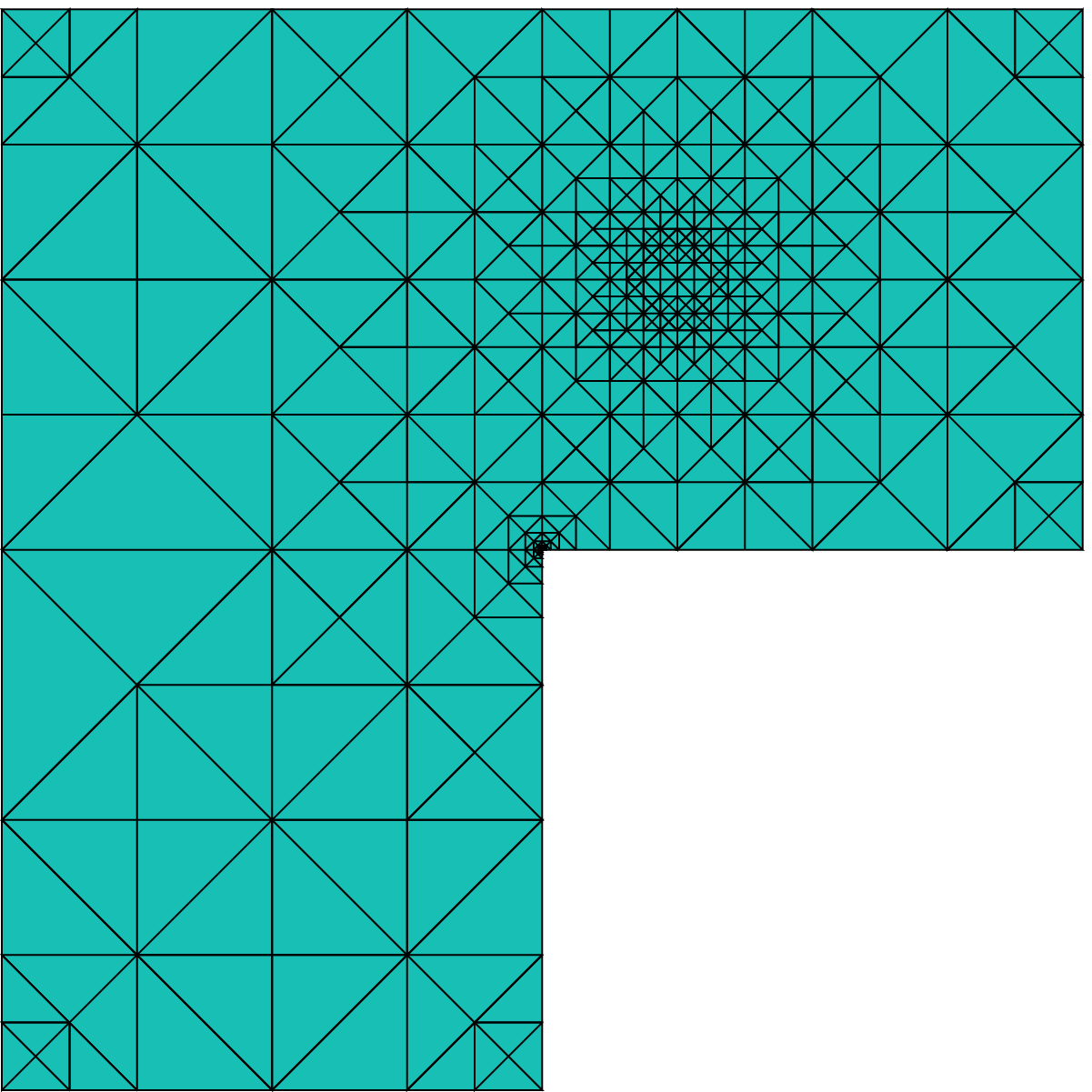}
\includegraphics[width=.3\textwidth]{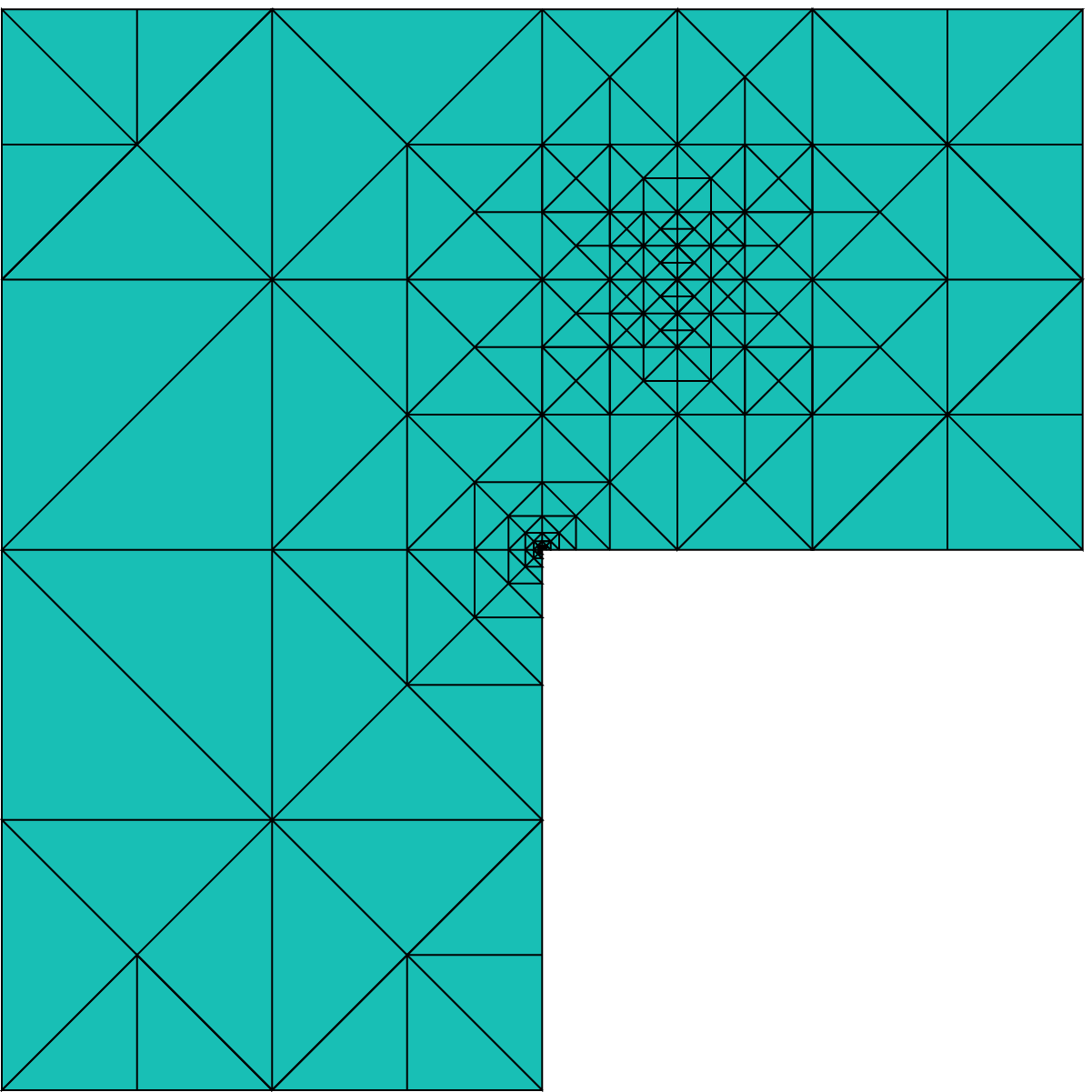}

\includegraphics[width=.3\textwidth]{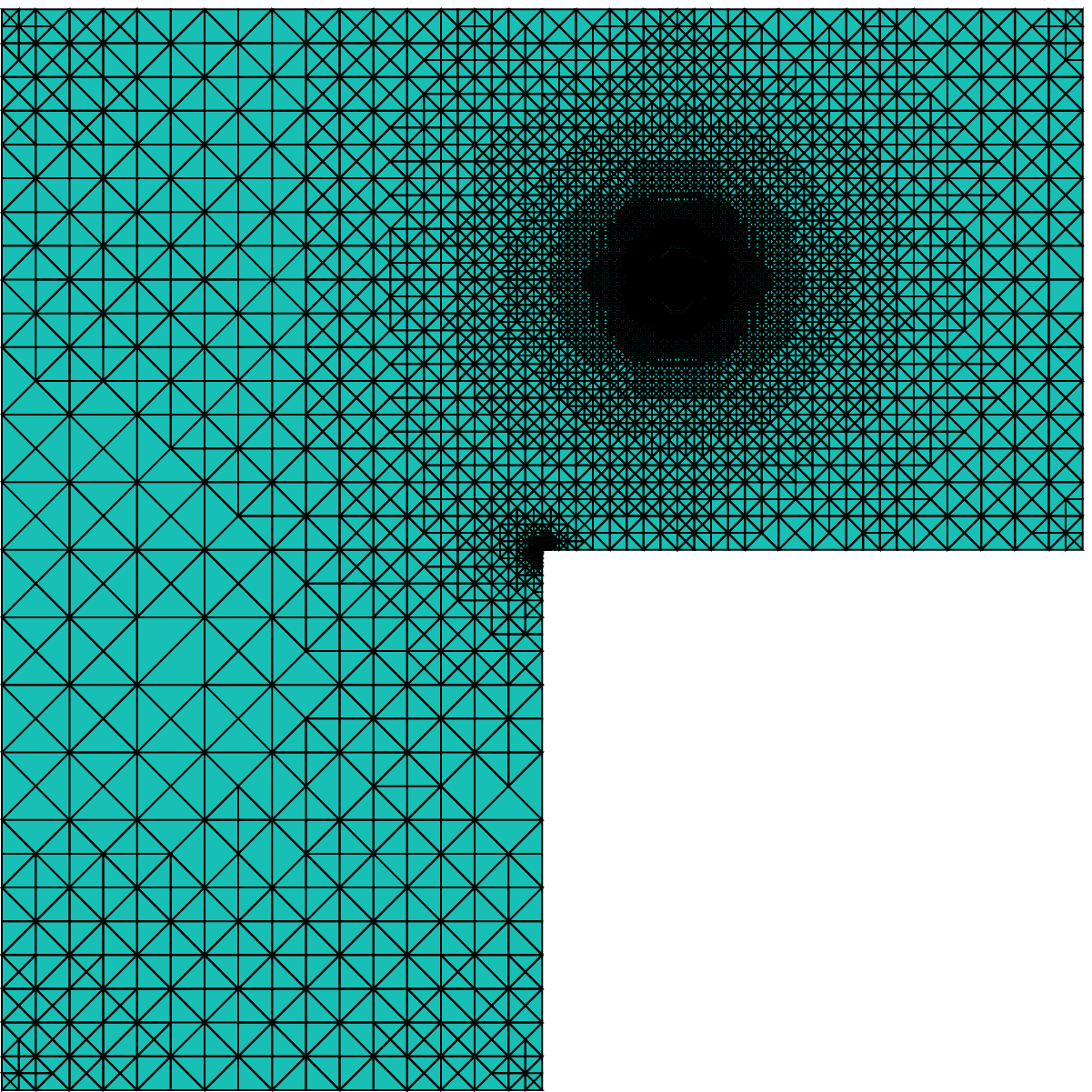}
\includegraphics[width=.3\textwidth]{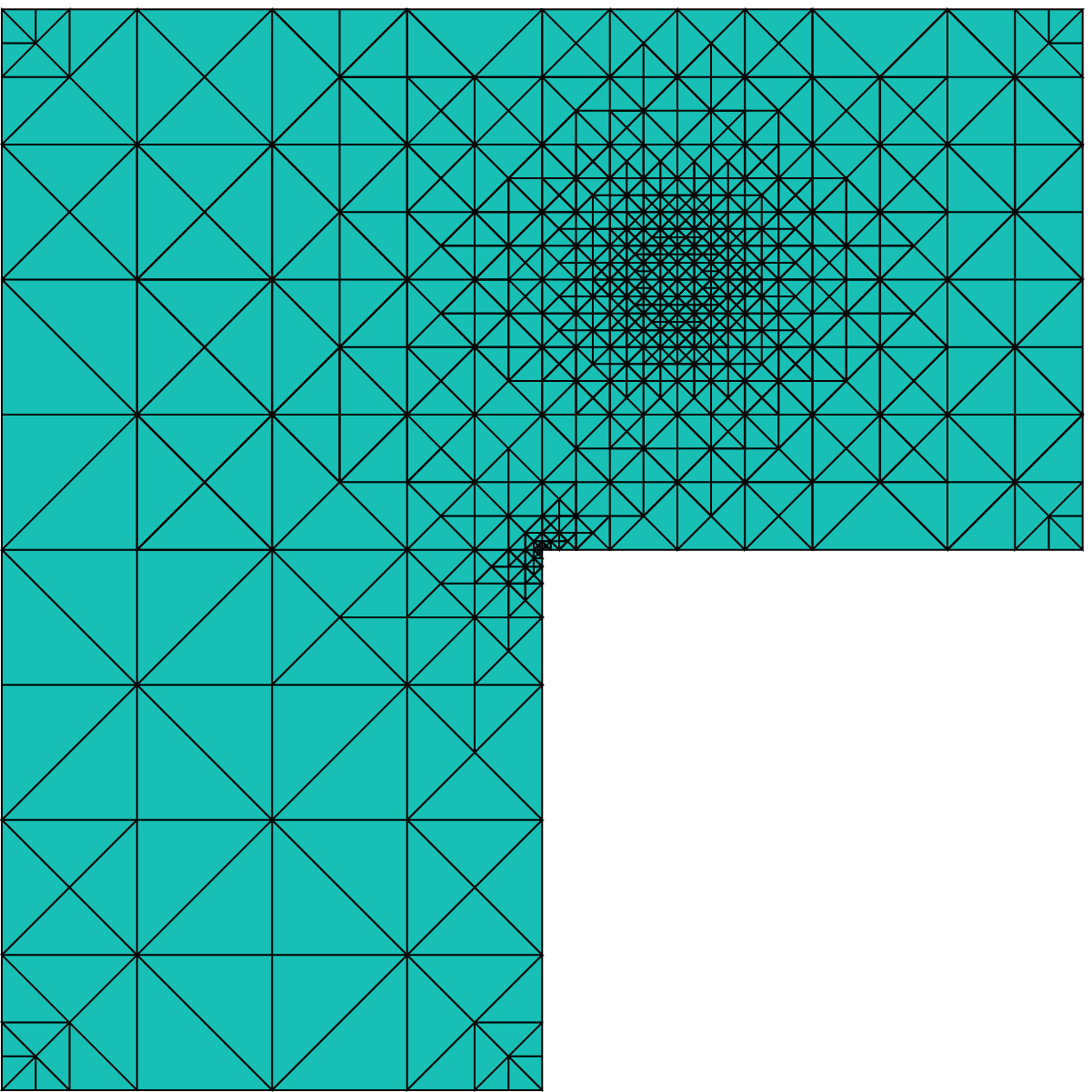}
\includegraphics[width=.3\textwidth]{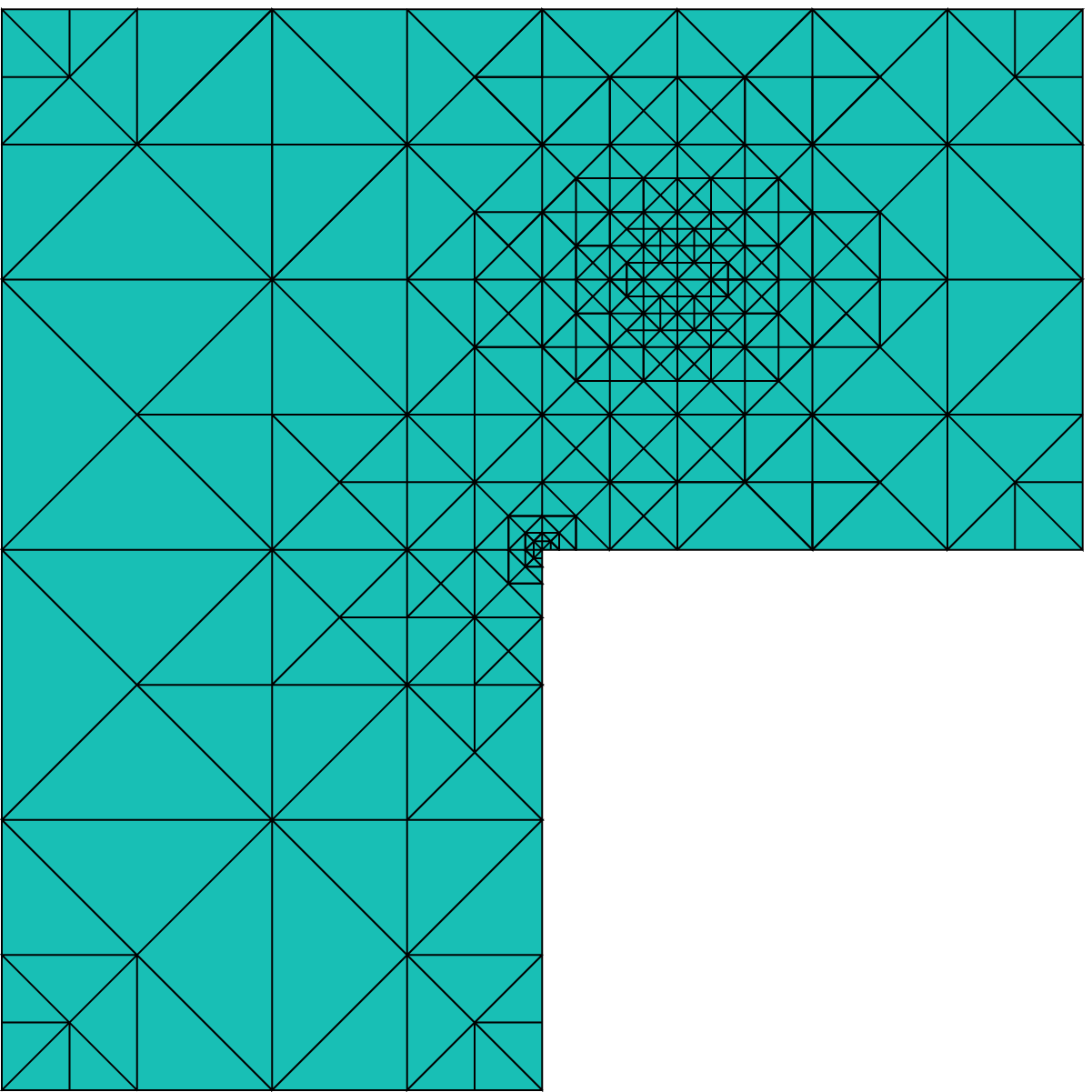}
\end{center}
\caption{Example 2: $s_2$  - Final meshes of the adaptive procedure associated {with} $p=1, 2$ and  $3$ from left to right and $\Theta=0.5$ (top) with meshes with $5119$, $614$ and $282$ triangles respectively
and a uniform error distribution with $\Delta_\text{tol}=10^{-4}$ (bottom) where the meshes have $9882$, $1086$ and $422$ triangles respectively.}
\label{fig:Lshaped_QoI_final_meshes} 
\end{figure}

The final meshes obtained in the adaptive procedures are shown in Figure \ref{fig:Lshaped_QoI_final_meshes}. 
As can be appreciated the adaptive procedure refines both in the corner singularity and in the area where the source term of the adjoint problem presents a large gradient and large data oscillation errors appear. It can also be seen that in this case, the bulk criterion yields coarser meshes for the same accuracy.

\section{Concluding remarks}
\label{sec:rem}

A general framework to compute guaranteed lower and upper bounds for quantities of interest from potential and equilibrated (or zero-oder equilibrated) flux reconstructions is presented. The bounds are guaranteed regardless of the size of the underlying finite element mesh and regardless of the kind of data (the source term and the Neumann boundary conditions are not required to be piecewise polynomial functions). 

In particular, bounds for quantities of interest from HDG approximations of both the primal and adjoint problems are obtained. Properly exploiting the superconvergence properties of local post-processed HDG approximations {yields} optimal convergence curves for the bound gap in the quantity of interest, both using uniform and adaptive mesh refinements.

Two numerical examples are presented to demonstrate the accuracy of the proposed technique when using HDG approximations. The obtained results seem to confirm the superconvergent properties of the bounds and show that using high-order HDG approximations yields very accurate bounds for the quantity of interest, even for very coarse meshes.

\section{Acknowledgements}

This work was partially supported by the Spanish Ministry of Economy and Competitiveness (Grant numbers: PGC2018-097257-B-C33 / DPI2017-85139-C2-2-R) and the Generalitat de Catalunya (Grant number: 2017-SGR-1278). 


\bibliography{ref_error}

\appendix
\section{Output bounds from potential and equilibrated flux reconstructions -- proof of Theorem \ref{th:property_QoI_01}}
\label{App:eq:property_QoI_01}

The key ingredient to prove Theorem \ref{th:property_QoI_01}  is the reformulation of the output of interest as a constrained minimization. This reasoning is similar to the approaches introduced for conforming non-mixed approximations \cite{MR1487947,MR1948321,MR2114293,MR2116366,MR2186143}.
We begin by writing the quantity of interest $s=\ell^O(u,\q)$ as a constrained minimization problem
\begin{equation}\label{eq:app01:aux01}
\begin{array}{ll}
\pm s &
\displaystyle
= \inf\limits_{(w,\bv)\in \test\times\testq} \pm\ell^O(w,\bv) + \kappa (a(w,\bv;w,\bv) - \ell(w,\bv))\\[1ex]
&\qquad\qquad  \text{s.t.} \quad a(w,\bv;\phi,\bvarphi) = \ell(\phi,\bvarphi) \quad \forall (\phi,\bvarphi)\in \test\times\testq,
\end{array}
\end{equation}
where $\kappa\in(0,+\infty)$ is an arbitrary parameter. The above statement {is}  easily verified by noting that, from \eqref{eq:general_variational_formulation}, the constraint $a(w,\bv;\phi,\bvarphi) = \ell(\phi,\bvarphi), \forall (\phi,\bvarphi)\in \test\times\testq$ is only satisfied when $(w,\bv)=(u,\q)$ due to the uniqueness of the solution and clearly $a(u,\q;\phi,\bvarphi) = \ell(\phi,\bvarphi) $. 
Now, the Lagrangian associated with the above constrained minimization problem is given by
\[
\Lag(w,\bv;\phi,\bvarphi)=\pm\ell^O(w,\bv) + \kappa (a(w,\bv;w,\bv) - \ell(w,\bv))+a(w,\bv;\phi,\bvarphi) - \ell(\phi,\bvarphi),
\]
and problem \eqref{eq:app01:aux01} becomes
\begin{equation}\label{eq:app01:aux02}
\pm s= \inf\limits_{(w,\bv)\in \test\times\testq}\sup\limits_{(\phi,\bvarphi)\in \test\times\testq} 
	\Lag(w,\bv;\phi,\bvarphi).
\end{equation}
Bounds for the output $s$ can be easily found using  the strong duality of the convex optimization problem and the saddle point property of the Lagrange multipliers as
\begin{equation}\label{eq:app01:aux03}
\pm s 
	= \sup\limits_{(\phi,\bvarphi)\in \test\times\testq} \inf\limits_{(w,\bv)\in \test\times\testq}
	\Lag(w,\bv;\phi,\bvarphi)
	\geq \inf\limits_{(w,\bv)\in \test\times\testq}
	\Lag(w,\bv;\uhPA^\mp,\qhPA^\mp)\equiv \pm s^\mp_h
	\quad \forall(\uhPA^\mp,\qhPA^\mp)
	\in \test\times\testq,
\end{equation}
where in order to obtain sharp bounds, it is important to use a good approximation $(\uhPA^\mp,\qhPA^\mp)$ of the Lagrange multipliers. Note that the explicit dependence of $s^\mp_h$ on $\kappa$ is omitted here for simplicity of presentation.

The explicit expression for the bounds $s^\mp_h$ associated {with} a particular choice of $(\uhPA^\mp,\qhPA^\mp)$ is found imposing the stationary conditions, that is, requiring the variations of $\Lag(w,\bv;\uhPA^\mp,\qhPA^\mp)$ with respect to $w$ and $\bv$ vanish.
From the definition of $\Lag(w,\bv;\uhPA^\mp,\qhPA^\mp)$
and taking into account \eqref{eq:a_wv_wv} it is easy to see that
\begin{equation}\label{eq:app01:aux04}
\begin{array}{l}
\Lag(w+\delta w,\bv+\delta\bv;\uhPA^\mp,\qhPA^\mp) - \Lag(w,\bv;\uhPA^\mp,\qhPA^\mp) 
\\[2ex]
\displaystyle\qquad=
\pm\ell^O(\delta w,\delta\bv)-\kappa\ell(\delta w,\delta\bv) + a(\delta w,\delta \bv;\uhPA^\mp,\qhPA^\mp)
+ 2\kappa (\nu^{-1}\bv,\delta \bv)
+ \kappa(\nu^{-1}\delta \bv,\delta \bv)
\\[2ex]
\displaystyle\qquad=
(\pm f^O-\kappa f-\nabla\cdot\qhPA^\mp,\delta w)
+\langle \pm \trac^O+\kappa \trac+\qhPA^\mp\cdot\bn,\delta w\rangle _{\GN}
\\[2ex]
\displaystyle\qquad\qquad
+  (2\kappa\nu^{-1}\bv+\nu^{-1}\qhPA^\mp-\grad  \uhPA^\mp,\delta \bv)
+\langle \pm \uD^O+\kappa \uD+\uhPA^\mp,\delta \bv\cdot\bn\rangle _{\GD}
+ \kappa(\nu^{-1}\delta \bv,\delta \bv),
\end{array}
\end{equation}
and therefore, denoting by $(w^\mp,\bv^\mp)$ the minimizers of $\Lag(w,\bv;\uhPA^\mp,\qhPA^\mp)$, 
the stationary conditions require the conditions given in equation \eqref{eq:reconstructionsComb} to hold.
\begin{equation}\label{eq:reconstructionsComb}
\begin{tabular}{| ll |}
\hline
Combined 
&$\uhPA^\mp\in\test$ \\[1ex]
primal/adjoint potential reconstruction:
&$\uhPA^\mp = \mp \uD^O-\kappa \uD \text{ on } \GD$ 
\\[1ex]\hline
Combined 
&$\qhPA^\mp\in\testq$ \\[1ex]
primal/adjoint equilibrated flux reconstruction:
&$\nabla\cdot\qhPA^\mp=\pm f^O-\kappa f \text{ in } \Omega$ \\[1ex]
&$\qhPA^\mp\cdot\bn=\mp \trac^O-\kappa \trac \text{ on } \GN$ 
\\\hline
Minimizer condition: & $w^\mp$ free\\ 
&$\bv^\mp = \dfrac{1}{2\kappa}\left(-\qhPA^\mp+\nu\grad  \uhPA^\mp\right)$
\\[1.5ex]\hline
\end{tabular}
\end{equation}
It is worth noting that the combined primal/adjoint potential and equilibrated flux reconstructions can be computed introducing the potential and equilibrated flux reconstructions of the primal and adjoint problems $(\uh,\qh)$ and $(\uhA,\qhA)$ satisfying \eqref{eq:reconstructions} as
\[
\uhPA^\mp = \mp \uhA -\kappa \uh \qquad,\qquad \qhPA^\mp = \pm \qhA -\kappa \qh.
\] 
Now, the expression for the bounds $s^\mp_h$ can be rewritten by first noting 
that the stationary condition \eqref{eq:app01:aux04} for the optimal values $(w^\mp,\bv^\mp)\in\test\times\testq$ implies that
\[
\pm\ell^O(\delta w,\delta\bv)-\kappa\ell(\delta w,\delta\bv) + a(\delta w,\delta \bv;\uhPA^\mp,\qhPA^\mp)
+ 2\kappa (\nu^{-1}\bv^\mp,\delta \bv)=0 \qquad\forall (\delta w,\delta\bv)\in\test\times \testq,
\]
which in particular holds for $(\delta w,\delta\bv)=(w^\mp,\bv^\mp)$, and using equation \eqref{eq:a_wv_wv} for $\delta\bv=\bv^\mp$. Inserting these expressions into the definition of $\Lag(w^\mp,\bv^\mp;\uhPA^\mp,\qhPA^\mp)$ yields, after some algebraic manipulations, to
\begin{equation}\label{eq:app01:aux05}
\pm s^\mp_h=\Lag(w^\mp,\bv^\mp;\uhPA^\mp,\qhPA^\mp)
=
-\dfrac{1}{4\kappa}{\norm{\qhPA^\mp-\nu\grad  \uhPA^\mp}}^2 - \ell(\uhPA^\mp,\qhPA^\mp).
\end{equation}
Also, using equation \eqref{eq:green}, it is easy to see that  the primal and adjoint equilibrated flux reconstructions satisfying \eqref{eq:reconstructions}
verify that  forall $(w,\bv)\in \test\times \testq$ 
\begin{subequations}\label{eq:primal_adjoint_weak}
\begin{align}
	&\displaystyle
\ell(w,\bv)
=-(\qh,\grad  w) + \langle \qh\cdot\bn ,w\rangle _{\GD}  - \langle \uD,\bv\cdot\bn\rangle _{\GD},
	\label{eq:primal_adjoint_weak_P}
	\\[1ex]&\displaystyle
\ell^O(w,\bv)
=-(\qhA,\grad  w) + \langle \qhA\cdot\bn ,w\rangle _{\GD}  + \langle \uD^O,\bv\cdot\bn\rangle _{\GD},	
	\label{eq:primal_adjoint_weak_A}
\end{align}
\end{subequations}
and therefore taking $(w,\bv)=(\uhPA^\mp,\qhPA^\mp)=(\mp \uhA -\kappa \uh, \pm \qhA -\kappa \qh)$ into \eqref{eq:primal_adjoint_weak_P} and $(w,\bv)=(\uh,\qh)$ into \eqref{eq:primal_adjoint_weak_A} yields, after some simplifications, 
\[
\ell(\uhPA^\mp,\qhPA^\mp)
=\mp\ell^O(\uh,\qh)
\mp(\qhA,\grad  \uh)
\pm(\qh,\grad  \uhA ) 
+\kappa(\qh,\grad   \uh).
\]
Finally, expanding ${\norm{\qhPA^\mp-\nu\grad  \uhPA^\mp}}^2$ and rearranging terms yields
\begin{equation}\label{eq:property_QoI_01_mod}
\begin{array}{ll}
\pm s^\mp_h
&
=\pm\ell^O(\uh,\qh)\mp (\qh+\nu\grad  \uh,\grad  \uhA)
-\dfrac{1}{4\kappa}{\norm{(\qhA +\nu\grad   \uhA)\mp\kappa(\qh+\nu\grad  \uh)}}^2
\\[2ex]
&
=\pm\ell^O(\uh,\qh)
-\dfrac{1}{4\kappa}{\norm{\qhA +\nu\grad   \uhA }}^2 
-\dfrac{\kappa}{4}{\norm{\qh+\nu\grad  \uh}}^2
\pm \dfrac{1}{2}(\nu^{-1}(\qh+\nu\grad  \uh),\qhA-\nu\grad  \uhA),
\end{array}
\end{equation}
and substituting the optimal value of $\kappa_{\rm opt} = {\norm{\qhA +\nu\grad   \uhA}/\norm{\qh+\nu\grad  \uh}}$ concludes the proof. Indeed, joining all the obtained expressions provides
\[
\pm s 
	\geq \pm s^\mp_h(\kappa_{\rm opt}) 
	=\pm\ell^O(\uh,\qh)
-\dfrac{1}{2}{\norm{\qh+\nu\grad  \uh}\,\norm{\qhA+\nu\grad  \uhA}}
\pm \dfrac{1}{2}(\nu^{-1}(\qh+\nu\grad  \uh),\qhA-\nu\grad  \uhA).
\]

\section{Output bounds from potential and zero-order equilibrated flux reconstructions
 -- proof of Theorem \ref{th:property_QoI_02}
}
\label{App:eq:bounds_projected}

Equation \eqref{eq:property_QoI_01_mod} shows that for any potential and equilibrated flux reconstructions of the primal and adjoint problems $(\uh,\qh)$ and $(\uhA,\qhA)$ then
\begin{equation}\label{eq:property_QoI_01_mod_mod}
\pm s \geq \pm s^\mp_h
=\pm\ell^O(\uh,\qh)\mp (\qh+\nu\grad  \uh,\grad  \uhA)
-\dfrac{1}{4\kappa}{\norm{ (\qhA +\nu\grad   \uhA)\mp\kappa(\qh+\nu\grad  \uh)}}^2.
\end{equation}
Moreover, the first two terms in 
\eqref{eq:property_QoI_01_mod_mod} can be rewritten  to yield
\begin{equation}\label{eq:app:smph_alternative}
\begin{array}{ll}
 s \geq\pm s^\mp_h
& =
\pm(\source^O,\uh)
\pm\langle \trac^O,\uh\rangle _{\GN}
\pm (\source,\uhA)
\mp\langle \trac,\uhA\rangle _{\GN}
\mp(\nu\grad  \uh,\grad  \uhA)
\\[1ex]
&
-\dfrac{1}{4\kappa}{\norm{(\qhA +\nu\grad   \uhA)\mp\kappa(\qh+\nu\grad  \uh)}}^2,
\end{array}
\end{equation}
which in particular holds for $\qh = \q =-\nu \grad  u$ and $\qhA = \qA=-\nu \grad  \xi$ yielding
\begin{equation}\label{eq:property_QoI_01_exact_fluxes}
\begin{array}{ll}
\pm s \geq 
&
\pm(\source^O,\uh)
\pm\langle \trac^O,\uh\rangle _{\GN}
\pm (\source,\uhA)
\mp\langle \trac,\uhA\rangle _{\GN}
\mp(\nu\grad  \uh,\grad  \uhA)
\\[1ex]
&
-\dfrac{1}{4\kappa}{\norm{\nu\grad(   \xi - \uhA\mp\kappa( u-  \uh))}}^2.
\end{array}
\end{equation}
Therefore, to compute bounds for the quantity of interest it is sufficient to be able to compute upper bounds for
\[
\dfrac{1}{4\kappa}{\norm{\nu\grad(   \xi - \uhA\mp\kappa( u-  \uh))}}^2 
=\dfrac{1}{4\kappa}{\norm{\nu\grad( \uPA^\pm - \uhPA^\pm )}}^2
\]
where $\uhPA^\pm= \pm \uhA -\kappa\uh$ and $\uPA^\pm= \pm \xi -\kappa u$
satisfies
\begin{equation}\label{eq:strongPA}
\begin{array}{rcll}
    -\nabla\cdot(\nu \grad \uPA^\pm)       &\!\!=\!\!& 
        \pm\source^O-\kappa\source =\source^\pm    & \text{ in } \Omega, \\
    \uPA^\pm                           &\!\!=\!\!& \pm\uD^O-\kappa\uD         & \text{ on } \GD, \\
    -\nu \grad  \uPA^\pm \cdot \bn  &\!\!=\!\!&\mp\trac^O-\kappa \trac  = \trac^\pm       & \text{ on }
    \GN.
\end{array}
\end{equation}

Upper bounds for the energy norm ${\norm{\nu\grad( \uPA^\pm - \uhPA^\pm )}}^2$ are computed introducing the zero-order equilibrated flux reconstruction of $\uPA^\pm$, namely  $(\qhPA^\mp)^0 =\pm \qhAZ -\kappa \qhZ\in\testq$ such that
\begin{equation}\label{eq:uPApm_ZO}
\begin{array}{l}
(\nabla\cdot(\qhPA^\mp)^0,1)=(\Pi_\elem^0(\pm f^O-\kappa  f),1)=(f^\pm,1) \text{ in } \Omega,
\\[1ex]
((\qhPA^\mp)^0\cdot\bn,1)=(\Pi_\edge^0(\mp \trac^O-\kappa \trac),1)=(\trac^\pm,1) \text{ on } \GN.
\end{array}
\end{equation}

Indeed, let $w\in\test$ be such that $\at{w}{\GD}=0$, that is, $w\in H^1_0(\Omega)$. Using equation \eqref{eq:green} for $\omega=\Omega$ and $\q = (\qhPA^\mp)^0$, namely
\[
(\nabla\cdot(\qhPA^\mp)^0, w)-\langle (\qhPA^\mp)^0\cdot\bn ,w\rangle _{\GN}+((\qhPA^\mp)^0,\grad  w)  = 0
\]
and equation \eqref{eq:strongPA} yields after some rearrangements
\begin{equation}\label{eq:appZO:aux01}
\begin{array}{l}
(\nu\grad( \uPA^\pm - \uhPA^\pm ),\grad w)
=
-\langle \trac^\pm-(\qhPA^\mp)^0\cdot\bn  , w\rangle _{\GN}
+(\source^\pm-\nabla\cdot(\qhPA^\mp)^0,w)
-((\qhPA^\mp)^0+\nu\grad \uhPA^\pm,\grad w)
\\[2ex]
\qquad
=\sum\limits_{\elem\in\triang} \Bigl[
(\source^\pm-\nabla\cdot(\qhPA^\mp)^0,w)_\elem 
-((\qhPA^\mp)^0+\nu\grad \uhPA^\pm,\grad w)_\elem 	
- \sum\limits_{\edge\in \GN \cap \partial \elem }\langle \trac^\pm-(\qhPA^\mp)^0\cdot\bn  , w\rangle _{\GN}		
	\Bigr].
\end{array}
\end{equation}

In order to bound the three terms in the previous summation, we need to introduce the following Poincaré and trace inequalities 
\[
\begin{array}{c}
||w-\Pi_\elem^0 w||_{\mathcal{L}^2(\elem )} \leq C_1 ||\grad w||_{\mathcal{L}^2(\elem )}
=C_1 \nu_\elem^{-1/2} {\norm{\nu_\elem  \grad w}}_\elem  
\\[1ex]
||w-\Pi_\edge^0 w||_{\mathcal{L}^2(\edge)}
\leq C_2||\grad w||_{\mathcal{L}^2(\elem )} 
=C_2 \nu_\elem^{-1/2} {\norm{\nu_\elem  \grad w}}_\elem,  
\end{array}
\] 
where{, recall that, } $||\cdot||_{\mathcal{L}^2(\elem)}$ denotes the $\mathcal{L}^2(\elem )$ norm both in $\mathbb{R}$ and $\mathbb{R}^d$, ${\norm{\cdot}}_\elem $ is the restriction of the energy norm defined in \eqref{eq:a_wv_wv} to element $\elem$ and 
\begin{equation}\label{eq:C1C2}
C_1=h_\elem /\pi \qquad,\qquad
C_2^2 = \dfrac{|\edge|}{d|\elem|}\dfrac{h_\elem }{\pi}
\left(2\max\limits_{\bm{x}\in\edge} |\bm{x} - \bm{x}_\edge|+\dfrac{d \,h_\elem }{\pi}\right),
\end{equation}
where $\bm{x}_\edge$ denotes the vertex of element $\elem$ opposite to the facet $\edge$, $|\bm{x} - \bm{x}_\edge|$ denotes the $\mathbb{R}^d$ Euclidean norm of the vector $\bm{x} - \bm{x}_\edge$, 
$|\edge|$ is the measure of the facet $\gamma$ and $h_k=\max_{\bm{x},\bm{y}\in\elem} |\bm{x} - \bm{y}|$ and $|\elem|$ are the diameter
and measure of element $\elem$ respectively. Note that 
$\max_{\bm{x}\in\edge} |\bm{x} - \bm{x}_\edge|$ can be replaced by $h_\elem $ and the inequalities still hold.
The proof of these results can be found in \cite{MR3335498,MR3577961,MR4007990,MR3262938}.

Then, since $(\qhPA^\mp)^0$ satisfies \eqref{eq:uPApm_ZO} it holds that
\begin{equation}\label{eq:lemma2_boundf}
\begin{array}{rl}
\displaystyle
\int_{\elem} (\source^\pm-\nabla\cdot(\qhPA^\mp)^0) w \ d\Omega
	&\displaystyle 
	=\int_{\elem} (\source^\pm-\nabla\cdot(\qhPA^\mp)^0) (w-\Pi_\elem^0 w) \ d\Omega
	\\[2ex]&\displaystyle
	\leq ||\source^\pm-\nabla\cdot(\qhPA^\mp)^0||_{\mathcal{L}^2(\elem )}
	||w-\Pi_\elem^0 w||_{\mathcal{L}^2(\elem )}
	\\[2ex]
	&\displaystyle 
	\leq  
	C_1 \nu_\elem^{-1/2}
	||\source^\pm-\nabla\cdot(\qhPA^\mp)^0||_{\mathcal{L}^2(\elem )}
	{\norm{\nu_\elem  \grad w}}_\elem 
\end{array}	
\end{equation}
and
\begin{equation}\label{eq:lemma2_boundg}
\begin{array}{rl}
\displaystyle
\int_{\edge} (\trac^\pm-(\qhPA^\mp)^0\cdot\bn) w \ d\Gamma 
&\displaystyle
 = \int_{\edge} (\trac^\pm-(\qhPA^\mp)^0\cdot\bn) (w-\Pi_\edge^0 w) \ d\Gamma 
\\[2ex]&\displaystyle
\leq 
	||\trac^\pm-(\qhPA^\mp)^0\cdot\bn||_{\mathcal{L}^2(\edge)}
	||w-\Pi_\edge^0 w||_{\mathcal{L}^2(\edge)}
\\[2ex]
&\displaystyle	
\leq
C_2 \nu_\elem^{-1/2} ||\trac^\pm-(\qhPA^\mp)^0\cdot\bn||_{\mathcal{L}^2(\edge)}
{\norm{\nu_\elem  \grad w}}_\elem . 
\end{array}
\end{equation}
Finally, it also holds that
\[
(\nu\grad \uhPA^\pm+(\qhPA^\mp)^0 ,\grad w)_\elem 
=(\nu^{-1}((\qhPA^\mp)^0+\nu\grad \uhPA^\pm) ,\nu \grad w)_\elem 
\leq 
{\norm{(\qhPA^\mp)^0+\nu\grad \uhPA^\pm}}_\elem {\norm{\nu\grad w}}_\elem,
\]
which introduced in \eqref{eq:appZO:aux01} along with the previous inequalities yields
\[
\hspace{-1.8cm}
\begin{array}{rl}
(\nu\grad( \uPA^\pm - \uhPA^\pm ),\grad w)
	&\displaystyle						
	\leq 
	\sum\limits_{\elem\in\triang} \Bigl[
	{\norm{(\qhPA^\mp)^0+\nu\grad \uhPA^\pm}}_\elem  
	+
	C_1 \nu_\elem^{-1/2}
	||\source^\pm-\nabla\cdot(\qhPA^\mp)^0||_{\mathcal{L}^2(\elem)}	
	\\[2ex]&\displaystyle
	\qquad\qquad+ \sum\limits_{\edge\in \GN \cap \partial\elem}
	\CTb \nu_\elem^{-1/2}||\trac^\pm-(\qhPA^\mp)^0\cdot\bn||_{\mathcal{L}^2(\edge)}
	\Bigr]{\norm{\nu \grad w}}_\elem 
%
%
	\\[2ex]&\displaystyle
	= 
	\sum\limits_{\elem\in\triang} \eta_K^{0\mp} {\norm{\nu \grad w}}_\elem 
	\leq 
	\sqrt{\sum\limits_{\elem\in\triang}(\eta_K^{0\mp})^2}
	\sqrt{\sum\limits_{\elem\in\triang}{\norm{\nu \grad w}}_\elem^2}
	=
	\sqrt{\sum\limits_{\elem\in\triang}(\eta_K^{0\mp})^2}
	{\norm{\nu \grad w}}.
\end{array}
\]

Finally, since $\uPA^\pm - \uhPA^\pm|_{\GD}=0$, we can substitute $w=\uPA^\pm - \uhPA^\pm$ in the previous inequality to yield
\[
{\norm{\nu\grad( \uPA^\pm - \uhPA^\pm )}}^2 
 = (\nu\grad( \uPA^\pm - \uhPA^\pm ),\grad( \uPA^\pm - \uhPA^\pm )) 
	\leq 
	\sqrt{\sum\limits_{\elem\in\triang} (\eta_K^{0\mp})^2}
	{\norm{\nu \grad ( \uPA^\pm - \uhPA^\pm )}}
\]
and therefore
\[
{\norm{\nu\grad( \uPA^\pm - \uhPA^\pm )}}^2 
\leq \sum\limits_{\elem\in\triang} (\eta_K^{0\mp})^2	
\]
yielding the desired bound
\begin{equation}\label{eq:app:property_QoI_01_exact_fluxes_DO}
\pm s \geq
\pm(\source^O,\uh)
\pm\langle \trac^O,\uh\rangle _{\GN}
\pm (\source,\uhA)
\mp\langle \trac,\uhA\rangle _{\GN}
\mp(\nu\grad  \uh,\grad  \uhA)
-\dfrac{1}{4\kappa}\sum\limits_{\elem\in\triang} (\eta_K^{0\mp})^2.
\end{equation}

Finally, the estimator $\eta_K^{0\mp}$ can be rewritten explicitly in terms of the primal and adjoint problems as
\[
\begin{array}{rl}
\eta_K^{0\mp}
&=
{\norm{\pm (\qhAZ +\nu\grad \uhA) -\kappa (\qhZ+\nu\grad \uh)}}_\elem  
	+
	C_1 \nu_\elem^{-1/2}
	||\pm(\source^O-\nabla\cdot \qhAZ )-\kappa(\source-\nabla\cdot \qhZ)||_{\mathcal{L}^2(\elem)}	
\\[1ex]
&\displaystyle	
	+ \sum\limits_{\edge\in \GN \cap \partial\elem}
	\CTb \nu_\elem^{-1/2}||\mp(\trac^O+\qhAZ\cdot\bn)-\kappa (\trac-\qhZ\cdot\bn)||_{\mathcal{L}^2(\edge)}
\end{array}
\]

\section{Exact representation for the quantity of interest -- Proof of Theorem \ref{th:property_QoI_03}}
\label{App:eq:property_QoI_02}

The bounds given by \eqref{eq:app01:aux03} are exact  if  
$(\uhPA^\mp, \qhPA^\mp)
=(\uPA^\mp, \qPA^\mp) = (\mp \xi -\kappa u,\pm \qA -\kappa \q)$ since the infimum is reached imposing $(w^\mp,\bv^\mp)$ in \eqref{eq:reconstructionsComb} to be $(w^\mp,\bv^\mp)=(u,\q)$  for all values of $\kappa$. Moreover, in this case, from equation \eqref{eq:app01:aux05} if holds that
\begin{equation}\label{eq:app_C1}
\pm s =
	\Lag(w^\mp,\bv^\mp;\uPA^\mp, \qPA^\mp)
	=-\dfrac{1}{4\kappa}{\norm{\qPA^\mp-\nu\grad  \uPA^\mp}}^2 - \ell(\uPA^\mp,\qPA^\mp).
\end{equation}
It is worth noting that this exact representation can also be algebraically derived by substituting $(\uPA^\mp, \qPA^\mp) = (\mp \xi -\kappa u,\pm \qA -\kappa \q)$ into the right-hand side of \eqref{eq:app_C1} and simplifying the terms appearing therein.

Let now $(\uh,\qh)$ and $(\uhA,\qhA)$ be two pair of approximations both in $\test\times\testq$ but not necessarily satisfying \eqref{eq:reconstructions}, and define the errors in the approximations as
\[
\begin{array}{l}
u = \uh + \eu \, , \,
\xi = \uhA + \euA \, , \,
\uPA^\mp=\uhPA^\mp + e_\phi^\mp\, , \,
\\[1ex]
\q = \qh + \eq \, , \,
\qA = \qhA + \eqA \, , \,
\qPA^\mp = \qhPA^\mp + \bm{e}_{\bm{\varphi}}^\mp.
\end{array}
\]
Then, it holds that 
\[
\begin{array}{l}
\dfrac{1}{4\kappa}{\norm{\qhPA^\mp-\nu\grad  \uhPA^\mp}}^2 
+ \ell(\uhPA^\mp,\qhPA^\mp)
\\[2ex]
\qquad=
\dfrac{1}{4\kappa}{\norm{\qPA^\mp-\nu\grad  \uPA^\mp - \bm{e}_{\bm{\varphi}}^\mp+\nu\grad  e_\phi^\mp}}^2 
+\ell(\uPA^\mp,\qPA^\mp)
-\ell(e_\phi^\mp,\bm{e}_{\bm{\varphi}}^\mp)
\\[2ex]
\qquad=
\mp s
+\dfrac{1}{4\kappa}{\norm{\bm{e}_{\bm{\varphi}}^\mp-\nu\grad  e_\phi^\mp}}^2 
+(u,\nabla\cdot\bm{e}_{\bm{\varphi}}^\mp)
-\langle u,\bm{e}_{\bm{\varphi}}^\mp\cdot\bn\rangle _{\GN}
-\langle  \q\cdot\bn, e_\phi^\mp\rangle _{\GD}
\\[2ex]
\qquad=
\mp s
+\dfrac{1}{4\kappa}{\norm{\bm{e}_{\bm{\varphi}}^\mp-\nu\grad  e_\phi^\mp}}^2 
+(u,\pm (\source^O - \nabla\cdot\qhA ) -\kappa(\source - \nabla\cdot\qh))
\\[1ex]
\qquad
-\langle u,\mp (\trac^O + \qhA  \cdot\bn) -\kappa(\trac - \qh \cdot\bn)\cdot\bn\rangle _{\GN}
-\langle  \q\cdot\bn, \mp(\uD^O - \uhA) - \kappa (\uD - \uh)\rangle _{\GD},
\end{array}	
\]
where we have used that $\qPA^\mp-\nu\grad \uPA^\mp =-2\kappa\q$, equation \eqref{eq:general_variational_formulation} with $(w,\bv)=(e_\phi^\mp,\bm{e}_{\bm{\varphi}}^\mp)$ and the fact that
\[
\begin{array}{ll}
e_\phi^\mp 
=\mp(\uD^O - \uhA) - \kappa (\uD - \uh) &\text{ on } \GD
\\[1ex]
\nabla\cdot\bm{e}_{\bm{\varphi}}^\mp 
= \pm (\source^O - \nabla\cdot\qhA ) -\kappa(\source - \nabla\cdot\qh)&\text{ in } \Omega
\\[1ex]
\bm{e}_{\bm{\varphi}}^\mp \cdot\bn
= \mp (\trac^O + \qhA  \cdot\bn) -\kappa(\trac - \qh \cdot\bn)
&\text{ on } \GN.
\end{array}	
\]

Finally, the Theorem is proved by noting that $\mp \hat{s}_h^\mp$ coincides with
\[
\begin{array}{l}
\dfrac{1}{4\kappa}{\norm{\qhPA^\mp-\nu\grad  \uhPA^\mp}}^2 
+ \ell(\uhPA^\mp,\qhPA^\mp)
= \dfrac{1}{4\kappa}{\norm{\pm (\qhA +\nu\grad   \uhA )-\kappa (\qh-\nu\grad  \uh)}}^2
 \mp\ell(\uhA,-\qhA )
-\kappa\ell(\uh,\qh).
\end{array}
\]

\section{Lower bounds for the energy norm of  ${\norm{\bm{e}_{\bm{\varphi}}^\mp-\nu\grad  e_\phi^\mp}}^2$ -- Proof of equation \eqref{eq:exactQoI_nodataoscillations_LB}}
\label{App:LB}

Let $(\uh,\qh)$ and $(\uhA,\qhA)$ be potential and equilibrated flux reconstructions of the primal and adjoint problems satisfying \eqref{eq:reconstructions}, and consider 
$e_\phi^\mp 
=\mp(\xi - \uhA) - \kappa (u - \uh)$ and 
$
\bm{e}_{\bm{\varphi}}^\mp 
= \pm (\qA - \qhA ) -\kappa(\q - \qh)
$. Then, in two and three dimensions, a lower bound for ${\norm{\bm{e}_{\bm{\varphi}}^\mp-\nu\grad  e_\phi^\mp}}^2$ can be computed using a Helmholtz decomposition of $\bm{e}_{\bm{\varphi}}^\mp$, see 
\cite{MR3850360,MR851383}. Indeed, since $e_\phi^\mp\in H^1_0(\Omega)$, 
 the error $\bm{e}_{\bm{\varphi}}^\mp-\nu\grad  e_\phi^\mp\in\testq\subset[\mathcal{L}^2(\Omega)]^d$ can be rewritten in the form $\btau -\nu\grad  e_\phi^\mp = \nu \grad (\btauZ-e_\phi^\mp) + \curl \btauC$ 
where $\btauZ\in\mathcal{H}^1_0(\Omega)$ satisfies
\[
(\nu\grad\btauZ,\grad w)= (\btau,\grad w) \quad\forall w\in H^1_0(\Omega)
\]
and $\btauC\in \mathcal{H}^1_{0\times}(\Omega) =\{\psi \in [H^1(\Omega)]^{2d-3} , \curl\psi\cdot\bn = 0 \text{ on } \GN\}$ satisfies
\[
(\nu^{-1}\curl\btauC,\curl\btauC)=(\nu^{-1}\btau,\curl\btauC),
\]
where $\grad\times$ is the standard curl operator, see \cite[Sec 2.3]{MR851383}.
Now, for any $w^\mp\in H^1_0(\Omega)$ and $\psi^\mp\in [H^1(\Omega)]^{2d-3}$, consider
\[
\ell_\times^\mp(w^\mp,\psi^\mp)=
(\nu^{-1}(\bm{e}_{\bm{\varphi}}^\mp-\nu\grad  e_\phi^\mp), \nu\grad w^\mp + \curl \psi^\mp),
\]
and the associated scalar parameter $\lambda^\mp=-\ell_\times^\mp(w^\mp,\psi^\mp)/{\norm{\nu\grad w^\mp + \curl \psi^\mp}}^2$. Then,
\[
\hspace{-0cm}
\begin{array}{l}
{\norm{\bm{e}_{\bm{\varphi}}^\mp-\nu\grad  e_\phi^\mp + \lambda^\mp(\nu\grad w^\mp + \curl \psi^\mp)}}^2 
\\[2ex]
\qquad
=
{\norm{\bm{e}_{\bm{\varphi}}^\mp-\nu\grad  e_\phi^\mp }}^2 
+(\lambda^\mp)^2{\norm{\nu\grad w^\mp + \curl \psi^\mp}}^2 
+2\lambda^\mp(\nu^{-1}(\bm{e}_{\bm{\varphi}}^\mp-\nu\grad  e_\phi^\mp), \nu\grad w^\mp + \curl \psi^\mp)
\\[2ex]
\qquad
=
{\norm{\bm{e}_{\bm{\varphi}}^\mp-\nu\grad  e_\phi^\mp }}^2 
+(\lambda^\mp)^2{\norm{\nu\grad w^\mp + \curl \psi^\mp}}^2 
+2\lambda^\mp\ell_\times^\mp(w^\mp,\psi^\mp)
\\[2ex]
\qquad
=
{\norm{\bm{e}_{\bm{\varphi}}^\mp-\nu\grad  e_\phi^\mp }}^2 
-\dfrac{(\ell_\times^\mp(w^\mp,\psi^\mp))^2}{{\norm{\nu\grad w^\mp + \curl \psi^\mp}}^2}.
\end{array}
\]
which yields to
\[
{\norm{\bm{e}_{\bm{\varphi}}^\mp-\nu\grad  e_\phi^\mp }}^2 =
{\norm{\bm{e}_{\bm{\varphi}}^\mp-\nu\grad  e_\phi^\mp + \lambda^\mp(\nu\grad w^\mp + \curl \psi^\mp)}}^2 
+\dfrac{(\ell_\times^\mp(w^\mp,\psi^\mp))^2}{{\norm{\nu\grad w^\mp + \curl \psi^\mp}}^2}
\geq 
\dfrac{(\ell_\times^\mp(w^\mp,\psi^\mp))^2}{{\norm{\nu\grad w^\mp + \curl \psi^\mp}}^2}.
\]
Moreover, for $w^\mp=\btauZ-e_\phi^\mp$ and $\psi^\mp = \btauC$, 
since $\nu\grad w^\mp + \curl \psi^\mp=\btau -\nu\grad  e_\phi^\mp 
$ then $\ell_\times^\mp(w^\mp,\psi^\mp)
={\norm{\btau -\nu\grad  e_\phi^\mp}}^2
$
and 
the previous inequality becomes an equality yielding to
\[
{\norm{\bm{e}_{\bm{\varphi}}^\mp-\nu\grad  e_\phi^\mp }}^2 =
\sup\limits_{\scriptsize\begin{array}{c}w^\mp\in H^1_0(\Omega)\\\psi^\mp\in [H^1(\Omega)]^{2d-3}\end{array}}
\dfrac{(\ell_\times^\mp(w^\mp,\psi^\mp))^2}{{\norm{\nu\grad w^\mp + \curl \psi^\mp}}^2}.
\]

Equation \eqref{eq:exactQoI_nodataoscillations_LB} is finally proved by noting that  if $\uh$ and $\uhA$ are potential reconstructions of the primal and adjoint problems, since
\[
(\nabla v ,\curl \psi) = 0 \qquad \forall v \in H_0^1 , \psi \in H^1_{0\times}.
\]
then
\[
\begin{array}{l}
\ell_\times^\mp(w^\mp,\psi^\mp)=
\mp(\nu^{-1} (\qhA +\nu\grad\uhA ),\nu\grad w^\mp + \curl \psi^\mp )
\\[1ex]
\qquad
+ \kappa \left( 2 (\source,w^\mp) - 2<\trac,w^\mp>_{\GN}
+(\qh -\nu\grad \uh,\grad w^\mp )
+(\nu^{-1}(\qh+\nu\grad \uh ),\curl \psi^\mp)\right).
\end{array}
\] 
Moreover, if $\qh$ and $\qhA$ are equilibrated flux reconstructions,
$\ell_\times^\mp(w^\mp,\psi^\mp)$ reduces to
\[
\ell_\times^\mp(w^\mp,\psi^\mp)=
\mp(\nu^{-1}(\qhA +\nu\grad  \uhA),\nu\grad w^\mp+\curl \psi^\mp)
- \kappa( \nu^{-1}(\qh+\nu\grad \uh),\nu\grad w^\mp-\curl \psi^\mp),
\]
proving the desired result.

\end{document}